\documentclass[11pt]{sat}

\title{Korovkin-type Theorems and Approximation by Positive Linear Operators}
\def\shorttitle{Korovkin-type theorems and positive operators}

\author{Francesco Altomare}
\def\shortauthor{F.~Altomare}

\def\versiondate{10 September 2010}

\def\abstracttext{
This survey paper contains a detailed self-contained introduction to
Korovkin-type theorems and to some of their applications concerning
the approximation of continuous functions as well as of
$L^p$-functions, by means of positive linear operators.

The paper  also contains several new results and applications.
Moreover, the organization of the subject follows a simple and direct
approach which quickly leads both to the main results of the theory
and to some new ones. }

\def\MSCnumbers{41A36, 46E05, 47B65} 

\def\keywords{ Korovkin-type theorem, positive
operator, approximation by positive operators, Stone-Weierstrass
theorem, (weighted) continuous function space, $L^p$-space.
} 

\input colordvi

\def\Appendix{Appendix}
\def\norm#1{\Vert#1\Vert}
\def\Norm#1{\left\Vert#1\right\Vert}
\def\dword#1{{\bf #1}}
\let\eword\emph
\def\dd{\,{\rm d}}  
\def\ee{{\rm e}}  
\def\ii{{\rm i}}  
\def\supp{\mathop{\rm supp}\nolimits}

\usepackage{amsfonts}
\usepackage{amsmath}
\usepackage{amsthm}
\usepackage{amssymb}
\usepackage{graphicx}%

\newtheorem{theorem}{Theorem}[section]
\newtheorem{proposition}[theorem]{Proposition}
\newtheorem{lemma}[theorem]{Lemma}
\newtheorem{corollary}[theorem]{Corollary}
\newtheorem{definition}[theorem]{Definition}
\newtheorem{examples}[theorem]{Examples}
\newtheorem{example}[theorem]{Example}
\theoremstyle{definition}
\newtheorem{remark}[theorem]{Remark}
\newtheorem{remarks}[theorem]{Remarks}
\newtheorem{problem}[theorem]{Problem}
\numberwithin{equation}{section} \numberwithin{figure}{section}
\numberwithin{table}{section} \numberwithin{equation}{section}
\renewenvironment{proof}[1][{\textbf Proof}]{\noindent\textbf{#1.} }{\hfill$\Box$\vskip8pt}

\def\startpagenumber{92}
\def\volumenumber{6}
\def\year{2010}

\newcommand{\beginddoc}{
\begin{document}
\maketitle \vskip-\baselineskip \centerline{\footnotesize\emph{To
the memory of my parents}} \centerline{\footnotesize\emph{Maria
Giordano (1915-1989) and Luigi Altomare (1898-1963)}}
\vskip\baselineskip\vskip\baselineskip
\begin{abstract}
\abstracttext \vskip1pt MSC: \MSCnumbers
\ifx\keywords\empty\else\vskip1pt Keywords: \keywords\fi
\end{abstract}
\insert\footins{\scriptsize
\medskip
\baselineskip 8pt \leftline{Surveys in Approximation Theory}
\leftline{Volume \volumenumber, \year.
pp.~\thepage--\pageref{endpage}.} \leftline{\copyright\ \year\
Surveys in Approximation Theory.} \leftline{ISSN 1555-578X}
\leftline{All rights of reproduction in any form reserved.}
\smallskip
\par\allowbreak}
\tableofcontents}
\renewcommand\rightmark{\ifodd\thepage{\it \hfill\shorttitle\hfill}\else {\it \hfill\shortauthor\hfill}\fi}
\markboth{{\it \shortauthor}}{{\it \shorttitle}} \markright{{\it
\shorttitle}}
\def\endddoc{\label{endpage}\end{document}}
\date{{\small \versiondate}}
\setlength\oddsidemargin{0pc} \setlength\evensidemargin{0pc}
\setlength\topmargin{0in} \setlength\textwidth{6.5in}
\setlength\textheight{8.6in}
\beginddoc



\section{Introduction}
Korovkin-type theorems furnish simple and useful tools for
ascertaining whether a given sequence of positive linear operators,
acting on some function space is an approximation process or,
equivalently, converges strongly to the identity operator.

Roughly speaking, these theorems exhibit a variety of test subsets
of functions which guarantee that the approximation (or the
convergence) property holds on the whole space provided it holds on
them.

The custom of calling these kinds of results ``Korovkin-type
theorems" refers to P.~P.~Korovkin who in 1953 discovered such a
property for the functions $\textbf{1}, x$ and $x^2$ in the space
$C([0,1])$ of all continuous functions on the real interval $[0,1]$
as well as for the functions $\textbf{1}, \cos$ and $\sin$ in the
space of all continuous $2\pi$-periodic functions on the real line
([77-78]).

After this discovery, several mathematicians have undertaken the
program of extending Korovkin's theorems in many ways and to several
settings, including function spaces, abstract Banach lattices,
Banach algebras, Banach spaces and so on. Such developments
delineated a theory which is nowadays referred to as Korovkin-type
approximation theory.

This theory has fruitful connections with real analysis, functional
analysis, harmonic analysis, measure theory and probability theory,
summability theory and partial differential equations. But the
foremost applications are concerned with constructive approximation
theory which uses it as a valuable tool.

Even today, the development of Korovkin-type approximation theory is
far from complete, especially for those parts of it that concern
limit operators different from the identity operator (see Problems
5.3 and 5.4 and the subsequent remarks).

A quite comprehensive picture of what has been achieved in this
field until 1994 is documented in the monographs of Altomare and
Campiti ([8], see in particular Appendix D), Donner ([46]), Keimel
and Roth ([76]), Lorentz, v. Golitschek and Makovoz ([83]). More
recent results can be found, e.g., in [1], [9-15], [22], [47-52],
[63], [71-74], [79], [114-116], [117] and the references therein.

The main aim of this survey paper is to give a detailed
self-contained introduction to the field as well as a secure entry
into a theory that provides useful tools for understanding and
unifying several aspects pertaining, among others, to real and
functional analysis and which leads to several applications in
constructive approximation theory and numerical analysis.

This paper, however, not only presents a survey on Korovkin-type
theorems but also contains several new results and applications.
Moreover, the organization of the subject follows a simple and
direct approach which quickly leads  both to the main results of the
theory and to some new ones.

In Sections 3 and 4, we discuss the first and the second theorem of
Korovkin.  We obtain both of them from a simple unifying result
which we state in the setting of metric spaces (see Theorem 3.2).

This general result also implies the multidimensional extension of
Korovkin's theorem due to Volkov ([118]) (see Theorem 4.1).
Moreover, a slight extension of it into the framework of locally
compact metric spaces allows to extend the Korovkin's theorems to
arbitrary real intervals or, more generally, to locally compact
subsets of $\mathbb{R}^d$, $d\geq 1$.

Throughout the two sections, we present some applications concerning
several classical approximation processes ranging from Bernstein
operators on the unit interval or on the canonical hypercube and the
multidimensional simplex, to Kantorovich operators, from Fej\'{e}r
operators to Abel-Poisson operators, from Sz\'{a}sz-Mirakjan
operators to Gauss-Weierstrass operators.

We also prove that the first and the second theorems of Korovkin are
actually equivalent to the algebraic and the trigonometric version,
respectively, of the classical Weierstrass approximation theorem.

Starting from Section 5, we enter into the heart of the theory by
developing some of the main results in the framework of the space
$C_0(X)$ of all real-valued continuous functions vanishing at
infinity on a locally compact space $X$ and, in particular, in the
space $C(X)$ of all real-valued continuous functions on a compact
space $X$.

We choose these continuous function spaces because they play a
central role in the whole theory and are the most useful for
applications. Moreover, by means of them it is also possible to
easily obtain some Korovkin-type theorems in weighted continuous
function spaces and in $L^p$-spaces, $1\leq p$.  These last aspects
are treated at the end of Section 6 and in Section 8.

We point out that we discuss Korovkin-type theorems not only with
respect to the identity operator but also with respect to a positive
linear operator on $C_0(X)$ opening the door to a variety of
problems some of which are still unsolved.

In particular, in Section 10, we present some results concerning
positive projections on $C(X), X$ compact, as well as their
applications to the approximation of the solutions of Dirichlet
problems  and of other similar problems.

In Sections 6 and 7, we present several results and applications
concerning Korovkin sets for the identity operator.  In particular,
we show that, if $M$ is a subset of $C_0(X)$ that separates the
points of $X$ and if $f_0\in C_0(X)$ is strictly positive, then
$\{f_0\}\cup f_0M\cup f_0M^2$ is a Korovkin set in $C_0(X)$.

This result is very useful because it furnishes a simple way to
construct Korovkin sets, but in addition, as we show in Section 9,
it turns out that it is equivalent to the Stone generalization to
$C_0(X)$-spaces of the Weierstrass theorem. This equivalence was
already established in [8, Section 4.4] (see also [12-13]) but here
we furnish a different, direct and more transparent proof.

We also mention that, at the end of Sections 7 and 10, we present
some applications concerning Bernstein-Schnabl operators associated
with a positive linear operator and, in particular, with a positive
projection. These operators are useful for the approximation of not
just continuous functions but also --- and this was the real reason
for the increasing interest in them --- positive semigroups and
hence the solutions of initial-boundary value evolution problems.
These aspects are briefly sketched at the end of Section 10.

Following the main aim of ``Surveys in Approximation Theory", this
paper is directed to the graduate student level and beyond. However,
some parts of it as well as some new methods developed here could
also be useful to expert readers.

A knowledge of the basic definitions and results concerning locally
compact Hausdorff spaces and continuous function spaces on them is
required as well as some basic properties of positive linear
functionals on these function spaces (Radon measures). However, the
reader who is not interested in this level of generality may replace
everywhere our locally compact spaces with the space $\mathbb{R}^d,
d\geq 1$, or with an open or a closed subset of it  or with the
intersection of an open subset and a closed subset of
$\mathbb{R}^d$. However, this restriction does not produce any
simplification of the proofs or of the methods.

For the convenience of the reader and to make the exposition
self-contained, we collect all these prerequisites in the \Appendix.
There, the reader can also find some new simple and direct proofs of
the main properties of Radon measures which are required throughout
the paper, so that no a priori knowledge of the theory of Radon
measures is needed.

This paper contains introductory materials so that many aspects of
the theory have been omitted.  We refer, e.g., to [8, Appendix D]
for further details about some of the main directions developed
during the last fifty years.

Furthermore, in the applications shown throughout the paper, we
treat only general constructive aspects (convergence of the
approximation processes) without any mention of quantitative aspects
(estimates of the rate of convergence, direct and converse results
and so on) nor to shape preserving properties. For such matters, we
refer, e.g., to [8], [26], [38], [41], [42], [44], [45], [64], [65],
[81-82], [83], [92], [109].

We also refer to [19, Proposition 3.7], [67], [69] and [84] where
other kinds of convergence results for sequences of positive linear
operators can be found. The results of these last papers do not
properly fall into the Korovkin-type approximation theory but they
can be fruitfully used to decide whether a given sequence of
positive linear operators is strongly convergent (not necessarily to
the identity operator).

Finally we wish to express our gratitude to Mirella Cappelletti
Montano, Vita Leonessa and Ioan Ra\c{s}a for the careful reading of
the manuscript and for many fruitful suggestions. We are also
indebted to Carl de Boor, Allan Pinkus and Vilmos Totik for their
interest in this work as well as for their valuable advice and for
correcting several inaccuracies. Finally we want to thank Mrs.
Voichita Baraian for her precious collaboration in preparing the
manuscript in LaTeX for final processing.

\section{Preliminaries and notation} In this section, we assemble the main
notation which will be used throughout the paper together with some
generalities.

Given a metric space $(X,d)$, for every $x_0\in X$ and $r>0$, we
denote by $B(x_0,r)$ and $B'(x_0,r)$ the open ball and the closed
ball with center $x_0$ and radius $r$, respectively, i.e.,
\begin{equation}\label{2.1}
B(x_0,r):=\{x\in X \mid d(x_0,x)<r\}
\end{equation} and
\begin{equation}\label{2.2}
B'(x_0,r):=\{x\in X \mid d(x_0,x)\leq r\}.
\end{equation}
The symbol $$F(X)$$ stands for the linear space of all real-valued
functions defined on $X$. If $M$ is a subset of $F(X)$, then by
$\mathcal{L}(M)$ we designate the linear subspace generated by $M$.
\noindent We denote by
$$B(X)$$
 the linear subspace of all functions
$f:X\longrightarrow\mathbb{R}$ that are bounded, endowed with the
norm of uniform convergence (briefly, the sup-norm) defined by
\begin{equation}\label{2.3}
\norm{f}_{\infty}:=\underset{x\in X}{\sup}|f(x)|\quad (f\in B(X)),
\end{equation}
\noindent with respect to which it is a Banach space.

The symbols $$C(X)\quad\text{and}\quad C_b(X)$$ denote the linear
subspaces of all continuous (resp.\  continuous and bounded)
functions in $F(X)$. Finally, we denote by $$UC_b(X)$$ the linear
subspace of all uniformly continuous and bounded functions in
$F(X)$. Both $C_b(X)$ and $UC_b(X)$ are closed in $B(X)$ and hence,
endowed with the norm (\ref{2.3}), they are Banach spaces.

A linear subspace $E$ of $F(X)$ is said to be a \dword{lattice
subspace} if
\begin{equation}\label{2.4}
|f|\in E\quad\quad \text{for every}\quad f\in E.
\end{equation} For
instance, the spaces $B(X), C(X), C_b(X)$ and $UC_b(X)$ are lattice
subspaces.

Note that from (\ref{2.4}), it follows that $\sup(f,g), \inf(f,g)\in
E$ for every $f,g\in E$ where
\begin{equation}\label{2.5}
\sup(f,g)(x):=\sup(f(x),g(x))\quad\quad (x\in X)
\end{equation} and
\begin{equation}\label{2.6}
\inf(f,g)(x):=\inf(f(x),g(x))\quad\quad (x\in X).
\end{equation} This
follows at once by the elementary identities
\begin{equation}\label{2.7}
\sup(f,g)=\frac{f+g+|f-g|}{2} \quad \textrm{ and } \quad
\inf(f,g)=\frac{f+g-|f-g|}{2}.
\end{equation} More generally, if
$f_1,\dots,f_n\in E,\;n\geq 3$, then $\underset{1\leq i\leq
n}{\sup}f_i, \underset{1\leq i\leq n}{\inf}f_i\in E.$

We say that a linear subspace $E$ of $F(X)$ is a \dword{subalgebra}
if
\begin{equation}\label{2.8}
f\cdot g\in E\quad\quad\text{for every}\quad f,g\in E
\end{equation} or,
equivalently, if $f^2\in E$ for every $f\in E$. In this case, if
$f\in E$ and $n\geq 1$, then $f^n\in E$ and hence for every real
polynomial $Q(x):=\alpha_1x+\alpha_2x^2+\dots+\alpha_nx^n\
(x\in\mathbb{R})$ vanishing at $0$, the function
\begin{equation}\label{2.9}
Q(f):=\alpha_1f+\alpha_2f^2+\dots+\alpha_nf^n
\end{equation} belongs to $E$ as well. If $E$ contains the constant functions,
then $P(f)\in E$ for every real polynomial $P$.

Note that a subalgebra is not necessarily a lattice subspace (for
instance, $C^1([a,b])$ is such an example).  However every closed
subalgebra of $C_b(X)$ is a lattice subspace (see Lemma 9.1).

Given a linear subspace $E$ of $F(X)$, a linear functional
$\mu:E\longrightarrow\mathbb{R}$ is said to be positive if
\begin{equation}\label{2.10}
\mu(f)\geq 0\quad\quad\text{for every}\quad f\in E,\ f\geq 0.
\end{equation}
The simplest example of a positive linear functional is the
so-called evaluation functional at a point $a\in X$ defined by
\begin{equation}\label{2.11}
\delta_a(f):=f(a)\quad\quad (f\in E).
\end{equation}

If $(Y,d')$ is another metric space, we say that a linear operator
$T:E\longrightarrow F(Y)$ is positive if
\begin{equation}\label{2.12}
T(f)\geq 0\quad\quad \text{for every}\quad f\in E,\ f\geq 0.
\end{equation} Every positive linear operator $T:E\longrightarrow F(Y)$
gives rise to a family $(\mu_y)_{y\in Y}$ of positive linear
functionals on $E$ defined by
\begin{equation}\label{2.13}
\mu _{y}(f):=T(f)(y)\quad(f\in E).
\end{equation} Below, we state some
elementary properties of both positive linear functionals and
positive linear operators.

In what follows, the symbol $\mathbb{F}$ stands either for the field
$\mathbb{R}$ or for a space $F(Y)$, $Y$ being an arbitrary metric
space.

Consider a linear subspace $E$ of $F(X)$ and a positive linear
operator $T:E\longrightarrow$ $\mathbb{F}$.  Then:
\begin{itemize}
\item [(i)] For every $f,g\in E, f\leq g,$
\begin{equation}\label{2.14}
  T(f)\leq T(g)
\end{equation}
\item [(ii)] If $E$ is a lattice subspace, then
\begin{equation}\label{2.15}
|T(f)|\leq T(|f|)\quad \text{for every }f\in E.
\end{equation}
\item [(iii)] (Cauchy-Schwarz inequality) If $E$ is both a lattice subspace and a subalgebra, then
\begin{equation}\label{2.16} T(|f\cdot
g|)\leq\sqrt{T(f^2)T(g^2)}\quad\quad (f,g\in E).
\end{equation}

In particular, if $\textbf{1}\in E$, then
\begin{equation}\label{2.17}
T(|f|)^2\leq T(\textbf{1})T(f^2)\quad\quad (f\in E).
\end{equation}
\item [(iv)]If $X$ is compact, $\textbf{1}\in E$ and $F$ is either
$\mathbb{R}$ or $B(Y)$, then $T$ is continuous and
\begin{equation}\label{2.18}
\norm{T}=\norm{T(\textbf{1}) }.
\end{equation} Thus, if
$\mu:E\longrightarrow\mathbb{R}$ is a positive linear functional,
then $\mu$ is continuous and $\norm{\mu}=\mu(\textbf{1})$.
\end{itemize}

\section{Korovkin's first theorem} Korovkin's theorem provides a very
useful and simple criterion for whether a given sequence
$(L_n)_{n\geq 1}$ of positive linear operators on $C([0,1])$ is an
\dword{approximation process}, i.e., $L_n(f)\longrightarrow f$
uniformly on $[0,1]$ for every $f\in C([0,1])$.

In order to state it, we need to introduce the functions
\begin{equation}\label{3.1}
e_m(t):=t^m\quad\quad (0\leq t\leq 1)
\end{equation} $(m\geq 1)$.
\begin{theorem}\sl (Korovkin ([77])) Let $(L_n)_{n\geq 1}$ be a sequence of
positive linear operators from $C([0,1])$ into $F([0,1])$ such that
for every $g\in\{\textbf{1},e_1,e_2\}$
$$\lim_{n\rightarrow\infty}L_n(g)=g\quad\text{uniformly on}\quad [0,1].$$
Then, for every $f\in C([0,1])$,
$$\lim_{n\rightarrow\infty}L_n(f)=f\quad\text{uniformly on}\quad [0,1].$$
\end{theorem}
Below, we present a more general result from which Theorem 3.1
immediately follows.

For every $x\in [0,1]$ consider the auxiliary function
\begin{equation}\label{3.2}
d_x(t):=|t-x|\quad\quad (0\leq t\leq 1).
\end{equation} Then
$$d_x^2=e_2-2xe_1+x^2\textbf{1}$$ and hence, if $(L_n)_{n\geq 1}$ is a
sequence of positive linear operators satisfying the assumptions of
Theorem 3.1, we get
\begin{equation}\label{3.3}
\lim_{n\rightarrow\infty}L_n(d^2_x)(x)=0
\end{equation} uniformly with
respect to $x\in[0,1]$, because for $n\geq 1$
$$L_n(d_x^2)=(L_n(e_2)-x^2)+2x(L_n(e_1)-x)+x^2(L_n(\textbf{1})-\textbf{1}).$$

After these preliminaries, the reader can easily realize that
Theorem 3.1 is a particular case of the following more general
result which, together with its modification (i.e., Theorem 3.5) as
well as the further consequences presented at the beginning of
Section 4, should also be compared with the simple but different
methods of [79].

Consider a metric space $(X,d)$.  Extending (\ref{3.2}), for any
$x\in X$ we denote by $d_x\in C(X)$ the function
\begin{equation}\label{3.4}
d_x(y):=d(x,y)\quad\quad (y\in X).
\end{equation}
\begin{theorem}\sl  Let
$(X,d)$ be a metric space and consider a lattice subspace $E$ of
$F(X)$ containing the constant functions and all the functions
$d_x^2$ ($x\in X$).  Let $(L_n)_{n\geq 1}$ be a sequence of positive
linear operators from $E$ into $F(X)$ and let $Y$ be a subset of $X$
such that
\begin{itemize}
\item [(i)]$\lim\limits_{n\rightarrow\infty}L_n(\textbf{1})=\textbf{1}$
uniformly on $Y$; \item
[(ii)]$\lim\limits_{n\rightarrow\infty}L_n(d_x^2)(x)=0$ uniformly with
respect to $x\in Y$.
\end{itemize} Then for every $f\in E\cap UC_b(X)$
$$\lim_{n\rightarrow\infty}L_n(f)=f\quad\text{uniformly on}\quad Y.$$
\end{theorem}
\begin{proof}
Consider $f\in E\cap UC_b(X)$ and $\varepsilon>0$.  Since $f$ is
uniformly continuous, there exists $\delta>0$ such that
$$|f(x)-f(y)|\leq\varepsilon\quad\quad\text{for every}\quad x,y\in
X,\,d(x,y)\leq\delta.$$ On the other hand, if $d(x,y)\geq\delta$,
then
$$|f(x)-f(y)|\leq2\norm{f}_\infty\leq\frac{2\norm{f}_\infty}{\delta^2}d^2(x,y).$$
Therefore, for $x\in X$ fixed, we obtain
$$|f-f(x)|\leq\frac{2\norm{f}_\infty}{\delta^2}d_x^2+\varepsilon\textbf{1} $$
and hence, for any $n\geq 1$,
\begin{eqnarray*}
|L_n(f)(x)-f(x)L_n(\textbf{1} )(x)|\;\leq\; L_n(|f-f(x)|)(x)
\;\leq\;\frac{2\norm{f}_\infty}{\delta^2}L_n(d_x^2)(x)+\varepsilon
L_n(\textbf{1})(x).
\end{eqnarray*} We may now easily conclude that
$\lim\limits_{n\rightarrow\infty}L_n(f)=f$ uniformly on $Y$ because
of the assumptions (i) and (ii).
\end{proof}

Theorem 3.2 has a natural generalization to completely regular
spaces (for more details, we refer to [15]).  Furthermore, the above
proof can be adapted to show the next result.

\begin{theorem}\sl
Consider $(X,d)$  and $E\subset F(X)$ as in Theorem 3.2.  Consider a
sequence $(L_n)_{n\geq 1}$ of positive linear operators from $E$
into $F(X)$ and assume that for a given $x\in X$
\begin{itemize}
\item [(i)]$\lim\limits_{n\rightarrow\infty}L_n(\textbf{1})(x)=1$;
\item [(ii)]$\lim\limits_{n\rightarrow\infty}L_n(d_x^2)(x)=0$.
\end{itemize} Then, for every bounded function $f\in E$ that is continuous
at $x$, $$\lim\limits_{n\rightarrow\infty}L_n(f)(x)=f(x).$$

\end{theorem} Adapting the proof of Theorem 3.2, we can show a further result.
We first state a preliminary lemma.
\begin{lemma}\sl
Let $(X,d)$ be a locally compact metric space.  Then for every
compact subset $K$ of $X$ and for every $\varepsilon>0$, there exist
$0<\overline{\varepsilon}<\varepsilon$ and a compact subset
$K_\varepsilon$ of $X$ such that
$$B'(x,\overline{\varepsilon})\subset K_\varepsilon\quad\text{for
every}\quad x\in K.$$
\end{lemma}

\begin{proof} Given $x\in K$, there exists $0<\varepsilon(x)<\varepsilon$
such that $B'(x,\varepsilon(x))$ is compact.  Since
$K\subset\bigcup\limits_{x\in K}B(x,\varepsilon(x)/2)$, there exist
$x_1,\dots,x_p\in K$ such that $K\subset\bigcup\limits_{i=1}^p
B(x_i,\varepsilon(x_i)/2)$.  Set
$\overline{\varepsilon}:=\underset{1\leq i\leq
p}{\min}\varepsilon(x_i)<\varepsilon$ and
$K_\varepsilon:=\bigcup\limits_{i=1}^pB'(x_i,\varepsilon(x_i))$.
Now, if $x\in K$ and $y\in X$ and if
$d(x,y)\leq\overline{\varepsilon}$, then there exists an
$i\in\{1,\dots,p\}$ such that $d(x,x_i)\leq\varepsilon(x_i)/2$, and
hence $d(y,x_i)\leq d(y,x)+d(x,x_i)\leq\varepsilon(x_i).$ Therefore
$y\in K_\varepsilon$.
\end{proof}

\begin{theorem}\sl
Let $(X,d)$ be a locally compact metric space and consider a lattice
subspace $E$ of $F(X)$ containing the constant function $\textbf{1}
$ and all the functions $d_x^2$ $(x\in X)$.  Let $(L_n)_{n\geq 1}$
be a sequence of positive linear operators from $E$ into $F(X)$ and
assume that

\begin{itemize}
  \item [(i)]$\lim\limits_{n\rightarrow\infty}L_n(\textbf{1})=\textbf{1}$
  uniformly on compact subsets of $X$; \item
  [(ii)]$\lim\limits_{n\rightarrow\infty}L_n(d_x^2)(x)=0$ uniformly on
  compact subsets of $X$.
\end{itemize} Then, for every $f\in E\cap
C_b(X)$, $$\lim_{n\rightarrow\infty}L_n(f)=f\quad\text{uniformly on
compact subsets of}\quad X.$$

\end{theorem}
\begin{proof} Fix $f\in E\cap C_b(X)$ and consider a
compact subset $K$ of $X$.  Given $\varepsilon>0$, consider
$0<\overline{\varepsilon}<\varepsilon$ and a compact subset
$K_\varepsilon$ of $X$ as in Lemma 3.4.

Since $f$ is uniformly continuous on $K_\varepsilon$, there exists
$0<\delta<\overline{\varepsilon}$ such that
$$|f(x)-f(y)|\leq\varepsilon\quad\text{for every}\quad x,y\in
K_\varepsilon,\, d(x,y)\leq\delta.$$ Given $x\in K$ and $y\in X$, if
$d(x,y)\leq\delta$, then $y\in B'(x,\overline{\varepsilon})\subset
K_\varepsilon$ and hence $|f(x)-f(y)|\leq\varepsilon.$

If $d(x,y)\geq\delta$, then
$$|f(x)-f(y)|\leq\frac{2\norm{f}_\infty}{\delta^2}d^2(x,y).$$ Therefore,
once again,
$$|f-f(x)|\leq\frac{2\norm{f}_\infty}{\delta^2}d_x^2+\varepsilon
\textbf{1}$$ so that, arguing as in the final part of the proof of
Theorem 3.2, we conclude that
$$\lim_{n\rightarrow\infty}L_n(f)(x)=f(x)$$ uniformly with respect
to $x\in K$.
\end{proof}

Theorem 3.1 was obtained by P.~P.~Korovkin in 1953 ([77], see also
[78]).  However, in [35], H.~Bohman showed a result like Theorem 3.1
by considering sequences of positive linear operators on $C([0,1])$
of the form $$L(f)(x)=\sum_{i\in I}f(a_i)\varphi_i(x)\quad(0\leq
x\leq 1),$$ where $(a_i)_{i\in I}$ is a finite family in $[0,1]$ and
$\varphi_i\in C([0,1])\;(i\in I)$. Finally, we point out that the
germ of the same theorem can be also traced back to a paper by
T.~Popoviciu ([95]).

Korovkin's theorem 3.1 (often called \dword{Korovkin's first
theorem}) has many important applications in the study of positive
approximation processes in $C([0,1])$.

One of them is concerned with the \dword{Bernstein operators} on
$C([0,1])$ which are defined by
\begin{equation}\label{3.5}
B_n(f)(x):=\sum_{k=0}^nf\Big(\frac{k}{n}\Big)\binom{n}{k}x^k(1-x)^{n-k}
\end{equation} $(n\geq 1, f\in C([0,1]),\; 0\leq x\leq 1)$.
Each $B_n(f)$ is a polynomial of degree not greater than $n$.  They
were introduced by S.~N.~Bernstein ([34]) to give the first
constructive proof of the Weierstrass approximation theorem
(algebraic version) ([119]).

Actually, we have that:
\begin{theorem}\sl  For every $f\in C([0,1])$,
$$\lim_{n\rightarrow\infty}B_n(f)=f \quad\text{uniformly on}\quad [0,1].$$

\end{theorem}
\begin{proof} Each $B_n$ is a positive linear operator on
$C([0,1])$.  Moreover, it is easy to verify that for any $n\geq 1$
$$B_n(\textbf{1})=\textbf{1},\quad B_n(e_1)=e_1$$ and
$$B_n(e_2)=\frac{n-1}{n}e_2+\frac{1}{n}e_1.$$ Therefore, the result follows
from Theorem 3.1.
\end{proof} The original proof of Bernstein's Theorem
3.6 is based on probabilistic considerations (namely, on the weak
law of large numbers).  For a survey on Bernstein operators, we
refer, e.g., to [82] (see also [42] and [8]).

Note that Theorem 3.6 furnishes a constructive proof of the
Weierstrass approximation theorem [119] which we state below. (For a
survey on many other alternative proofs of Weierstrass' theorem, we
refer, e.g., to [93-94].)
\begin{theorem}\sl  For every $f\in C([0,1])$, there exists a sequence of
algebraic polynomials that uniformly converges to $f$ on $[0,1]$.
\end{theorem} Using modern language, Theorem 3.7 can be restated as
follows
\begin{center} \textsl{``The subalgebra of all algebraic polynomials
is dense in $C([0,1])$\\ with respect to the uniform norm".}
\end{center}
By means of Theorem 3.6, we have seen that the Weierstrass
approximation theorem can be obtained from Korovkin's theorem.

It seems to be not devoid of interest to point out that, from the
Weierstrass theorem, it is possible to obtain a special version of
Korovkin's theorem which involves only positive linear
operators $L_n$, $n\geq 1$, such
that $L_n(C([0,1]))\subset B([0,1])$ for every $n\geq 1$.  This special version will be
referred to as the \dword{restricted version of Korovkin's theorem.}

\begin{theorem}\sl  The restricted version of Korovkin's theorem and Weierstrass' Approximation
Theorem are equivalent.
\end{theorem}

\begin{proof} We have to furnish a proof of the restricted version of Korovkin's theorem based
solely on the Weierstrass Theorem.

Consider a sequence of positive linear operators $(L_n)_{n\geq 1}$
from $C([0,1])$ into $B([0,1])$ such that
$\lim\limits_{n\rightarrow\infty}L_n(g)=g$ uniformly on $[0,1]$ for
every $g\in\{\textbf{1},e_1,e_2\}$.  As in the proof of Theorem 3.1,
we then get
$$\lim\limits_{n\rightarrow\infty}L_n(d_x^2)(x)=0$$ uniformly with respect
to $x\in[0,1]$.

For $m\geq 1$ and $x,y\in[0,1]$, we have $$|x^m-y^m|\leq m|y-x|$$
and hence, recalling the function $e_m(x)=x^m\quad (0\leq x\leq 1)$,
$$|e_m-y^m\textbf{1}|\leq m|e_1-y\textbf{1}| \quad (y\in[0,1]).$$

An application of the Cauchy-Schwarz inequality (\ref{2.16})
implies, for any $n\geq 1$ and $y\in[0,1]$, \\ $$\begin{array}{l}
|L_n(e_m)-y^mL_n(\textbf{1})|\leq mL_n(|e_1-y\textbf{1}|)\\
\\\kern3.5truecm  \leq
m\sqrt{L_n(\textbf{1})}\sqrt{L_n((e_1-y\textbf{1})^2)}=m\sqrt{L_n(\textbf{1})}\sqrt{L_n(d_y^2)}.
\end{array}$$ Therefore, $\lim\limits_{n\rightarrow\infty}L_n(e_m)=e_m$
uniformly on $[0,1]$ for any $m\geq 1$ and hence
$\lim\limits_{n\rightarrow\infty}L_n(P)=P$ uniformly on $[0,1]$ for
every algebraic polynomial $P$ on $[0,1]$.

We may now conclude the proof because, setting $M:=\underset{n\geq
1}{\sup}\,\norm{L_n}=\underset{n\geq
1}{\sup}\,\norm{L_n(\textbf{1})}<+\infty$ and fixing $f\in C([0,1])$
and $\varepsilon>0$, there exists an algebraic polynomial $P$ on
$[0,1]$ such that $\norm{f-P}\leq\varepsilon$, and an integer
$r\in\mathbb{N}$ such that $\norm{L_n(P)-P}\leq\varepsilon$ for
every $n\geq r$, so that
\begin{eqnarray*}
\norm{L_n(f)-f}\leq \norm{L_n(f)-L_n(P)}+\norm{L_n(P)-P}+\norm{P-f}\\
\leq M\norm{f-P}+\norm{L_n(P)-P}+\norm{P-f}\leq(M+2)\varepsilon.
\end{eqnarray*}
\vskip-\baselineskip
\end{proof}

For another proof of Korovkin's first theorem which involves
Weierstrass theorem, see [117].

It is well-known that there are ''trigonometric" versions of both
Korovkin's theorem and Weierstrass' theorem (see Theorems 4.3 and
4.6). Also, these versions are equivalent (see Theorem 4.7).  In the
sequel, we shall also prove that the generalizations of these two
theorems to compact and to locally compact settings are equivalent
as well (see Theorem 9.4).

We proceed now to illustrate another application of Korovkin's
theorem that concerns the approximation of functions in
$L^p([0,1]),\;1\leq p<+\infty$, by means of positive linear
operators.  Note that Bernstein operators are not suitable to
approximate Lebesgue integrable functions (see, for instance, [82,
Section 1.9]).

The space $C([0,1])$ is dense in $L^p([0,1])$ with respect to the
natural norm
\begin{equation}\label{3.6}
\norm{f}_p:=\Big(\int_0^1|f(t)|^p\dd t\Big)^{1/p}\quad (f\in
L^p([0,1]))
\end{equation} and
\begin{equation}\label{3.7}
\norm{f}_p\leq\norm{f}_\infty\quad\quad\text{if}\quad f\in C([0,1]).
\end{equation} Therefore, the subalgebra of all algebraic polynomials on
$[0,1]$ is dense in $L^p([0,1])$.

The \dword{Kantorovich polynomials} introduced by L.~V.~Kantorovich
([75]) furnish the first constructive proof of the above mentioned
density result. They are defined by
\begin{equation}\label{3.8}
K_n(f)(x):=\sum_{k=0}^n\Big[(n+1)\int_{\frac{k}{n+1}}^{\frac{k+1}{n+1}}f(t)\dd
t\Big]\binom{n}{k}x^k(1-x)^{n-k}
\end{equation} for every $n\geq 1, f\in L^p([0,1]), 0\leq x\leq 1$.
Each $K_n(f)$ is a polynomial of degree not greater than $n$ and
every $K_n$ is a positive linear operator from $L^p([0,1])$ (and, in
particular, from $C([0,1])$ into $C([0,1])$). For additional
information on these operators, see [8, Section 5.3.7], [82], [42,
Chapter 10].

\begin{theorem}\sl  If $f\in C([0,1])$, then
$$\lim_{n\rightarrow\infty}K_n(f)=f\quad\text{uniformly on}\quad [0,1].$$
\end{theorem}
\begin{proof} A direct calculation which involves the
corresponding formulas for Bernstein operators gives for $n\geq 1$
$$K_n(\textbf{1})=\textbf{1},\, K_n(e_1)=\frac{n}{n+1}e_1+\frac{1}{2(n+1)}$$
and
$$K_n(e_2)=\frac{n(n-1)}{(n+1)^2}e_2+\frac{2n}{(n+1)^2}e_1+\frac{1}{3(n+1)^2}.$$
Therefore, the result follows at once from Korovkin's Theorem 3.1.
\end{proof}

Before showing a result similar to Theorem 3.9 for $L^p$-functions,
we need to recall some properties of convex functions.

Consider a real interval $I$ of $\mathbb{R}$.  A function
$\varphi:I\longrightarrow\mathbb{R}$ is said to be \dword{convex} if
$$\varphi(\alpha x+(1-\alpha)y)\leq\alpha\varphi(x)+(1-\alpha)\varphi(y)$$
for every $x,y\in I$ and $0\leq\alpha\leq 1$. If $I$ is open and
$\varphi$ is convex, then, for every finite family $(x_k)_{1\leq
k\leq n}$ in $I$ and $(\alpha_k)_{1\leq k\leq n}$ in $[0,1]$ such
that $\sum\limits_{k=1}^n\alpha_k=1$,
$$\varphi\Big(\sum_{k=1}^n\alpha_kx_k\Big)\leq\sum_{k=1}^n\alpha_k\varphi(x_k)$$
(\dword{Jensen's inequality}).

The function $|t|^p\;(t\in\mathbb{R})$, $1\le p<\infty$, is convex.
Given a probability space $(\Omega,\mathcal{F},\mu)$, an open
interval $I$ of $\mathbb{R}$ and a $\mu$-integrable function
$f:\Omega\longrightarrow I$, then
$$\int_\Omega f\dd\mu\in I.$$
Furthermore, if $\varphi:I\longrightarrow\mathbb{R}$ is convex and
$\varphi\circ f:\Omega\longrightarrow\mathbb{R}$ is
$\mu$-integrable, then
\pagebreak
$$\varphi\Big(\int_\Omega f\dd\mu\Big)\leq\int_\Omega\varphi\circ f\dd\mu$$
(\dword{Integral Jensen inequality}).

In particular, if
$f\in\mathcal{L}^p(\Omega,\mu)\subset\mathcal{L}^1(\Omega,\mu)$,
then
\begin{equation}\label{3.9} \left|\int f\dd\mu\right|^p\leq\int|f|^p\dd\mu.
\end{equation} (For more details see, e.g., [29, pp.18--21].)

 After these
preliminaries, we now proceed to show the approximation property of
$(K_n)_{n\geq 1}$ in $L^p([0,1])$.

\begin{theorem}\sl  If $f\in L^p([0,1]),\;1\leq p<+\infty$, then
$$\lim_{n\rightarrow\infty}K_n(f)=f\quad\text{in}\quad L^p([0,1]).$$
\end{theorem}
\begin{proof}
For every $n\geq 1$, denote by $\norm{K_n}$ the operator norm of
$K_n$ considered as an operator from $L^p([0,1])$ into $L^p([0,1])$.

To prove the result, it is sufficient to show that there exists an
$M\geq 0$ such that $\norm{K_n}\leq M$ for every $n\geq 1$.  After
that, the result will follow immediately because, for a given
$\varepsilon>0$, there exists $g\in C([0,1])$ such that
$\norm{f-g}_p\leq\varepsilon$ and there exists $\nu\in\mathbb{N}$
such that, for $n\geq \nu$,
$$\norm{K_n(g)-g}_\infty\leq\varepsilon$$ so that $$\norm{K_n(f)-f}_p\leq
M\norm{f-g}_p+\norm{K_n(g)-g}_p+\norm{g-f}_p\leq(2+M)\varepsilon.$$
Now, in order to obtain the desired estimate, we shall use the
convexity of the function $|t|^p$ on $\mathbb{R}$ and inequality
(3.9).

Given $f\in L^p([0,1])$, for every $n\geq 1$ and $0\leq k\leq n$, we
have indeed
\begin{equation*}
\left((n+1)\int_{\frac{k}{n+1}}^{\frac{k+1}{n+1}}\mid f(t)\mid \dd
t\right)^{p}\leq (n+1)\int_{\frac{k}{n+1}}^{\frac{k+1}{n+1}}\mid
f(t)\mid^{p}\dd t
\end{equation*} and hence, for every $x\in \lbrack 0,1]$,
\begin{gather*}
|K_n(f)(x)|^p\leq\sum\limits_{k=0}^n\binom{n}{k}x^k(1-x)^{n-k}\Big[(n+1)\int_{\frac{k}{n+1}}^{\frac{k+1}{n+1}}
|f(t)|\dd t\Big]^p\\
\kern1.9truecm\leq\sum\limits_{k=0}^n\binom{n}{k}x^k(1-x)^{n-k}(n+1)\int_{\frac{k}{n+1}}^{\frac{k+1}{n+1}}|f(t)|^p
\dd t.
\end{gather*} 
Therefore,
$$\int_0^1|K_n(f)(x)|^p\dd x\leq\sum_{k=0}^n\binom{n}{k}\Big(\int_0^1x^k(1-x)^{n-k}\dd x\Big)\Big((n+1)
\int_{\frac{k}{n+1}}^{\frac{k+1}{n+1}}|f(t)|^p\dd t\Big).$$
On the other hand, by considering the beta function
$$B(u,v):=\int_0^1t^{u-1}(1-t)^{v-1}\dd t\quad\quad (u>0, v>0),$$ it is not
difficult to show that, for $0\leq k\leq n$,
$$\int_0^1x^k(1-x)^{n-k}\dd x=B(k+1,n-k+1)=\frac{1}{(n+1)\binom{n}{k}},$$ and
hence
$$\int_0^1|K_n(f)(x)|^p\dd x\leq\sum_{k=0}^n\int_{\frac{k}{n+1}}^{\frac{k+1}{n+1}}|f(t)|^p\dd t=\int_0^1|f(t)|^p
\dd t.$$
Thus, $\norm{K_n(f)}_p\leq\norm{f}_p$ for every $f\in L^p([0,1])$,
i.e., $\norm{K_n}\leq 1$.
\end{proof}
\begin{remarks}\rm
\begin{itemize}
\item[]
\item [1.]For every $f\in L^p([0,1]),\;1\leq p<+\infty$, it can also be
shown that $$\lim_{n\rightarrow\infty}K_n(f)=f\quad\text{almost
everywhere on}\quad [0,1]$$ (see [82, Theorem 2.2.1]). \item [2.] If
$f\in C([0,1])$ is continuously differentiable in [0,1], then, by
referring again to Bernstein operators (\ref{3.5}), it is not
difficult to show that, for $n\geq 1$ and $x\in[0,1]$,
$$B_{n+1}(f)'(x)=\sum_{h=0}^n(n+1)\Big[f\Big(\frac{h+1}{n+1}\Big)-f
\Big(\frac{h}{n+1}\Big)\Big]\binom{n}{h}x^h(1-x)^{n-h}=K_n(f')(x).$$
Therefore, by Theorem 3.9, we infer that
\begin{equation}\label{3.10}
  \lim_{n\rightarrow\infty}B_n(f)'=f'\quad\text{uniformly on}\quad [0,1].
\end{equation}
More generally, if $f\in C([0,1])$ possesses continuous derivatives
in [0,1] up to the order $m\geq 1$, then for every $1\leq k\leq m$,
\begin{equation}\label{3.11}
\lim_{n\rightarrow\infty}B_n(f)^{(k)}=f^{(k)}\quad\text{uniformly
on}\quad [0,1]
\end{equation}
([82, Section 1.8]).
\item [3.] Another
example of positive approximating operators on $L^p([0,1]),\;1\leq
p<+\infty$, is
furnished by the \dword{Bernstein-Durrmeyer operators} defined by
\begin{equation}\label{3.12}
D_n(f)(x):=\sum\limits_{k=0}^n\Big(\int\limits_0^1(n+1)\binom{n}{k}t^k(1-t)^{n-k}f(t)
  \dd t\Big)\binom{n}{k}x^k(1-x)^{n-k}
\end{equation}
$(f\in L^p([0,1]),\;0\leq x\leq 1)$\quad(\text{see [53], [40], [8,
  Section 5.3.8]}).
\end{itemize} We also refer the interested reader to [16] where a
generalization of Kantorovich operators is introduced and studied.
\end{remarks}

\section{Korovkin's second theorem and something else} In this section, we
shall consider the space $\mathbb{R}^d$, $d\geq 1$, endowed with the
Euclidean norm
\begin{equation}\label{4.1}
\norm{x}=\Big(\sum_{i=1}^dx_i^2\Big)^{1/2}\quad(x=(x_i)_{1\leq i\leq
d}\in\mathbb{R}^d).
\end{equation} For every $j=1,\dots,d$, we shall denote
by $$pr_j:\mathbb{R}^d\longrightarrow\mathbb{R}$$ the $j$-th
coordinate function which is defined by
\begin{equation}\label{4.2} pr_j(x):=x_j\quad
(x=(x_i)_{1\leq i\leq d}\in\mathbb{R}^d).
\end{equation} By a common abuse
of notation, if $X$ is a subset of $\mathbb{R}^d$, the restriction
of each $pr_j$ to $X$ will be again denoted by $pr_j$.  In  this
framework, for the functions $d_x$ $(x\in X)$ defined by
(\ref{3.4}), we get
\begin{equation}\label{4.3}
d_x^2=\norm{x}^2\textbf{1}-2\sum_{i=1}^dx_ipr_i+\sum_{i=1}^dpr_i^2.
\end{equation} Therefore, from Theorem 3.5, we then obtain
\begin{theorem}\sl
Let $X$ be a locally compact subset of $\mathbb{R}^d$, $d\geq 1$,
i.e.,  $X$ is the intersection of an open subset and a closed subset
of $\mathbb{R}^d$ (see \Appendix).  Consider a lattice subspace $E$
of $F(X)$ containing
$\{\textbf{1},pr_1,\dots,pr_d,\sum\limits_{i=1}^dpr_i^2\}$ and let
$(L_n)_{n\geq 1}$ be a sequence of positive linear operators from
$E$ into $F(X)$ such that for every
$g\in\{\textbf{1},pr_1,\dots,pr_d,\sum\limits_{i=1}^dpr_i^2\}$
$$\lim_{n\rightarrow\infty}L_n(g)=g\quad\text{uniformly on compact subsets
of}\quad X.$$ Then, for every $f\in E\cap C_b(X)$ $$
\lim_{n\rightarrow\infty}L_n(f)=f \quad\text{uniformly on compact
subsets of}\quad X.$$
\end{theorem} The special case of Theorem 4.1 when $X$ is
compact follows indeed from Theorem 3.2 and is worth being stated
separately.  It is due to Volkov ([118]).
\begin{theorem}\sl  Let $X$ be a
compact subset of $\mathbb{R}^d$ and consider a sequence
$(L_n)_{n\geq 1}$ of positive linear operators from $C(X)$ into
$F(X)$ such that for every
$g\in\{\textbf{1},pr_1,\dots,pr_d,\sum\limits_{i=1}^dpr_i^2\}$
$$\lim_{n\rightarrow\infty}L_n(g)=g\quad\text{uniformly on} \quad X.$$ Then
for every $f\in C(X)$ $$
\lim_{n\rightarrow\infty}L_n(f)=f\quad\text{uniformly on} \quad X.$$
\end{theorem} Note that, if $X$ is contained in some sphere of
$\mathbb{R}^d$, i.e.,  $\sum\limits_{i=1}^dpr_i^2$ is constant on
$X$, then the test subset in Theorem 4.2 reduces to $\{\textbf{1},
pr_1,\dots,pr_d\}$.  (In [8, Corollary 4.5.2], the reader can find a
complete characterization of those subsets $X$ of $\mathbb{R}^d$ for
which $\{\textbf{1},pr_1,\dots,pr_d\}$ satisfies Theorem 4.2.)

This remark applies in particular for the unit circle of
$\mathbb{R}^2$
\begin{equation}\label{4.4}
\mathbb{T}:=\{(x,y)\in\mathbb{R}^2\mid x^2+y^2=1\}.
\end{equation} On the
other hand, the space $C$($\mathbb{T}$) is isometrically (order)
isomorphic to the space
\begin{equation}\label{4.5}
C_{2\pi}(\mathbb{R}):=\{f\in C(\mathbb{R})\mid f\quad\text{is}\quad
2\pi\text{-periodic}\}
\end{equation} (endowed with the sup-norm and
pointwise ordering) by means of the isomorphism
$\Phi:C(\mathbb{T})\longrightarrow C_{2\pi}(\mathbb{R})$ defined by
\begin{equation}\label{4.6} \Phi(F)(t):=F(\cos t, \sin t)\quad
(t\in\mathbb{R})\,.
\end{equation} Moreover,
\begin{equation}\label{4.7}
\Phi(\textbf{1})=\textbf{1},\; \Phi(pr_1)=\cos,\; \Phi(pr_2)=\sin
\end{equation}
and so we obtain \dword{Korovkin's second theorem}.

\begin{theorem}\sl  Let $(L_n)_{n\geq 1}$ be a sequence of positive linear
operators from $C_{2\pi}(\mathbb{R})$ into $F(\mathbb{R})$ such that
$$\lim_{n\rightarrow\infty}L_n(g)=g\quad\text{uniformly on}\quad
\mathbb{R}$$ for every  $g\in\{\textbf{1},\cos,\sin\}$.  Then
$$\lim_{n\rightarrow\infty}L_n(f)=f\quad\text{uniformly on}\quad
\mathbb{R}$$ for every $f\in C_{2\pi}(\mathbb{R})$.
\end{theorem} Below, we
discuss some applications of Theorem 4.3.

For $1\leq p<+\infty$, we shall denote by $$L_{2\pi}^p(\mathbb{R})$$
the Banach space of all (equivalence classes of) functions
$f:\mathbb{R}\longrightarrow\mathbb{R}$ that are Lebesgue integrable
to the $p$-th power over $[-\pi,\pi]$ and that satisfy
$f(x+2\pi)=f(x)$ for a.e.\ $x\in\mathbb{R}$. The space
$L_{2\pi}^p(\mathbb{R})$ is endowed with the norm
\begin{equation}\label{4.8}
\norm{f}_p:=\Big(\frac{1}{2\pi}\int_{-\pi}^\pi|f(t)|^p\dd
t\Big)^{1/p}\quad (f\in L_{2\pi}^p(\mathbb{R})).
\end{equation}

A family $(\varphi_n)_{n\geq 1}$ in $L_{2\pi}^1(\mathbb{R})$ is said
to be \dword{a positive periodic kernel} if every $\varphi_n$ is
positive, i.e., $\varphi_n\geq 0$ a.e.\  on $\mathbb{R}$, and
\begin{equation}\label{4.9}
\lim_{n\rightarrow\infty}\frac{1}{2\pi}\int_{-\pi}^\pi\varphi_n(t)\dd
t=1.
\end{equation} Each positive kernel $(\varphi_n)_{n\geq 1}$ generates a
sequence of positive linear operators on $L_{2\pi}^1(\mathbb{R})$.
For every $n\geq 1, f\in L_{2\pi}^1(\mathbb{R})$ and
$x\in\mathbb{R}$, set
\begin{eqnarray}\label{4.10}
L_n(f)(x):=(f*\varphi_n)(x)=\frac{1}{2\pi}\int_{-\pi}^{\pi}f(x-t)\varphi_n(t)\dd
t\\\nonumber =\frac{1}{2\pi}\int_{-\pi}^{\pi}f(t)\varphi_n(x-t)\dd
t.
\end{eqnarray} From
Fubini's theorem and H\"{o}lder's inequality, it follows that
$L_n(f)\in L_{2\pi}^p(\mathbb{R})$ if $f\in L_{2\pi}^p(\mathbb{R}),
1\leq p<+\infty$.

Moreover, if $f\in C_{2\pi}(\mathbb{R})$, then the Lebesgue
dominated convergence theorem implies that $L_n(f)\in
C_{2\pi}(\mathbb{R})$. Furthermore,
\begin{equation}\label{4.11}
\norm{L_n(f)}_p\leq\norm{\varphi_n}_1\norm{f}_p\quad\quad (f\in
C_{2\pi}(\mathbb{R}))
\end{equation} and
\begin{equation}\label{4.12}
\norm{L_n(f)}_{\infty}\leq\norm{\varphi_n}_1\norm{f}_{\infty}\quad\quad
(f\in C_{2\pi}(\mathbb{R})).
\end{equation} A positive kernel
$(\varphi_n)_{n\geq 1}$ is called an \dword{approximate identity} if
for every $\delta\in]0,\pi[$
\begin{equation}\label{4.13}
\lim_{n\rightarrow\infty}\int_{-\pi}^{-\delta}\varphi_n(t)\dd
t+\int_{\delta}^{\pi}\varphi_n(t)\dd t=0\,.
\end{equation}
\begin{theorem}\sl
Consider a positive kernel $(\varphi_n)_{n\geq 1}$ in
$L_{2\pi}^1(\mathbb{R})$ and the corresponding sequence
$(L_n)_{n\geq 1} $ of positive linear operators defined by
(\ref{4.10}).  For every $n\geq 1$, set
\begin{equation}\label{4.14}
\beta_n:=\frac{1}{2\pi}\int_{-\pi}^{\pi}\varphi_n(t)\sin^2\frac{t}{2}\dd
t.
\end{equation} Then the following properties are equivalent:
\begin{itemize}
\item [a)] For every $1\leq p<+\infty$ and $f\in L_{2\pi}^p(\mathbb{R})$
$$\lim_{n\rightarrow\infty}L_n(f)=f\quad \text{in}\quad
L_{2\pi}^p(\mathbb{R})$$ as well as
$$\lim_{n\rightarrow\infty}L_n(f)=f\quad \text{in}\quad
C_{2\pi}(\mathbb{R})$$ provided $f\in C_{2\pi}(\mathbb{R})$.
\item [b)] $\lim\limits_{n\rightarrow\infty}\beta_n=0$.
\item [c)] $(\varphi_n)_{n\geq 1}$ is an approximate identity.
\end{itemize}
\end{theorem}
\begin{proof}
To show the implication $(a)\Rightarrow (b)$, it is sufficient to
point out that for every $n\geq 1$ and $x\in\mathbb{R}$
\begin{eqnarray*}
\beta_n&=&\frac{1}{2\pi}\int_{-\pi}^{\pi}\varphi_n(u-x)\sin^2\frac{u-x}{2}\dd u\\
&=&\frac{1}{2}\Big(\frac{1}{2\pi}\int_{-\pi}^{\pi}\varphi_n(t)\dd
t-(\cos x) L_n(\cos)(x)-(\sin x) L_n(\sin)(x)\Big)
\end{eqnarray*} and hence
$\beta_n\rightarrow 0$ as $n\rightarrow\infty$.

Now assume that (b) holds.  Then, for $0<\delta<\pi$ and $n\geq 1$,\\
\begin{gather*}
\displaystyle{\frac{\sin^2(\delta/2)}{2\pi}\Big(\int\limits_{-\pi}^{-\delta}
\varphi_n(t)\dd t+\int\limits_{\delta}^{\pi}\varphi_n(t)\dd t\Big)}
\displaystyle{\;\leq\;\frac{1}{2\pi}\Big(\int\limits_{-\pi}^{-\delta}\varphi_n(t)\sin^2\frac{t}{2}\dd
t+\int\limits_{\delta}^{\pi} \varphi_n(t)\sin^2\frac{t}{2}\dd
t\Big)\leq\beta_n}
\end{gather*} and hence
(c) follows.

We now proceed to show the implication $(c)\Rightarrow(a)$.  Set
$M:=\underset{n\geq 1}{\sup}\int\limits_{-\pi}^{\pi}\varphi_n(t)\dd
t$.  For a given $\varepsilon>0$ there exists $\delta\in]0,\pi[$
such that $|\cos t-1|\leq\displaystyle{\frac{\varepsilon}{6(M+1)}}$
and $|\sin t|\leq\displaystyle{\frac{\varepsilon}{6(M+1)}}$ for any
$t\in\mathbb{R}, |t|\leq\delta$, and hence, for sufficiently large
$n\geq 1$
$$\Big|\frac{1}{2\pi}\int_{-\pi}^{\pi}\varphi_n(t)\dd t-1\Big|\leq\frac{\varepsilon}{3}\quad\text{and}\quad
\int_{\delta\leq |t|\leq\pi}\varphi_n(t)\dd
t\leq\pi\frac{\varepsilon}{3}.$$

Therefore, for every $x\in\mathbb{R}$,\\
\begin{equation*} \begin{split}
&|L_n(\sin)(x)-\sin
x|\leq\frac{1}{2\pi}\int\limits_{-\pi}^{\pi}|\sin(x-t)-\sin
x|\varphi_n(t)\dd t\\&
\kern4truecm+\Big|\Big(\frac{1}{2\pi}\int_{-\pi}^{\pi}\varphi_n(t)\dd
t-1\Big)\Big||\sin x|\\&
\leq\frac{1}{2\pi}\int\limits_{\delta\leq|t|\leq\pi}|\sin(x-t)-\sin
x|\varphi_n(t)\dd t\\&\kern1truecm
+\frac{1}{2\pi}\int\limits_{|t|<\delta}|(\cos t-1)\sin x-\cos x\sin
t|\varphi_n(t)\dd t+\varepsilon/3\\&
\leq\frac{1}{\pi}\int\limits_{\delta\leq|t|\leq\pi}\varphi_n(t)\dd
t+\frac{\varepsilon}{3(M+1)}
\frac{1}{2\pi}\int\limits_{-\pi}^{\pi}\varphi_n(t)\dd
t+\frac{\varepsilon}{3} \;\;\leq\;\;\varepsilon.
\end{split}
\end{equation*}\\

Therefore, $\lim\limits_{n\rightarrow\infty}L_n(\sin)=\sin$
uniformly on $\mathbb{R}$.  The same method can be used to show that
$\lim\limits_{n\rightarrow\infty}L_n(\cos)=\cos$ uniformly on
$\mathbb{R}$ and hence, by Korovkin's second theorem 4.3, we obtain
$\lim\limits_{n\rightarrow\infty}L_n(f)=f$ in $C_{2\pi}(\mathbb{R})$
for every $f\in C_{2\pi}(\mathbb{R})$.

By reasoning as in the proof of Theorem 3.10, it is a simple matter
to get the desired convergence formula in $L_{2\pi}^p(\mathbb{R})$
by using the previous one on $C_{2\pi}(\mathbb{R})$, the denseness
of $C_{2\pi}(\mathbb{R})$ in $L_{2\pi}^p(\mathbb{R})$ and formula
(\ref{4.11}) which shows that the operators $L_n$, $n\geq 1$, are
equibounded from $L_{2\pi}^p(\mathbb{R})$ into
$L_{2\pi}^p(\mathbb{R})$.
\end{proof} Two simple applications of Theorem 4.4 are particularly worthy
of mention.  For other applications, we refer to [8, Section 5.4],
[38], [41], [78].

We begin by recalling that a trigonometric polynomial of degree
$n\in\mathbb{N}$ is a real function of the form
\begin{equation}\label{4.15}
u_n(x)=\frac{1}{2}a_0+\sum_{k=1}^n a_k\cos kx+b_k\sin kx\quad
(x\in\mathbb{R})
\end{equation} where $a_0, a_1,\dots, a_n, b_1,\dots,
b_n\in\mathbb{R}$.  A series of the form
\begin{equation}\label{4.16}
\frac{1}{2}a_0+\sum_{k=1}^\infty a_k\cos kx+b_k\sin kx \quad
(x\in\mathbb{R})
\end{equation} $(a_k,b_k\in\mathbb{R})$ is called a
\dword{trigonometric series}.

If $f\in L_{2\pi}^1(\mathbb{R})$, the trigonometric series
\begin{equation}\label{4.17} \frac{1}{2}a_0(f)+\sum_{k=1}^\infty a_k(f)\cos
kx+b_k(f)\sin kx \quad (x\in\mathbb{R})
\end{equation} where
\begin{equation}\label{4.18} a_0(f):=\frac{1}{\pi}\int_{-\pi}^{\pi}f(t)\dd t,
\end{equation}
\begin{equation}\label{4.19}
a_k(f):=\frac{1}{\pi}\int_{-\pi}^{\pi}f(t)\cos kt\dd t,\quad k\geq
1,
\end{equation}
\begin{equation}\label{4.20}
b_k(f):=\frac{1}{\pi}\int_{-\pi}^{\pi}f(t)\sin kt\dd t,\quad k\geq
1,
\end{equation} is called the \dword{Fourier series} of $f$.  The $a_n$'s and
$b_n$'s are called the \dword{real Fourier coefficients} of $f$.

For any $n\in\mathbb{N}$, denote by $$S_n(f)$$ the $n$-th partial
sum of the Fourier series of $f$, i.e.,
\begin{equation}\label{4.21}
S_0(f)=\frac{1}{2}a_0(f)
\end{equation} and, for $n\geq 1$,
\begin{equation}\label{4.22}
S_n(f)(x)=\frac{1}{2}a_0(f)+\sum_{k=1}^na_k(f)\cos kx+b_k(f)\sin kx.
\end{equation} Each $S_n(f)$ is a trigonometric polynomial; moreover,
considering the functions
\begin{equation}\label{4.23}
D_n(t):=1+2\sum_{k=1}^n\cos kt\quad (t\in\mathbb{R}),
\end{equation} we
also get
\begin{equation}\label{4.24}
S_n(f)(x)=\frac{1}{2\pi}\int_{-\pi}^{\pi}f(t)D_n(x-t)\dd
t\quad(x\in\mathbb{R}).
\end{equation} The function $D_n$ is called the \dword{$n$-th Dirichlet
kernel}.

By multiplying (\ref{4.23}) by $\sin t/2$, we obtain\\
\begin{eqnarray*}
\sin\frac{t}{2}D_n(t)&=&\sin\frac{t}{2}+\sum_{k=1}^n\sin\Big(\frac{1+2k}{2}t\Big)-\sin
\Big(\frac{2k-1}{2}t\Big)\\
&=&\sin\frac{1+2n}{2}t,
\end{eqnarray*} so that
\begin{equation}\label{4.25} D_n(t)=\begin{cases}
\displaystyle{\frac{\sin(1+2n)t/2}{\sin t/2}}&\text{if}\;t\;\text{is
not a multiple of}\;\pi,\\ 2n+1&\text{if}\;t\;\text{is  a multiple
of}\;\pi.
\end{cases}
\end{equation} $D_n$ is not positive and $(D_n)_{n\geq 1}$ is
not an approximate identity ([38, Prop.~1.2.3]).  Moreover, there
exists $f\in C_{2\pi}(\mathbb{R})$ such that $(S_n(f))_{n\geq 1}$
does not converge uniformly (nor pointwise) to $f$, i.e.,  the
Fourier series of $f$ does not converge uniformly (nor pointwise) to
$f$.

For every $n\in\mathbb{N}$, put
\begin{equation}\label{4.26}
F_n(f):=\frac{1}{n+1}\sum_{k=0}^nS_k(f).
\end{equation} $F_n(f)$ is a
trigonometric polynomial.  Moreover, from the identity
\begin{equation}\label{4.27}
(\sin\frac{t}{2})\sum_{k=0}^{n-1}\sin\frac{2k+1}{2}
t=\sin^2\frac{n}{2}t\quad (t\in\mathbb{R}),
\end{equation} it follows that for every
$x\in\mathbb{R}$\\
\begin{eqnarray}\label{4.28}
F_n(f)(x)&=&\frac{1}{2\pi}\int_{-\pi}^{\pi}f(t)\frac{1}{(n+1)}\sum_{k=0}^n\frac{\sin((2k+1)(x-t)/2)}{\sin((x-t)/2)}\dd t\\
\nonumber &=&\frac{1}{2\pi}
\int_{-\pi}^{\pi}f(t)\frac{1}{(n+1)}\frac{\sin^2((n+1)(x-t)/2)}{\sin^2((x-t)/2)}\dd t\\
\nonumber &=&\frac{1}{2\pi} \int_{-\pi}^{\pi}f(t)\varphi_n(x-t)\dd
t,
\end{eqnarray}

where
\begin{equation}\label{4.29}
\varphi_n(x):=\begin{cases}
\frac{\sin^2((n+1)x/2)}{(n+1)\sin^2(x/2)}&\text{if}\;x\;\text{is not
a multiple of}\;2\pi,\\ n+1&  \text{if}\;x\;\text{is  a multiple
of}\;2\pi.
\end{cases}
\end{equation} Actually, the sequence $(\varphi_n)_{n\geq 1}$
is a positive kernel which is called the \dword{Fej\'{e}r kernel},
and the corresponding operators $F_n$, $n\geq 1$, are called the
\dword{Fej\'{e}r convolution operators}.
\begin{theorem}\sl  For every $f\in
L_{2\pi}^p(\mathbb{R}),\; 1\leq p<+\infty$,
$$\lim_{n\rightarrow\infty}F_n(f)=f\quad \text{in}\quad
L_{2\pi}^p(\mathbb{R})$$ and, if $f\in C_{2\pi}(\mathbb{R})$,
$$\lim_{n\rightarrow\infty}F_n(f)=f\quad \text{in}\quad
C_{2\pi}(\mathbb{R}).$$
\end{theorem}
\begin{proof} Evaluating the Fourier
coefficients of $\textbf{1}$, $\cos$ and $\sin$, and by using
(\ref{4.26}), we obtain, for $n\geq 1$,
$$F_n(\textbf{1})=\textbf{1},\ F_n(\cos)=\frac{n}{n+1}\cos,\
F_n(\sin)=\frac{n}{n+1}\sin$$ and
$$\beta_n=\frac{1}{2(n+1)}.$$ The result now follows from Theorem 4.4 or,
more directly, from Theorem 4.3.
\end{proof}

Theorem 4.5 is due to Fej\'{e}r ([56-57]).  It furnishes the first
constructive proof of \eword{the Weierstrass approximation theorem
for periodic functions}.

\begin{theorem}\sl  If $f\in L_{2\pi}^p(\mathbb{R})$, $1\leq p<+\infty$ (resp.\
$f\in C_{2\pi}(\mathbb{R})$) then  there exists a sequence of
trigonometric polynomials that converges to $f$ in
$L_{2\pi}^p(\mathbb{R})$ (resp. uniformly on $\mathbb{R}$).
\end{theorem} As in the ``algebraic" case, we
shall now show that from Weierstrass' approximation theorem, it is
possible to deduce a ``restricted" version of Theorem 4.3, where in
addition it is required that each operator $L_n$ maps
$C_{2\pi}(\mathbb{R})$ into $B(\mathbb{R})$.  We shall also refer to
this version as the \dword{restricted version of Korovkin's second
theorem}.

\begin{theorem}\sl  The restricted version of Korovkin's second Theorem 4.3 and Weierstrass' Theorem
4.6 are equivalent.
\end{theorem}
\begin{proof} An inspection of the proof
of Theorem 4.5 shows that Theorem 4.3 implies Theorem 4.5 and,
hence, Theorem 4.6.

Conversely, assume that Theorem 4.6 is true and consider a sequence
$(L_n)_{n\geq 1}$ of positive linear operators from
$C_{2\pi}(\mathbb{R})$ into $B(\mathbb{R})$ such that
$L_n(g)\longrightarrow g$ uniformly on $\mathbb{R}$ for every
$g\in\{\textbf{1},\cos,\sin\}$. For every $m\geq 1$, set
$f_m(x):=\cos mx$ and $g_m(x):=\sin mx$ $(x\in\mathbb{R})$. Since
the subspace of all trigonometric polynomials is dense in
$C_{2\pi}(\mathbb{R})$ and $$\underset{n\geq
1}{\sup}\,\norm{L_n}=\underset{n\geq
1}{\sup}\,\norm{L_n(\textbf{1})}<+\infty,$$ it is enough to show
that $L_n(f_m)\rightarrow f_m$ and $L_n(g_m)\rightarrow g_m$
uniformly on $\mathbb{R}$ for every $m\geq 1$.

Given $x\in\mathbb{R}$, consider the function
$\Phi_x(y)=\sin\frac{x-y}{2}\;(y\in\mathbb{R})$.  Then
$$\Phi^2(y):=\sin^2\frac{x-y}{2}=\frac{1}{2}(1-\cos x\cos y-\sin x\sin
y)\quad (y\in\mathbb{R})$$ and hence $$L_n(\Phi_x)(x)\rightarrow
0\quad\text{uniformly with respect to}\quad x\in\mathbb{R}.$$ On the
other hand, for $m\geq 1$ and $x,y\in\mathbb{R}$, we get
$$|f_m(x)-f_m(y)|=2\Big|\sin m\Big(\frac{x+y}{2}\Big)\Big|\Big|\sin
m\Big(\frac{x-y}{2}\Big)\Big|\leq
c_m\Big|\sin\Big(\frac{x-y}{2}\Big)\Big|$$ where $c_m:=2\,
\sup_{\alpha\in\mathbb{R}}\,\frac{\sin m\alpha}{\sin\alpha}$, and
hence
$$|f_m(x)L_n(\textbf{1})-L_n(f_m)|\leq c_mL_n(|\Phi_x|)\leq
c_m\sqrt{L_n(\textbf{1})}\sqrt{L_n(\Phi_x^2)}.$$ Therefore,
$L_n(f_m)(x)-f_m(x)\rightarrow 0$ uniformly with respect to
$x\in\mathbb{R}$.

A similar reasoning can be applied also to the functions $g_m, m\geq
1$, because $$|g_m(x)-g_m(y)|=2\Big|\cos
m\Big(\frac{x+y}{2}\Big)\sin m\Big(\frac{x-y}{2}\Big)\Big|\leq K_m
\Big|\sin\Big(\frac{x-y}{2}\Big)\Big|$$ and this completes the
proof.
\end{proof} For another short proof of Korovkin's second theorem
that uses the trigonometric version of Weierstrass' theorem, see
[117].

Fej\'{e}r's Theorem 4.5 is noteworthy because it reveals an
important property of the Fourier series.  Actually, it shows that
Fourier series are always \eword{Cesaro-summable} to $f$ in
$L_{2\pi}^p(\mathbb{R})$ or in $C_{2\pi}(\mathbb{R})$, provided that
$f\in L_{2\pi}^p(\mathbb{R})$ or $f\in C_{2\pi}(\mathbb{R})$.
Another deeper theorem, ascribed to Fej\'{e}r and Lebesgue, states
that, if $f\in L_{2\pi}^1(\mathbb{R})$, then its Fourier series is
Cesaro-summable to $f$ a.e.\ on $\mathbb{R}$ ([112, Theorem 8.35]).

Below, we further discuss another regular summation method, namely
the \dword{Abel-summation method}, which applies to Fourier series.

We begin with the following equality
\begin{equation}\label{4.30}
\frac{1+z}{1-z}=1+2\sum_{k=1}^\infty z^k\quad (z\in\mathbb{C},
|z|<1)
\end{equation} which holds uniformly on any compact subset of
$\{z\in\mathbb{C}\mid |z|<1\}$. Given $x\in\mathbb{R}$ and $0\leq
r<1$, by applying (\ref{4.30}) to $z=r\ee^{\ii x}$ and by taking the
real parts of both sides, we get
\begin{equation}\label{4.31} 1+2\sum_{k=1}^\infty r^k\cos
kx=\frac{1-r^2}{1-2r\cos x+r^2}
\end{equation} and the identity holds
uniformly when $r$ ranges in a compact interval of $[0,1[$.

The family of functions
\begin{equation}\label{4.32}
P_r(t):=\frac{1-r^2}{1-2r\cos t+r^2}\quad (t\in\mathbb{R})
\end{equation}
$(0\leq r<1)$ is called the \dword{Abel-Poisson kernel} and the
corresponding operators
\begin{equation}\label{4.33}
P_r(f)(x):=\frac{1-r^2}{2\pi}\int_{-\pi}^{\pi}\frac{f(t)}{1-2r\cos(x-t)+r^2}\dd
t\quad (x\in\mathbb{R})
\end{equation} $(0\leq r<1,\;f\in L_{2\pi}^1(\mathbb{R}))$
are called the \dword{Abel-Poisson convolution operators}. Taking
(\ref{4.31}) into account, it is not difficult to show that
\begin{equation}\label{4.34} P_r(f)(x)=\frac{1}{2}a_0(f)+\sum_{k=1}^\infty
r^k(a_k(f)\cos kx+b_k(f)\sin kx)
\end{equation} where the coefficients
$a_k(f)$ and $b_k(f)$ are the Fourier coefficients of $f$ defined by
(\ref{4.18})--(\ref{4.20}).

\begin{theorem}\sl  If $f\in
L_{2\pi}^p(\mathbb{R}), 1\leq p<+\infty$, then  $$\lim_{r\rightarrow
1^{-}}P_r(f)=f\quad\text{in}\quad L_{2\pi}^p(\mathbb{R})$$ and, if
$f\in C_{2\pi}(\mathbb{R})$, $$\lim_{r\rightarrow
1^{-}}P_r(f)=f\quad\text{uniformly on}\quad \mathbb{R}.$$
\end{theorem}
\begin{proof} The kernels $p_r$, $0\leq r<1$, are positive.  Moreover, from
(\ref{4.34}), we get
$$P_r(\textbf{1})=\textbf{1},\;P_r(\cos)=r\cos,\;P_r(\sin)=r\sin$$ and
hence $$\beta_r=\frac{1-r}{2}.$$ Therefore, the result follows from
Theorem 4.4 (or from Theorem 4.5).
\end{proof} According to (\ref{4.34}), Theorem
4.8 claims that the Fourier series of a function $f\in
L_{2\pi}^p(\mathbb{R})$ (resp.  $f\in C_{2\pi}(\mathbb{R})$) is Abel
summable to $f$ in $L_{2\pi}^p(\mathbb{R})$ (resp.  uniformly on
$\mathbb{R}$).

For further applications of Korovkin's second theorem to
approximation by convolution operators and summation processes, we
refer, e.g., to [8, Section 5.4], [38], [41], [78].

We finally point out the relevance of Theorem 4.8 in the study of
\dword{the Dirichlet problem for the unit disk}
$\mathbb{D}:=\{(x,y)\in\mathbb{R}^2\mid x^2+y^2\leq 1\}$. Given
$F\in C(\partial\mathbb{D})=C(\mathbb{T})\equiv
C_{2\pi}(\mathbb{R})$, this problem consists in finding a function
$U\in C(\mathbb{D})$ possessing second partial derivatives on the
interior of $\mathbb{D}$ such that
\begin{equation}\label{4.35}\begin{cases}
\displaystyle{\frac{\partial^2 U}{\partial
x^2}(x,y)+\frac{\partial^2 U}{\partial y^2}(x,y)}=0& (x^2+y^2<1),\\
U(x,y)=F(x,y)& (x^2+y^2=1).
\end{cases}\end{equation}
By using polar coordinates $x=r\cos\theta$ and $y=r\sin\theta
\;(0\leq r< 1,\, \theta\in\mathbb{R})$ and the functions
$f(\theta):=F(\cos\theta, \sin\theta)\;(\theta\in\mathbb{R})$ and
$u(r,\theta ):=U(r\cos \theta ,r\sin \theta)\;(0\leq r< 1,\,
\theta\in\mathbb{R})$, problem (\ref{4.35}) turns into
\begin{equation}\label{4.36}\begin{cases}
\displaystyle{\frac{\partial^2u}{\partial r^2}(r, \theta)
+\frac{1}{r}\frac{\partial u}{\partial r}(r,
\theta)+\frac{1}{r^2}\frac{\partial^2u}{\partial\theta^2}}(r,
\theta)=0&0<r<1,\;\theta\in\mathbb{R},\\
u(0,\theta)=\frac{1}{2\pi}\int\limits_{-\pi}^{\pi}f(t)\dd t&\theta\in\mathbb{R},\\
\lim\limits_{r\rightarrow
1^{-}}u(r,\theta)=f(\theta)&\text{uniformly \quad w.r.t.}
\quad\theta\in\mathbb{R}.
\end{cases}\end{equation}.

With the help of Theorem 4.8 it is not difficult to show that a
solution to problem (\ref{4.36}) is given by
\begin{equation}\label{4.37}
u(r,\theta)=P_r(f)(\theta)\quad\quad(0\leq r<1,\;
\theta\in\mathbb{R}).
\end{equation}
Accordingly, the function
\[
U(x,y):=\left\{
\begin{array}{l}
u(r,\theta )\text{ \ \ \ if \ }x=r\cos \theta ,\text{ }y=r\sin
\theta
\text{ and }x^{2}+y^{2}<1, \\
F(x,y)\text{ \ \ if \ }x^{2}+y^{2}=1,%
\end{array}%
\right.
\]
is a solution to problem (\ref{4.35}). Furthermore, as it is
well-known, the solution to (\ref{4.35}) is unique.
A similar result also holds if $F:\mathbb{T\rightarrow R}$ is a Borel-measurable function such that $\frac{1%
}{2\pi }\int_{-\pi }^{\pi }\mid F(\cos t,\sin t)\mid ^{p}dt<+\infty ,$ $%
1\leq p<+\infty $. For more details we refer, e.g., to [38,
Proposition 1.2.10 and Theorem 7.1.3] and to [96, Section 1.2].

\indent We end this section by presenting an application of Theorem
4.1. Consider the $d$-dimensional simplex
\begin{equation}\label{4.38}
K_d:=\{x=(x_i)_{1\leq i\leq d}\in\mathbb{R}^d\mid x_i\geq 0, 1\leq
i\leq d,\text{and} \sum_{i=1}^dx_i\leq 1\}
\end{equation}
and for every $n\geq 1, f\in C(K_d)$ and $x=(x_i)_{1\leq i\leq d}\in
K_d$, set
\begin{equation}\label{4.39}
\begin{split}
B_n(f)(x):&=\sum_{\begin{subarray}{l}h_1,\dots,h_d=0,\dots,n\\h_1+\dots+h_d\leq
n\end{subarray}}f\Big(\frac{h_1}{n},\dots,\frac{h_d}{n}\Big)\frac{n!}{h_1!\cdots
h_d!(n-h_1-\dots-h_d)!} \\&\kern3truecm \times x_1^{h_1}\cdots
x_d^{h_d}(1-x_1-\dots -x_d)^{n-h_1-\dots-h_d}\,.
\end{split}
\end{equation}
$B_n(f)$ is a polynomial and it is usually called the $n$-th
\dword{Bernstein polynomial on the $d$-dimensional simplex}
associated with $f$. These polynomials were first studied by Dinghas
([43]).
\begin{theorem}\sl
For every $f\in C(K_d)$,
$$\lim_{n\rightarrow\infty}B_n(f)=f\quad\text{uniformly on}\quad
K_d.$$
\end{theorem}
\begin{proof}
From the multinomial theorem, it follows that
$B_n(\textbf{1})=\textbf{1}$. Consider now the first coordinate
function $pr_1$ (see (\ref{4.2})). Then, for $x=(x_i)_{1\leq i\leq
d}\in K_d$ and $n\geq 1$,
\begin{equation*}
\begin{split}
B_n(pr_1)(x)&=\sum\limits_{\begin{subarray}{l}h_1,\dots,h_d=0,\dots,n\\h_1+\dots+h_d\leq
n\end{subarray}}
\displaystyle{\frac{h_1}{n}\frac{n!}{h_1!\cdots h_d!(n-h_1-\dots-h_d)!}}\\
&\kern3truecm\times x_1^{h_1}\cdots x_d^{h_d}(1-x_1-\dots-x_d)^{n-h_1-\dots-h_d}\\
&=x_1\sum\limits_{\begin{subarray}{l}k_1,h_2,\dots,h_d=0,\dots,n\\
k_1+h_2+\dots+h_d\leq n-1\end{subarray}}
\displaystyle{\frac{(n-1)!}{k_1!h_2!\cdots h_d!((n-1)-k_1-h_2-\dots-h_d)!}}\\
&\kern3truecm\times x_1^{k_1}x_2^{h_2}\cdots
x_d^{h_d}(1-x_1-\dots-x_d)^{(n-1)-k_1-h_2-\dots-h_d}\\ &=x_1.
\end{split}
\end{equation*}
On the other hand, for $n\geq 2, $\\
\begin{equation*}
\begin{split}
B_n(pr_1^2)(x)&=\sum\limits_{\begin{subarray}{l}h_1,\dots,h_d=0,\dots,n\\h_1+h_2+\dots+h_d\leq
n\end{subarray}}\frac{h_1^2}{n^2}\frac{n!}{h_1!\cdots
h_d!(n-h_1-\dots-h_d)!}\\&\kern1truecm\times x_1^{h_1}\cdots
x_d^{h_d}(1-x_1-\dots-x_d)^{n-h_1-\dots-h_d}\\&
=\sum\limits_{h_1+\dots +h_d\leq
n}\frac{h_1}{n}\frac{(n-1)!}{(h_1-1)!h_2!\cdots
h_d!(n-h_1-\dots-h_d)!}\\&\kern1truecm \times x_1^{h_1}\cdots
x_d^{h_d}(1-x_1-\dots-x_d)^{n-h_1-\dots-h_d}\\&
=\sum\limits_{h_1+\dots +h_d\leq
n}\frac{n-1}{n}\frac{h_1-1}{n-1}\frac{(n-1)!}{(h_1-1)!h_2!\cdots
h_d!(n-h_1-\dots-h_d)!}\\& \kern1truecm\times x_1^{h_1}\cdots
x_d^{h_d}(1-x_1-\dots-x_d)^{n-h_1-\dots-h_d}
\;\;+\;\;\frac{x_1}{n}\\&=\frac{n-1}{n}x_1^2\;\;+\;\;\frac{x_1}{n}.
\end{split}
\end{equation*}\\
Similar considerations apply to the other coordinate functions. Thus
for every $j=1,\dots,d$,
\begin{equation}\label{4.40}
B_n(pr_j)=pr_j\quad\text{and}\quad
B_n(pr_j^2)=\frac{n-1}{n}pr_j^2+\frac{1}{n}pr_j.
\end{equation}
Applying Theorem 4.2 to the sequence of positive linear operators
$(B_n)_{n\geq 1}$, we get the result.
\end{proof}

\begin{remark}\rm Theorem 4.9
gives a constructive proof of the multidimensional versions of the
Weierstrass approximation theorem, namely

\textsl{For every $f\in C(K_d)$, there exists a sequence of
algebraic polynomials that converges to $f$ uniformly on $K_d$.}
\end{remark}

Another useful generalization of the one-dimensional Bernstein
polynomial is discussed below. Consider the hypercube $Q_d=[0,1]^d$
and for every $n\geq 1, f\in C(Q_d)$ and $x=(x_i)_{1\leq i\leq d}$,
set
\begin{equation}\label{4.41}
B_n(f)(x):=\sum_{h_1,\dots,h_d=0}^nf\Big(\frac{h_1}{n},\dots,\frac{h_d}{n}\Big)\binom{n}{h_1}\cdots
\binom{n}{h_d}x_1^{h_1}(1-x_1)^{n-h_1}\cdots
x_d^{h_d}(1-x_d)^{n-h_d}\,.
\end{equation} The polynomials $B_n(f)$, $n\geq 1$,
are called the \dword{Bernstein polynomials on the hypercube}
associated with $f$. They were first studied by Hildebrandt and
Schoenberg ([70]) and Butzer ([37]). The proof of Theorem 4.9 works
also for these operators giving the same formula (\ref{4.40}).
Therefore, again by Theorem 4.2, we obtain

\begin{theorem}\sl  For every $f\in
C(Q_d)$, $$\lim_{n\rightarrow\infty}B_n(f)=f\quad\text{uniformly
on}\quad Q_d.$$
\end{theorem} We end the section by discussing some applications of
Theorem 4.1 for noncompact subsets of $\mathbb{R}^d$ .

The first application is concerned with the interval $[0,+\infty[.$
Set
\begin{equation}\label{4.42}
\begin{split} E:=\{f\in
&C([0,+\infty[)\;|\;\text{there exist}\;\alpha\geq
0\;\text{and}\;M\geq 0\;\\&\text{such that}\;|f(x)|\leq M
\exp(\alpha x)\;(x\geq 0)\},
\end{split}
\end{equation} and for every $f\in E, n\geq 1,\text{and}\ x\geq
0$, define
\begin{equation}\label{4.43}
M_n(f)(x):=\exp(-nx)\sum_{k=0}^\infty
f\Big(\frac{k}{n}\Big)\frac{n^kx^k}{k!}\,.
\end{equation} The operator $M_n$
is linear and positive and is called the $n$-th
\dword{Sz\'{a}sz-Mirakjan operator} The sequence $(M_n)_{n\geq 1}$
was first introduced and studied by Mirakjan ([86]), Favard ([55])
and Sz\'{a}sz ([113]) and is one of the most studied sequences of
positive linear operators on function spaces on $[0,+\infty[$.

A simple calculation shows that $$M_n(\textbf{1})=\textbf{1},\quad
M_n(e_1)=e_1$$ and $$M_n(e_2)=e_2+\frac{1}{n}e_1\,.$$ Therefore,
from Theorem 4.1, we get that
\begin{theorem}\sl  For every $f\in C_b([0,+\infty[)$,
$$\lim_{n\rightarrow\infty}M_n(f)=f$$ uniformly on compact subsets of
$[0,+\infty[$.
\end{theorem} The next application concerns the case
$X=\mathbb{R}^d$, $d\geq 1$. For every $n\geq 1$ and for every
real-valued Borel-measurable function
$f:\mathbb{R}^d\rightarrow\mathbb{R}$, set
\begin{equation}\label{4.44}
G_n(f)(x):=\Big(\frac{n}{4\pi}\Big)^{d/2}\int_{\mathbb{R}^d}f(t)\exp(-\frac{n}{4}\norm{t-x}^2)\dd
t
\end{equation} $(x\in\mathbb{R}^d)$.

We shall consider the operators $G_n, n\geq 1$, defined on the
lattice subspace of all Borel-measurable functions $f\in
F(\mathbb{R}^d)$ for which the integral (\ref{4.44}) is absolutely
convergent. They are referred to as the \dword{Gauss-Weierstrass
convolution operators on $\mathbb{R}^d$}. Among other things, they
are involved in the study of the heat equation on $\mathbb{R}^d$
(see, e.g., [8], [38]).

By using the following formulae:
\begin{equation}\label{4.45}
\Big(\frac{1}{2\pi\sigma^2}\Big)^{1/2}\int_{\mathbb{R}}\exp\Big(-\frac{(t-\alpha)^2}{2\sigma^2}\Big)\dd
t=1,
\end{equation}
\begin{equation}\label{4.46}
\Big(\frac{1}{2\pi\sigma^2}\Big)^{1/2}\int_{\mathbb{R}}t
\exp\Big(-\frac{(t-\alpha)^2}{2\sigma^2}\Big)\dd t=\alpha,
\end{equation} and
\begin{equation}\label{4.47}
\Big(\frac{1}{2\pi\sigma^2}\Big)^{1/2}\int_{\mathbb{R}}t^2
\exp\Big(-\frac{(t-\alpha)^2}{2\sigma^2}\Big)\dd
t\;=\;\alpha^2+\sigma^2,
\end{equation} $(\alpha\in\mathbb{R}, \sigma>0)$ (see [30, formulae
(\ref{4.14}) and (\ref{4.15})]) and appealing to Fubini's theorem,
we immediately obtain for every $i=1,\dots,d$ and
$x\in\mathbb{R}^d$,
\begin{equation}\label{4.48}
G_n(\textbf{1})(x):=\prod_{j=1}^d\Big(\frac{n}{4\pi}\Big)^{1/2}\int_{\mathbb{R}}\exp(-\frac{n}{4}(t_j-x_j)^2)\dd
t_j\;=\;1,
\end{equation}
\begin{eqnarray}\label{4.49}
G_n(pr_i)(x)&=&\prod_{\substack{j=1\\j\neq
i}}^d\Big(\frac{n}{4\pi}\Big)^{1/2}\int_{\mathbb{R}}\exp\Big(-\frac{n}{4}(t_j-x_j)^2\Big)\dd
t_j\times\\\nonumber
&\times&\Big(\frac{n}{4\pi}\Big)^{1/2}\int_{\mathbb{R}}t_i
\exp\Big(-\frac{n}{4}(t_i-x_i)\Big)^2\dd t_i\;=\;x_i\;=\;pr_i(x)
\end{eqnarray} and
\begin{eqnarray}\label{4.50} G_n(pr_i^2)(x)&=&\prod_{\substack{j=1\\j\neq
i}}^d\Big(\frac{n}{4\pi}\Big)^{1/2}\int_{\mathbb{R}}\exp\Big(-\frac{n}{4}(t_j-x_j)^2\Big)\dd
t_j\times\\\nonumber
&\times&\Big(\frac{n}{4\pi}\Big)^{1/2}\int_{\mathbb{R}}t_i^2
\exp\Big(-\frac{n}{4}(t_i-x_i)^2\Big)\dd
t_i\;=\;x_i^2+\frac{2}{n}\;=\;pr_i^2(x)+\frac{2}{n}\,.
\end{eqnarray}
From Theorem 4.1, we then obtain
\begin{theorem}\sl  For every $f\in
C_b(\mathbb{R}^d)$,
$$\lim_{n\rightarrow\infty}G_n(f)=f\quad\text{uniformly on compact
subsets of}\quad\mathbb{R}^d.$$
\end{theorem} For additional
properties of Gauss-Weierstrass operators, we refer, e.g., to [38]
(see also a recent generalization given in [21-23]).

\section{Korovkin-type theorems for positive linear operators}

After the discovery of Korovkin's theorem, several mathematicians
tried to extend it in several directions with the aim, for instance:
\begin{itemize}
\item [(i)] to find other subsets of functions satisfying the same
property as $\{\textbf{1}, e_1,e_2\}$; \item [(ii)] to establish
results like Theorem 3.1 in other function spaces or in abstract
Banach spaces;
\item [(iii)] to establish results like Theorem 3.1 for other classes of
linear operators.
\end{itemize}

As a consequence of these investigations, a new theory was created
which is nowadays called ``Korovkin-type Approximation Theory". This
theory has strong and fruitful connections not only with classical
approximation theory but also with real analysis, functional
analysis, harmonic analysis, probability theory and partial
differential equations. We refer to [8] for a rather detailed
description of the development of this theory.

We shall next discuss some of the main results of the theory
obtained in the framework of the spaces $C_0(X)$ ($X$ locally
compact noncompact space), $C(X)$ ($X$ compact space), occasionally
in $L^p(X,\widetilde{\mu})$ spaces, $1\leq p<+\infty$, and in
weighted continuous function spaces. These spaces play a central
role in the theory and they are the most useful in applications.

In addition, it transpires that the elementary methods used in the
previous sections are not appropriate to give a wider possibility to
determine other test functions (like $\{\textbf{1},e_1,e_2\}$ in
$C([0,1]))$, let alone to characterize them. To this aim, it seems
to be unavoidable to use a more precise analysis involving the point
topology of the underlying space and the positive linear functionals
on the relevant continuous function spaces.

It turns out that the right framework is the local compactness which
we choose to consider in its generality because the class of locally
compact spaces includes, other than the Euclidean spaces
$\mathbb{R}^d$, $d\geq 1$, and their open or closed subsets, many
other topological spaces which are important in their own right. The
reader who is not interested in this level of generality, may
replace everywhere our locally compact spaces with a space
$\mathbb{R}^d$, $d\geq 1$, or with an open or a closed subset of it,
or with the intersection of an open subset and a closed subset of
$\mathbb{R}^d$.  However, this restriction does not produce any
simplification in the proofs or in the methods.

In the sequel, given a locally compact Hausdorff space $X$, we shall
denote by $$K(X)$$ the linear subspace of all real-valued continuous
functions on $X$ having compact support. Then $K(X)\subset C_b(X)$.
We shall denote by
$$C_0(X)$$ the closure of $K(X)$ with respect to the sup-norm
$\norm{\cdot}_\infty$ (see (\ref{2.3})). Thus, $C_0(X)$ is a closed
subspace of $C_b(X)$ and hence, endowed with the norm
$\norm{\cdot}_\infty$, is a Banach space.

If $X$ is compact, then $C_0(X)=C(X)$.  If $X$ is not compact, then
a function $f\in C(X)$ belongs to $C_0(X)$ if and only if
\begin{gather*}
\textsl{ for every } \varepsilon>0 \textsl{ there exists a compact
subset } K \textsl{ of } X\\ \textsl{ such that }
|f(x)|\leq\varepsilon \textsl{ for every } x\in X \backslash K.
\end{gather*} For additional topological and
analytical properties of locally compact spaces and of some relevant
continuous function spaces on them, we refer to the \Appendix.

The spaces $C_0(X)$ and $C(X)$, ($X$ compact), endowed with the
natural pointwise ordering and the sup-norm, become Banach lattices.
Similarly $L^p(X,\widetilde{\mu})$, endowed with the natural norm
$\norm{\cdot}_p$ and the ordering $$f\leq g\quad \text{if}\quad
f(x)\leq g(x)\quad\text{for}\quad \widetilde{\mu}{\text{-a.e.}}\;
x\in X,$$ is a Banach lattice (for more details on
($L^p(X,\widetilde{\mu})$,$\norm{\cdot}_p$) spaces, we refer to the
\Appendix\ (formulae (\ref{11.11}) and (\ref{11.12}))).

For the reader's convenience, we recall that a \dword{Banach
lattice} $E$ is a vector space endowed with a norm $\norm{\cdot}$
and an ordering $\leq$ on $E$ such that
\begin{itemize}
\item [(i)] $(E,\norm{\cdot})$ is a Banach space; \item [(ii)] $(E,\leq)$ is
a vector lattice; \item [(iii)] If $f,g\in E$ and $|f|\leq|g|$ then
$\norm{f}\leq\norm{g}$.
\end{itemize}
(where $|f| :=\sup (f,-f)$ for every $f\in E$).

Actually it is convenient to state the main definitions of the
theory in the framework of Banach lattices. However, the reader not
accustomed to this terminology may replace our abstract spaces with
the concrete ones such as $C_0(X)$, $C(X)$ or
$L^p(X,\widetilde{\mu})$. For more details on Banach lattices, we
refer, e.g., to [2].

If $E$ and $F$ are Banach lattices, a linear operator
$L:E\longrightarrow F$ is said to be \dword{positive} if $$L(f)\geq
0\quad\text{for every}\quad f\in E,\; f\geq 0.$$ Every positive
linear operator $L:E\longrightarrow F$ is continuous, (see, e.g.,
[2, Theorem 12.3]). Moreover, if $E=C(X)$, $X$ compact, then
$\norm{L}=\norm{L(\textbf{1})}$.

A \dword{lattice homomorphism} $S:E\longrightarrow F$ is a linear
operator satisfying $|S(f)|=S(|f|)$ for every $f\in E$.
Equivalently, this means that $S$ preserves the finite lattice
operation, i.e., for every $f_1,\dots,f_n\in E,\; n\geq 2$,\;
\begin{center}
$S\left(\underset{1\leq i\leq n}{\inf} f_i\right)=\underset{1\leq i\leq
n}{\inf}S(f_i)$ and $S\left(\underset{1\leq i\leq
n}{\sup}f_i\right)=\underset{1\leq
i\leq n}{\sup}S(f_i)$.\\
\end{center}
For instance, if $X$ is a locally
compact Hausdorff space, $\widetilde{\mu}$ a regular finite Borel
measure on $X$ and $1\leq p<+\infty$, then the natural embedding
$J_p:C_0(X)\longrightarrow L^p(X,\widetilde{\mu})$ defined by
$J_p(f):=f \ (f\in C_0(X))$ is a lattice homomorphism.

Analogously, if $X$ and $Y$ are compact spaces and
$\varphi:Y\longrightarrow X$ is a continuous mapping, then the
composition operator $T_{\varphi}(f):=f\circ\varphi,(f\in C(X)),$ is
a lattice homomorphism from $C(X)$ into $C(Y)$. Every lattice
homomorphism is positive and hence continuous. A linear bijection
$S:E\longrightarrow F$ is a lattice homomorphism if and only if $S$
and its inverse $S^{-1}$ are both positive [2, Theorem 7.3]. In this
case, we also say that $S$ is a \dword{lattice isomorphism}. When
there exists a lattice isomorphism between $E$ and $F$, then we say
that $E$ and $F$ are \dword{lattice isomorphic}.

The following definition, which is one of the most important of the
theory, is clearly motivated by the Korovkin theorem and was first
formulated by V.~A.~Baskakov ([28]).

\begin{definition}\sl
A subset $M$ of a Banach lattice $E$ is said to be a \dword{Korovkin
subset} of $E$ if for every sequence $(L_n)_{n\geq 1}$ of positive
linear operators from $E$ into $E$ satisfying
\\ \indent (i) \ \  $\underset{n\geq 1}{\sup}\norm{L_n}<+\infty$,\\ and\\
\indent (ii) \ \  $\lim\limits_{n\rightarrow\infty}L_n(g)=g$ for
every $g\in M$,\\ it turns out that
$$\lim\limits_{n\rightarrow\infty}L_n(f)=f\quad\text{for every}\quad f\in
E.$$
\end{definition}

Note that, if $E=C(X)$, $X$ compact space, and the constant function
$\textbf{1}$ belongs to the linear subspace $\mathcal{L}(M)$
generated by $M$, then condition (i) is superfluous because it is a
consequence of (ii).

According to Definition 5.1, we may restate Korovkin's Theorem 3.1
by saying that $\{\textbf{1},e_1,e_2\}$ is a Korovkin set in
$C([0,1])$.

We also point out that a subset $M$ is a Korovkin subset of $E$ if
and only if the linear subspace $\mathcal{L}(M)$ generated by $M$ is
a Korovkin subset.  In the sequel, a linear subspace that is a
Korovkin subset will be referred to as a \dword{Korovkin
subspace} of $E$. If $E$ and $F$ are lattice isomorphic and if
$S:E\longrightarrow F$ is a lattice isomorphism, then a subset $M$
of $E$ is a Korovkin subset in $E$ if and only if $S(M)$ is a
Korovkin subset in $F$.

Korovkin sets (when they exist) are useful for investigating the
convergence of equibounded sequences of positive linear operators
towards the identity operator or, from the point of view of
approximation theory, the approximation of every element $f\in E$ by
means of $(L_n(f))_{n\geq 1}$.

According to Lorentz ([80]), who first proposed a possible
generalization, it seems to be equally interesting to study the
following more general concept.
\begin{definition}\sl  Let $E$ and $F$ be Banach lattices and consider a
positive linear operator $T:E\longrightarrow F$.  A subset $M$ of
$E$ is said to be a \dword{Korovkin subset} of $E$ for $T$ if for
every sequence $(L_n)_{n\geq 1}$ of positive linear operators from
$E$ into $F$ satisfying
\smallskip\\
\indent (i) \ \ $\underset{n\geq 1}{\sup}\norm{L_n}<+\infty$ \\ and\\
\indent (ii) \ \ $\lim\limits_{n\rightarrow\infty}L_n(g)=T(g)$ for
every $g\in M$, \\ it turns out that
$$\lim\limits_{n\rightarrow\infty}L_n(f)=T(f)\quad\text{for every}\quad
f\in E.$$
\end{definition} Thus, such subsets can be used to investigate
the convergence of equibounded sequences of positive linear
operators towards a given positive linear operator
$T:E\longrightarrow F$ or to approximate weakly $T$ by means of
(generally, simpler) linear operators $L_n$, $n\geq 1$.

In the light of the above definition, two problems arise quite
naturally:

\begin{problem} Given a positive linear operator $T:E\longrightarrow F$,
find conditions under which there exists a nontrivial (i.e., the
linear subspace generated by it is not dense) Korovkin subset for
$T$.  In this case, try to determine some or all of them.
\end{problem}
\begin{problem} Given a subset $M$ of $E$, try to determine some or all of the
positive linear operators $T:E\rightarrow F$ (if they exist) for
which $M$ is a Korovkin subset.
\end{problem}

In the next sections, we shall discuss some aspects related to
Problem 5.3 (for further details, we refer to [8, Sections 3.3 and
3.4]). As regards Problem 5.4, very few results are available (see,
e.g., [71], [72-73], [74], [114-116]).

The next result furnishes a complete characterization of Korovkin
subsets for positive linear operators in the setting of $C_0(X)$
spaces.  It was obtained by Yu.~A.~Shashkin ([110]) in the case when
$X=Y$, $X$ compact metric space and $T=I$ the identity operator, by
H.~Berens and G.~G.~Lorentz ([33]) when $X=Y$, $X$ topological
compact space, $T=I$, by H.~Bauer and K.~Donner ([31]) when $X=Y$,
$X$ locally compact space, $T=I$, by C.~A.~Micchelli ([85]) and
M.~D.~Rusk ([105]) when $X=Y$, $X$ compact, and by F.~Altomare ([4])
in the general form below.

We recall that $\mathcal{M}_b^+(X)$ denotes the cone of all bounded
Radon measures on $X$ (see \Appendix).
\begin{theorem}\sl  Let $X$ and $Y$ be
locally compact Hausdorff spaces. Further, assume that $X$ has a
countable base and $Y$ is metrizable.  Given a positive linear
operator $T:C_0(X)\longrightarrow C_0(Y)$ and a subset $M$ of
$C_0(X)$, the following statements are equivalent:
\begin{enumerate}
\item [(i)] $M$ is a Korovkin subset of $C_0(X)$ for $T$. \item [(ii)] If
$\mu\in \mathcal{M}_b^+(X)$ and $y\in Y$ satisfy $\mu(g)=T(g)(y)$
for every $g\in M$, then $\mu(f)=T(f)(y)$ for every $f\in C_0(X)$.
\end{enumerate}
\end{theorem}
\begin{proof} (i)$\Rightarrow$(ii).  Fix
$\mu\in\mathcal{M}_b^+(X)$ and $y\in Y$ satisfying $\mu(g)=T(g)(y)$
for every $g\in M$.  Consider a decreasing countable base
$(U_n)_{n\geq 1}$ of open neighborhoods of $y$ in $Y$ and, for every
$n\geq 1$, choose $\varphi_n\in K(Y)$ such that: $0\leq\varphi_n\leq
1, \varphi_n(y)=1$ and $\supp(\varphi_n)\subset U_n$ (see Theorem
11.1 of the \Appendix).

Accordingly, define $L_n:C_0(X)\longrightarrow C_0(Y)$ by
$$L_n(f):=\mu(f)\varphi_n+T(f)(1-\varphi_n)\quad (f\in C_0(X)).$$ Each
$L_n$ is linear, positive and $\norm{L_n}\leq\norm{\mu}+\norm{T}$.
Moreover, if $g\in M$, then
$$\lim_{n\rightarrow\infty}L_n(g)=T(g)\quad\text{in}\quad C_0(Y),$$
because, given $\varepsilon>0$, there exists $v\in\mathbb{N}$ such
that $$|T(g)(z)-T(g)(y)|\leq\varepsilon\quad\text{for every}\quad
z\in U_v.$$ Hence, since for every $n\geq v$ (thus $U_n\subset U_v)$
and for every $z\in Y$
$$|L_n(g)(z)-T(g)(z)|=\varphi_n(z)|T(g)(z)-T(g)(y)|,$$ we
get $$|L_n(g)(z)-T(g)(z)|=\begin{cases} 0&\text{if}\quad z\not\in U_n,\\
\leq\varepsilon&\text{if}\quad z\in U_n,
\end{cases}$$ and so
$\norm{L_n(g)-T(g)}\leq\varepsilon$.

Since $M$ is a Korovkin subset for $T$, it turns out that, for every
$f\in C_0(X), \lim\limits_{n\rightarrow\infty}L_n(f)=T(f)$ and hence
$\lim\limits_{n\rightarrow\infty}L_n(f)(y)=T(f)(y)$. But, for every
$n\geq 1$, $L_n(f)(y)=\mu(f)$ and this completes the proof of (ii).

(ii)$\Rightarrow$(i).  Our proof starts with the observation that
from statement (ii), it follows that
\begin{equation}\label{5.1}
\text{if}\quad \mu\in\mathcal{M}_b^+(X)\quad\text{and}\quad
\mu(g)=0\quad\text{for every}\quad g\in M, \quad\text{then}\quad
\mu=0.
\end{equation} Moreover, since $X$ has a countable base, every bounded
sequence in $\mathcal{M}_b^+(X)$ has a vaguely convergent
subsequence (see Theorem 11.8 of \Appendix).  Consider now a
sequence $(L_n)_{n\geq 1}$ of positive linear operators from
$C_0(X)$ into $C_0(Y)$ satisfying properties (i) and (ii) of
Definition 5.2 and assume that for some $f_0\in C_0(X)$ the sequence
$(L_n(f_0))_{n\geq 1}$ does not converge uniformly to $T(f_0)$.

Therefore, there exist $\varepsilon_0>0$, a subsequence
$(L_{k(n)})_{n\geq 1}$ of $(L_n)_{n\geq 1}$ and a sequence
$(y_n)_{n\geq 1}$ in $Y$ such that
\begin{equation}\label{5.2}
|L_{k(n)}(f_0)(y_n)-T(f_0)(y_n)|\geq\varepsilon_0\quad\text{for
every}\quad n\geq 1.
\end{equation} We discuss separately the two cases when
$(y_n)_{n\geq 1}$ is converging to the point at infinity of $Y$ or
not (see the \Appendix).

In the first case (which can only occur when $Y$ is noncompact), we
have $\lim\limits_{n\rightarrow\infty}h(y_n)=0$ for every $h\in
C_0(Y)$.

For every $n\geq 1$, define $\mu_n\in\mathcal{M}_b^+(X)$ by
$$\mu_n(f):=L_{k(n)}(f)(y_n)\quad (f\in C_0(X)).$$ Since
$\norm{\mu_n}\leq\norm{L_{k(n)}}\leq M:=\underset{n\geq
1}{\sup}\norm{L_n}$, replacing, if necessary, the sequence
$(\mu_n)_{n\geq 1}$ with a suitable subsequence, we may assume that
there exists $\mu\in\mathcal{M}_b^+(X)$ such that
$\mu_n\longrightarrow\mu$ vaguely.  But, if $g\in M$, then
\begin{eqnarray*}
|\mu_n(g)|&\leq&|L_{k(n)}(g)(y_n)-T(g)(y_n)|+|T(g)(y_n)|\leq\\
&\leq&\norm{L_{k(n)}(g)-T(g)}+|T(g)(y_n)|
\end{eqnarray*} so that
$$\mu(g)=\lim_{n\rightarrow\infty}\mu_n(g)=0.$$
From (\ref{5.1}) it turns out that $\mu(f_0)=0$ as well and hence
$$|L_{k(n)}(f_0)(y_n)-T(f_0)(y_n)|=|\mu_n(f_0)-T(f_0)(y_n)|\longrightarrow
0$$ which contradicts (\ref{5.2}).

Consider now the case where the sequence $(y_n)_{n\geq 1}$ does not
converge to the point at infinity of $Y$.  Then, by replacing it
with a suitable subsequence, we may assume that it converges to some
$y\in Y$.

Again, consider for every $n\geq 1$
$$\mu_n(f):=L_{k(n)}(f)(y_n)\quad (f\in C_0(X)).$$ As in the
previous reasoning, we may assume that there exists
$\mu\in\mathcal{M}_b^+(X)$ such that $\mu_n\longrightarrow\mu$
vaguely.  Then, for $g\in M$, since
$$|\mu_n(g)-T(g)(y_n)|\leq\norm{L_{k(n)}(g)-T(g)}\longrightarrow 0,$$ we get
$\mu(g)=T(g)(y)$.  Therefore, assumption (ii)\;implies
$\mu(f_0)=T(f_0)(y)$, or equivalently,
$$\lim_{n\rightarrow\infty}\Bigl(L_{k(n)}(f_0)(y_n)-T(f_0)(y_n)\Bigr)=0$$ which is
impossible because of (\ref{5.2}).
\end{proof}
Some applications of Theorem 4.5 will be shown in subsequent
sections.
\section{Korovkin-type
theorems for the identity operator in $C_0(X)$} In this section, we
discuss more closely those subsets of $C_0(X)$ that are Korovkin
subsets in $C_0(X)$ (see Definition 4.1), i.e., that are Korovkin
subsets for the identity operator on $C_0(X)$.

Throughout the whole section, we shall fix a locally compact
Hausdorff space with a countable base, which is then metrizable as
well. The next result immediately follows from Theorem 4.5.

\begin{theorem}\sl  ([31]).  Given a subset $M$ of $C_0(X)$, the following
statements are equivalent:
\begin{itemize}
\item [(i)] $M$ is a Korovkin subset of $C_0(X)$. \item [(ii)] If
$\mu\in\mathcal{M}_b^+(X)$ and $x\in X$ satisfy $\mu(g)=g(x)$ for
every $g\in M$, then $\mu(f)=f(x)$ for every $f\in C_0(X)$, i.e.,
$\mu=\delta_x$ (see (2.11)).
\end{itemize} \end{theorem} In order to
discuss a first application of Theorem 6.1, we recall that a mapping
$\varphi:Y\longrightarrow X$ between two locally compact Hausdorff
spaces $Y$ and $X$ is said to be \dword{proper} if for every compact
subset $K\in X$, the inverse image $\varphi^{-1}(K):=\{y\in Y\mid
\varphi(y)\in K\}$ is compact in $Y$. In this case,
$f\circ\varphi\in C_0(Y)$ for every $f\in C_0(X)$.

\begin{corollary}\sl  Let $Y$ be a metrizable locally compact Hausdorff space.
If $M$ is a Korovkin subset of $C_0(X)$, then $M$ is a Korovkin
subset for any positive linear operator $T:C_0(X)\longrightarrow
C_0(Y)$ of the form
$$T(f):=\lambda(f\circ\varphi)\quad (f\in C_0(X))$$ where $\lambda\in
C_b(Y), \lambda\geq 0$, and $\varphi:Y\longrightarrow X$ is a proper
mapping.
\end{corollary}
\begin{proof} According to Theorem 5.5, we have
to show that, if $\mu\in\mathcal{M}_b^+(X)$ and $y\in Y$ satisfy
$\mu(g)=\lambda(y)g(\varphi(y))$ for every $g\in M$, then
$\mu(f)=\lambda(y)f(\varphi(y))$ for every $f\in C_0(X)$.

If $\lambda(y)=0$, then $\mu=0$ on $M$ and hence $\mu=0$ by Theorem
6.1 and property (1) in the proof $(ii)\Longrightarrow(i)$ of
Theorem 5.5. If $\lambda(y)>0$, it suffices to apply Theorem 6.1 to
$\displaystyle{\frac{1}{\lambda(y)}}\mu$ and $\varphi(y)$.
\end{proof}

Note that, if $Y$ is a closed subset of $X$, then the canonical
mapping $\varphi:Y\longrightarrow X$ defined by
$\varphi(y):=y\;(y\in Y)$ is proper.  Therefore, from Corollary 5.2,
we get:

\begin{corollary}\sl  Let $M$ be a
Korovkin subset of $C_0(X)$.  Consider an equibounded sequence
$(L_n)_{n\geq 1}$ of positive linear operators from $C_0(X)$ into
$C_0(X)$. Then the following properties hold. \begin{itemize}
\item [1)] Given a closed subset $Y$ of $X$ and $\lambda\in C_b(Y)$,
$\lambda\geq 0$, if $\lim\limits_{n\rightarrow\infty}L_n(g)=\lambda
g$ uniformly on $Y$ for every $g\in M$, then
$$\lim_{n\rightarrow\infty}L_n(f)=\lambda f\quad\text{uniformly on}\quad
Y\quad\text{for every}\quad f\in C_0(X).$$ \item [2)] If
$\lim\limits_{n\rightarrow\infty}L_n(g)=g$ uniformly on compact
subsets of $X$ for every $g\in M$, then
$\lim\limits_{n\rightarrow\infty}L_n(f)=f$ uniformly on compact
subsets of $X$ for every $f\in C_0(X)$.
\end{itemize} \end{corollary} Next, we
proceed to investigate some useful criteria to explicitly determine
Korovkin subsets of $C_0(X)$. We begin with the following result
which is at the root of all subsequent results. We recall that, if
$M$ is a subset of $C_0(M)$, then the symbol $\mathcal{L}(M)$
denotes the linear subspace of $C_0(X)$ generated by $M$.

\begin{proposition}\sl  Let $M$ be a subset of $C_0(X)$ and assume that for
every $x,y\in X, x\neq y$, there exists $h\in\mathcal{L}(M), h\geq
0$, such that $h(x)=0$ and $h(y)>0$.  Then $M$ is a Korovkin subset
of $C_0(X)$.
\end{proposition}
\begin{proof} We shall verify condition (ii) of Theorem
6.1. Therefore, consider $\mu\in\mathcal{M}_b^+(X)$ and $x\in X$
satisfying $\mu(g)=g(x)$ for every $g\in M$ and, hence, for every
$g\in\mathcal{L}(M)$.

If $y\in X, y\neq x$, then there exists $h\in\mathcal{L}(M)$, $h\geq
0$, such that $h(y)>0$ and $0=h(x)=\mu(h)$.  Therefore, by Theorem
11.7 of the \Appendix, there exists $\alpha\geq 0$ such that
$\mu=\alpha\delta_x$.  Choosing $h\in\mathcal{L}(M)$ such that
$h(x)>0$, we get $\alpha h(x)=\mu(h)=h(x)$, so $\alpha=1$ and
$\mu=\delta_x$.
\end{proof} Below, we state an important
consequence of Proposition 6.4.  In the sequel, if $M$ is a subset
of $C(X)$ and for $f_0\in C(X)$, we shall set
\begin{equation}\label{6.1}
f_0M:=\{f_0\cdot f\mid f\in M\}
\end{equation} and
\begin{equation}\label{6.2}
f_0M^2:=\{f_0\cdot f^2\mid f\in M\}.
\end{equation}
\begin{theorem}\sl
Consider a strictly positive function $f_0\in C_0(X)$ and a subset
$M$ of $C(X)$ that separates the points of $X$. Furthermore, assume
that $f_0M\cup f_0M^2\subset C_0(X)$.  Then $\{f_0\}\cup f_0M\cup
f_0M^2$ is a Korovkin subset of $C_0(X)$.

If, in addition, $M$ is finite, say $M=\{f_1,\dots,f_n\}, n\geq 1$,
then
$$\left\{f_0,f_0f_1,\dots,f_0f_n, f_0\sum\limits_{i=1}^nf_i^2\right\}$$ is a Korovkin
subset of $C_0(X)$.
\end{theorem}
\begin{proof} If $x,y\in X, x\neq y$, then
there exists $f\in M$ such that $f(x)\neq f(y)$.  Therefore, the
function $h:=f_0(f-f(x))^2$ belongs to $\mathcal{L}(\{f_0\}\cup
f_0M\cup f_0M^2)$, it is positive and $h(x)=0<h(y)$.  Thus the
result follows from Proposition 6.4.

If $M=\{f_1,\dots,f_n\}$ is finite, then in the above reasoning one
can consider the function $h:=f_0\sum\limits_{i=1}^n(f_i-f_i(x))^2$.
\end{proof}

The next result is an obvious consequence of Theorem 6.5 but is
worth being stated explicitly because of its connection with the
Stone-Weierstrass theorem (see Section 9).

\begin{theorem}\sl  Consider a
strictly positive function $f_0\in C_0(X)$ and a subset $M$ of
$C_0(X)$ that separates the points of $X$.  Then $\{f_0\}\cup
f_0M\cup f_0M^2$ is a Korovkin subset of $C_0(X)$.

Moreover, if $M$ is finite, say $M=\{f_1,\dots,f_n\}, n\geq1$, then
$$\left\{f_0,f_0f_1,\dots,f_0f_n,f_0\sum\limits_{i=1}^nf_i^2\right\}$$ is a Korovkin
subset of $C_0(X)$.

Finally, if $f_0$ is also injective, then $\{f_0,f_0^2,f_0^3\}$ is a
Korovkin subset of $C_0(X)$.
\end{theorem} From Theorem 6.6, the following
result immediately follows.
\begin{corollary}\sl  \begin{itemize}
\item [(1)] $\{e_1,e_2,e_3\}$ is a Korovkin subset of $C_0(]0,1]).$ \item
[(2)] $\{e_{-1}, e_{-2}, e_{-3}\}$ is a Korovkin subset of
$C_0([1,+\infty[)$, where $e_{-k}(x):=x^{-k}$ for every $\ x\in
\lbrack 1,+\infty \lbrack $ and $k=1,2,3$.  \item [(3)] $\{f_{1},
f_{2}, f_{3}\}$ is a Korovkin subset of $C_0([0,+\infty[)$ where
$f_{k}(x):=\exp{(-kx)}$ for every $\ x\in \lbrack 0,+\infty \lbrack
$ and $k=1,2,3$.  \item [(4)]
$\{\Phi,pr_1\Phi,\dots,pr_d\Phi,\norm{\cdot}^2\Phi\}$ is a Korovkin
subset of $C_0(\mathbb{R}^d)$, $d\geq 1$, where
$\Phi(x):=exp(-\norm{x}^2)$ for every $\ x\in\mathbb{R}^d$.
\end{itemize}
\end{corollary}

A useful generalization of the previous result is presented
below.

\begin{proposition}\sl  Given
$\lambda_1,\lambda_2,\lambda_3\in\mathbb{R}$,
$0<\lambda_1<\lambda_2<\lambda_3$, then
\begin{itemize}
\item [(1)] $\{e_{\lambda _{1}},e_{\lambda _{2}},e_{\lambda _{3}}\}$ is a Korovkin
subset of $C_0(]0,1])$ where $e_{\lambda _{k}}(x):=x^{\lambda _{k}}$
for every $\ x\in \rbrack 0,1\rbrack $ and $k=1,2,3$.  \item [(2)]
$\{e_{-\lambda _{1}},e_{-\lambda _{2}},e_{-\lambda _{3}}\}$ is a
Korovkin subset of $C_0([1,+\infty[)$ where $e_{-\lambda
_{k}}(x):=x^{-\lambda _{k}}$ for every $\ x\in \lbrack 1,+\infty
\lbrack $ and $k=1,2,3$.  \item [(3)] $\{f_{\lambda _{1}},
f_{\lambda _{2}}, f_{\lambda _{3}}\}$ is a Korovkin subset of
$C_0([0,+\infty[)$ where $f_{\lambda _{k}}(x):=\exp (-\lambda
_{k}x)$ for every $\ x\in \lbrack 0,+\infty \lbrack $ and $k=1,2,3$.
\end{itemize}
\end{proposition}

\begin{proof} We give the proof only for
(3); the proofs of the other statements are left to the reader (see,
e.g., [8, Proposition 4.2.4]).  We shall apply Proposition 6.4 and
to this end fix $x_0\in[0,+\infty[$.  Then, by using
differential calculus, it is not difficult to show that the function
$$h(x):=\exp(-\lambda_1 x)+ \alpha \exp(-\lambda_2 x)+\beta
\exp(-\lambda_3 x)\quad (x\geq 0)$$ where
$\alpha:=\displaystyle{\frac{\lambda_1-\lambda_3}{\lambda_3-\lambda_2}}\exp((\lambda_2-\lambda_1)x_0)$
and
$\beta:=\displaystyle{\frac{\lambda_2-\lambda_1}{\lambda_3-\lambda_2}}\exp((\lambda_3-\lambda_1)x_0)$,
satisfies $h(x_0)=0<h(y)$ for every $y\in[0,+\infty[$, $y\neq x_0$.
\end{proof}

Here, we discuss some applications of Proposition 6.8.  The first
application is taken from [31, Proposition 4.1].

We begin by recalling that, if $\varphi\in L^1([0,+\infty[)$, then
the \dword{Laplace transform} of $\varphi$ on $[0,+\infty[$ is
defined by
\begin{equation}\label{6.3}
\mathcal{L}(\varphi)(\lambda):=\int_0^{+\infty}\ee^{-\lambda
t}\varphi(t)\dd t\quad (\lambda\geq 0).
\end{equation} By Lebesgue's dominated
convergence theorem, $\mathcal{L}(\varphi)\in C_0([0,+\infty[)$ and
$\norm{\mathcal{L}(\varphi)}\leq\norm{\varphi}_1$. By means of
$\varphi$, we may naturally define a linear operator
$L_{\varphi}:C_0([0,+\infty[)\rightarrow C_0([0,+\infty[)$ by
setting for every $f\in C_0([0,+\infty[)$ and $x\geq 0$
\begin{equation}\label{6.4}
L_{\varphi}(f)(x):=\int_0^{+\infty}f(x+y)\varphi(y)\dd
y=\int_x^{+\infty}f(u)\varphi (x-u)\dd u.
\end{equation} Note that $L_{\varphi}(f)\in C_0([0,+\infty[)$ by
virtue of Lebesgue's dominated convergence theorem.  Moreover
$L_{\varphi}$ is bounded and
$\norm{L_{\varphi}}\leq\norm{\varphi}_1$.  Finally, if $\varphi\geq
0$, then $L_{\varphi}$ is positive as well.

For every $\lambda>0$, denoting by $f_{\lambda}$ the function
\begin{equation}\label{6.5} f_{\lambda}(x):=\exp(-\lambda x)\quad (x\geq 0),
\end{equation} we get
\begin{equation}\label{6.6}
L_{\varphi}(f_{\lambda})=\mathcal{L}(\varphi)(\lambda)f_{\lambda}.
\end{equation}
\begin{proposition}\sl  Consider a sequence $(\varphi_n)_{n\geq
1}$ of positive functions in $L^1([0,+\infty[)$ such that
$\underset{n\geq 1}{\sup}\norm{\varphi_n}_1<+\infty$ and, for every
$n\geq 1$, denote by
$$L_n:C_0([0,+\infty[)\rightarrow C_0([0,+\infty[)$$ the positive linear
operator associated with $\varphi_n$ defined by (\ref{6.4}). Then
the following statements are equivalent:
\begin{itemize}
\item [(i)] $\lim\limits_{n\rightarrow\infty}L_n(f)=f$ uniformly on
$[0,+\infty[$ for every $f\in C_0([0,+\infty[)$.  \item [(ii)]
$\lim\limits_{n\rightarrow\infty}\mathcal{L}(\varphi_n)(\lambda)=1$ for
every $\lambda>0$.  \item [(iii)] There exist
$\lambda_1,\lambda_2,\lambda_3\in\mathbb{R},\;0<\lambda_1<\lambda_2<\lambda_3$
such that
$$\lim_{n\rightarrow\infty}\mathcal{L}(\varphi_n)(\lambda_k)=1\quad\text{for}\quad k=1,2,3.$$
\end{itemize} \end{proposition}
\begin{proof} The
implication (i)$\Rightarrow$(ii) follows from (\ref{6.6}). The
implication (ii)$\Rightarrow$(iii) being obvious, the only point
remaining concerns (iii)$\Rightarrow$(i).

Since $(\varphi_n)_{n\geq 1}$ is bounded in $L^1([0,+\infty[)$, we
have $\underset{n\geq 1}{\sup}\norm{L_n}<+\infty$.
Moreover, (iii) means that
$\lim\limits_{n\rightarrow\infty}L_n(f_{\lambda_k})=f_{\lambda_k}$
uniformly on $[0,+\infty[)$ by (\ref{6.6}).  Therefore, (i) follows
by applying part (3) of Proposition 6.8.
\end{proof}

As another application of the previous results, we shall study the
behaviour of the Sz\'{a}sz-Mirakjan operators (see (\ref{4.43})) on
$C_0([0,+\infty[)$ and on continuous function spaces on
$[0,+\infty[$ with polynomial weights.

We recall that these operators are defined by
\begin{equation}\label{6.7}
M_n(f)(x):=\exp(-nx)\sum_{k=0}^\infty
f\Big(\frac{k}{n}\Big)\frac{n^kx^k}{k!}
\end{equation} for $n\geq 1, x\geq
0$ and $f\in C([0,+\infty[)$ such that $|f(x)|\leq M\exp(\alpha x)\;
(x\geq 0)$ for some $M\geq 0$ and $\alpha>0$.
\begin{lemma}\sl  If $f\in
C_0([0,+\infty[)$, then $M_n(f)\in C_0([0,+\infty[)$ and
$\norm{M_n(f)}\leq\norm{f}$ for every $n\geq 1$.
\end{lemma}
\begin{proof}
The function $M_n(f)$ is continuous because the series (\ref{6.7})
is uniformly convergent on compact subsets of $[0,+\infty[$.
Moreover, for every $x\geq 0$
$$|M_n(f)(x)|\leq\norm{f}\exp(-nx)\sum_{k=0}^\infty\frac{n^kx^k}{k!}=\norm{f}.$$
In order to show that $M_n(f)\in C_0([0,+\infty[)$, given
$\varepsilon>0$, choose $v\in \mathbb{N}$ such that
$|f\Big(\frac{k}{n}\Big)|\leq\varepsilon$ for every $k\geq v$. For
sufficiently large $x\geq 0$ we get
$$\exp(-nx)\sum_{k=0}^v\Big|f\Big(\frac{k}{n}\Big)\Big|\frac{(nx)^k}{k!}\leq\varepsilon$$
so that
\begin{eqnarray*} |M_n(f)(x)\Big|&\leq&
\exp(-nx)\sum\limits_{k=0}^v\Big|f\Big(\frac{k}{n}\Big)\Big|\frac{(nx)^k}{k!}+
\exp (-nx)\sum_{k=v+1}^\infty\Big|f
\Big(\frac{k}{n}\Big)\Big|\frac{(nx)^k}{k!}\\ &\leq& 2\varepsilon.
\end{eqnarray*}
\vskip-\baselineskip
\end{proof}
\begin{theorem}\sl  For every $f\in
C_0([0,+\infty[)$,
$$\lim_{n\rightarrow\infty}M_n(f)=f\quad\text{uniformly on}\quad
[0,+\infty[.$$
\end{theorem}
\begin{proof} Since $(M_n)_{n\geq 1}$
is an equibounded sequence of positive linear operators from
$C_0([0,+\infty[)$ into $C_0([0,+\infty[)$, by Proposition 6.8, (3),
it is sufficient to show that, for every $\lambda\geq 0$,
$$\lim_{n\rightarrow\infty}M_n(f_\lambda)=f_\lambda \quad\text{uniformly
on}\quad  [0,+\infty[$$ where $f_\lambda$ is defined by (\ref{6.5}).

A simple calculation indeed shows that
$$M_n(f_{\lambda})(x)=\exp\Big[-\lambda
x\Big(\frac{1-\exp(-\lambda/n)}{\lambda/n}\Big)\Big]$$ for every
$x\geq 0$. Since
$\lim\limits_{n\rightarrow\infty}\frac{1-\exp(-\lambda/n)}{\lambda/n}=1$,
the sequence $(M_n(f_{\lambda}))_{n\geq 1}$ converges pointwise to
$f_{\lambda}$ and is decreasing.  Moreover, each $M_n(f_{\lambda})$
and $f_{\lambda}$ vanishes at $+\infty$ and so, by Dini's theorem
applied in the framework of the compactification $[0,+\infty]$, we
obtain the uniform convergence as well.
\end{proof}

Our next aim is to discuss the behaviour of the operators $M_n$ on
continuous function spaces with polynomial weights. This naturally
leads to investigate some Korovkin-type results in weighted
continuous function spaces. For further results in this respect, we
also refer to the next Section 8 and to [10-11], [14], [60-62],
[104].

Consider again an arbitrary locally compact Hausdorff space $X$
having a countable base.
Let $w$ be a continuous weight on $X$, i.e., $w\in C(X)$ and
$w(x)>0$ for every $x\in X$, and set
\begin{equation}\label{6.8} C_0^w(X):=\{f\in
C(X)\;|\;wf\in C_0(X)\}.
\end{equation} The space $C_0^w(X)$, endowed with
the natural (pointwise) order and the weighted norm
\begin{equation}\label{6.9} \norm{f}_w:=\norm{wf}_\infty\quad (f\in C_0^w(X))
\end{equation} is a Banach lattice

If $w\in C_b(X)$, then $C_0(X)\subset C_0^w(X)$ and, if $w\in
C_0(X)$, then $C_b(X)\subset C_0^w(X)$. The spaces $C_0^w(X)$ and
$C_0(X)$ are lattice isomorphic, a lattice isomorphism between
$C_0(X)$ and $C_0^w(X)$ being the linear operator
$S:C_0(X)\rightarrow C_0^w(X)$ defined by
\begin{equation}\label{6.10}
S(f):=w^{-1}f\quad\quad(f\in C_0(X)).
\end{equation} Therefore, \textsl{if
$M$ is a subset of $C_0^w(X)$, then $M$ is a Korovkin subset
of $C_0^w(X)$ if and only if
\begin{equation}\label{6.11} wM:=\{wf\mid f\in M\}
\end{equation} is a Korovkin subset of $C_0(X)$}.

Accordingly, from Proposition 6.4 and Theorem 6.5, we immediately
obtain the following result.
\begin{corollary}\sl  Given a subset $M$ of $C_0^w(X)$, the
following statements hold:
\begin{itemize}
\item [(1)] If for every $x,y\in X$, $x\neq y$, there exists $h\in
\mathcal{L}(M), h\geq 0$, such that $h(x)=0<h(y)$, then $M$ is a
Korovkin subset of $C_0^w(X)$.  \item [(2)] Assume that $w\in
C_0(X)$. If $M^2\subset C_0^w(X)$ and $M$ separates the points of
$X$, then $\{\textbf{1}\}\cup M\cup M^2$ is a Korovkin subset of
$C_0^w(X)$. Moreover, if $M=\{f_1,\dots,f_n\}, n\geq 1$, is finite,
then
$\left\{\textbf{1},f_1,\dots,f_n,\sum\limits_{i=1}^nf_i^2\right\}$
is a Korovkin subset of $C_0^w(X)$.
\end{itemize} \end{corollary}

Applying the above part (2) to a subset $X$ of $\mathbb{R}^d$,
$d\geq 1$, and $M=\{pr_1,\dots,pr_d\}$, we obtain:

\begin{corollary}\sl  Consider a locally compact subset $X$ of
$\mathbb{R}^d$, $d\geq 1$, and $w\in C_0(X)$ such that
$\norm{\cdot}^2w\in C_0(X)$.  Then
$\{\textbf{1},pr_1,\dots,pr_d,\norm{\cdot}^2\}$ is a Korovkin subset
of $C_0^w(X)$.
\end{corollary}

From the previous Corollary 6.13, it follows in particular that, if
$X$ is a noncompact real interval, then $\{\textbf{1},e_1,e_2\}$ is
a Korovkin subset of $C_0^w(X)$ for any $w\in C_0(X)$ such that $e_2
w\in C_0(X).$ When $X\subset[0,+\infty[$, this result can be
considerably generalized.

\begin{corollary}\sl  Let $X$ be a noncompact subinterval of $[0,+\infty[$ and
let $w\in C_0(X)$ be a weight on $X$.  Consider
$\lambda_1,\lambda_2\in\mathbb{R}$ such that $0<\lambda_1<\lambda_2$
and $e_{\lambda _{2}}\in C_0^w(X)$.  Then $\{\textbf{1},e_{\lambda
_{1}},e_{\lambda _{2}}\}$ is a Korovkin subset of $C_0^w(X)$.
\end{corollary}

\begin{proof} First note that
$e_{\lambda _{1}}\in C_0^w(X)$ because $x^{\lambda_1}\leq
1+x^{\lambda_2}$ for every $x\geq 0$.  In order to get the result,
we shall apply part (1) of Corollary 6.12 by showing that for a
given $x_0\in X$ there exists
$h\in\mathcal{L}(\{\textbf{1},e_{\lambda _{1}},e_{\lambda _{2}})$
such that $h(x_0)=0<h(y)$ for any $y\in X, y\neq x_0$.

Set $a:=\inf X\geq 0$.  If $x_0=a$, then it is sufficient to
consider $h(x)=x^{\lambda_1}-a^{\lambda_1}\;(x\in X)$.  If $X$ is
upper bounded and $x_0=\sup X$, then we may consider
$h(x)=x_0^{\lambda_1}-x^{\lambda_1}\;(x\in X)$. Finally if $x_0$
belongs to the interior of $X$, then by means of differential
calculus it is not difficult to show that the function
$$h(x):=(\lambda_2-\lambda_1)x_0^{\lambda_2}-\lambda_2x_0^{\lambda_2-\lambda_1}
x^{\lambda_1}+\lambda_1x^{\lambda_2}\quad(x\in
X)$$ satisfies the required properties.
\end{proof}

From (\ref{6.11}) and
Proposition (\ref{6.8}), (3), we also get:
\begin{corollary}\sl
Consider the functions $g_{\beta }(x):=\exp (\beta x)$ and
$g_{\gamma }(x):=\exp (\gamma x)$ $(x\in \lbrack 0,+\infty \lbrack
),$ where $0<\beta
<\gamma .$ Then for every $\alpha >\gamma ,$ the subset $\left\{ \mathbf{1}%
,g_{\beta },g_{\gamma }\right\} $ is a Korovkin subset of the space
$$E_0^\alpha:=\{f\in
C([0,+\infty[)\;|\;\lim\limits_{x\rightarrow +\infty}\exp(-\alpha
x)f(x)=0\}$$ endowed with the weighted norm
$\norm{f}_\alpha:=\underset{x\geq 0}{\sup}\;\exp(-\alpha
x)|f(x)|\;(f\in E_0^\alpha)$).
\end{corollary}

In the case when $w\in C_b(X)$, another easy method of finding
Korovkin subsets of $C_0^w(X)$ is indicated below.

\begin{proposition}\sl  Consider a continuous bounded weight $w\in
C_b(X)$.  Then every Korovkin subset of $C_0(X)$ is a Korovkin
subset of $C_0^w(X)$ as well.
\end{proposition}

\begin{proof} Let $M$ be a Korovkin subset of $C_0(X)$.
In order to show that $M$ is a Korovkin subset of $C_0^w(X)$, we
shall show that $wM$ is a Korovkin subset of $C_0(X)$.  To this end,
fix $\mu\in\mathcal{M}^+(X)$ and $x\in X$ such that
$\mu(wg)=w(x)g(x)$ for every $g\in M$, and consider the positive
linear form $\nu:C_0(X)\longrightarrow\mathbb{R}$ defined by
$\nu(f):=\mu(wf)/w\;(f\in C_0(X))$.  Then $\nu\in\mathcal{M}_b^+(X)$
and $\nu(g)=g(x)$ for every $g\in M$. From Theorem 6.1, we then
conclude that $\nu(f)=f(x)$ for every $f\in C_0(X)$, i.e.,
$\mu(wf)=w(x)f(x)$ for every $f\in C_0(X)$.

Note that if $\varphi\in K(X)$, after setting $f:=\varphi/w\in
K(X)$, we get $\varphi=wf$ and hence
$\mu(\varphi)=\varphi(x)$.  By density, we conclude that
$\mu(f)=f(x)$ for every $f\in C_0(X)$ and hence the desired result
follows from Theorem 6.1.
\end{proof}

Two simple applications of Corollary 6.13 and Proposition 6.16 are
shown below. The first one is concerned again with the
Sz\'{a}sz-Mirakjan operators (6.7) on the weighted continuous
function spaces
\begin{equation}\label{6.12}
C_0^{w_m}([0,+\infty[):=\left\{f\in C([0,+\infty[)\mid
\lim_{x\rightarrow\infty}\frac{f(x)}{1+x^m}=0\right\},
\end{equation} $m\geq 1$, where $w_m(x):=\frac{1}{1+x^m}\; (x\geq 0).$

\begin{theorem}\sl
For every $m\geq 1$ and $n\geq 1$ and for every $f\in
C_0^{w_m}([0,+\infty[)$, $M_n(f)\in C_0^{w_m}([0,+\infty[)$ and
$\lim\limits_{n\rightarrow\infty}M_n(f)=f$ on $[0,+\infty[$ with
respect to the weighted norm $\norm{\cdot}_{w_m}$ and, hence,
uniformly on compact subsets of $[0,+\infty[$.
\end{theorem}

\begin{proof}In [32, Lemma 5], it was shown that every $M_n$ is a bounded
linear operator from $C_0^{w_m}$ into $C_0^{w_m}$ and that
$\underset{x\geq 0}{\sup}\norm{M_n}<+\infty$. For every $\lambda>0$
consider again the function $f_{\lambda}(x):=\exp(-\lambda
x)\;(x\geq 0)$.  By Propositions 6.8 and 6.16, the subset
$\{f_{\lambda_1}, f_{\lambda_2},f_{\lambda_3}\}$ is a Korovkin
subset of $\mathcal{C}_0^{w_m}([0,+\infty[)$ provided that
$0<\lambda_1<\lambda_2<\lambda_3$.  On the other hand, for each
function $f_\lambda$, we have already shown that
$\lim\limits_{n\rightarrow\infty} M_n(f_\lambda)=f_\lambda$
uniformly on $[0,+\infty[$ (see the proof of Theorem 6.11) and hence
the same limit relationship holds with respect to
$\norm{\cdot}_{w_m}$ because
$\norm{\cdot}_{w_m}\leq\norm{\cdot}_\infty$ on $C_0([0,+\infty[)$
and hence the result follows.
\end{proof}
\begin{remarks}\rm \begin{itemize}\item []
\item [1.]The behaviour of Sz\'{a}sz-Mirakjan operators on locally
convex weighted function spaces has been further investigated in
[11].
\item [2.]A result similar to Theorem 6.17 can be also proved for the
Baskakov operators [27]
\begin{equation}\label{6.13}
B_n(f)(x):=\frac{1}{(1+x)^n}\sum_{h=0}^n\binom{n+h-1}{h}f\Big(\frac{h}{n}\Big)\Big(
\frac{x}{1+x}\Big)^h
\end{equation} ($f\in C_0^{w_m}([0,+\infty[)$, $x\geq
0$, $m\geq 1$, $n\geq 1$).
\end{itemize} \end{remarks}

For the details, we
refer to [20, Theorem 2.3] (see also [8, pp.  341-344] and [32]).

As another application consider $X=\mathbb{R}$ and the weight
$w_m(x):=\frac{1}{1+x^m}\;(x\in\mathbb{R})$ with $m\geq 4$, $m$
even. Consider the sequence of the Gauss-Weierstrass operators
defined by (\ref{4.44}) for $p=1$.  Hence
\begin{equation}\label{6.14}
G_n(f)(x):=\Big(\frac{n}{4\pi}\Big)^{1/2}\int_\mathbb{R}f(t)\exp\Big(-\frac{n}{4}(t-x)^2\Big)\dd
t
\end{equation}
for $n\geq 1$ and $x\in\mathbb{R}$ where
$f:\mathbb{R}\longrightarrow\mathbb{R}$ is any Borel measurable
function for which the integral (\ref{6.14}) is absolutely
convergent. In particular, the operators $G_n$ are well-defined on
the function space
\begin{equation}\label{6.15} C_0^{w_m}(\mathbb{R}):=\left\{f\in
C(\mathbb{R})\mid \lim_{|x|\rightarrow\infty}\frac{f(x)}{1+x^m}=0\right\}.
\end{equation}

\begin{theorem}\sl  Under the above hypotheses,
$G_n(C_0^{w_n}(\mathbb{R}))\subset C_0^{w_m}(\mathbb{R})$ and
$$\lim_{n\rightarrow\infty}G_n(f)=f\quad\text{with respect to}\
\norm{\cdot}_{w_m}$$ (and hence uniformly on compact subsets of
$\mathbb{R}$) for any $f\in C_0^{w_m}(\mathbb{R})$.
\end{theorem}

\begin{proof}
Consider the function $e_m(x)=x^m\;(x\in\mathbb{R})$.  Then, for any
$n\geq 1$ and $x\in\mathbb{R}$,
\begin{eqnarray*}
G_n(e_m)(x)&=&\Big(\frac{n}{4\pi}\Big)^{1/2}\int_\mathbb{R}t^m\exp\Big(-\frac{n}{4}(t-x)^2\Big)\dd t\\
&=&\sum_{k=0}^m\binom{m}{k}M_k\Big(\frac{2}{n}\Big)^{k/2}x^{m-k}
\end{eqnarray*} where $M_0:=1$ and $$M_k:=\begin{cases}
0&\text{if}\;k\;\text{is odd},\\ 1\cdot
3\cdots(k-1)&\text{if}\;k\;\text{is even},
\end{cases}$$ (see [29,
formulae (\ref{4.20}) and (\ref{4.21})]).  Therefore,
\begin{eqnarray*}
G_n(e_m)(x)&\leq&\sum_{k=0}^m\binom{m}{k}M_k\Big(\frac{2}{n}\Big)^{k/2}|x|^{m-k}
\leq  x^m+\sum_{k=1}^m\binom{m}{k}M_k2^{k/2}|x|^{m-k}
\end{eqnarray*}
and
$$w_m(x)G_n(1+e_m)(x)\leq\frac{1}{1+x^m}\Big(x^m+\sum_{k=1}^m\binom{m}{k}M_k2^{k/2}|x|^{m-k}\Big)\leq
M,$$ $M$ being a suitable positive constant independent of $m$.

Consider now $f\in C_0^{w_m}(\mathbb{R})$.  Thus to each
$\varepsilon>0$ there corresponds $\delta_1>0$ such that for any
$x\in\mathbb{R}, |x|\geq\delta_1$,
$$w_m(x)|f(x)|\leq\frac{\varepsilon}{2M}.$$ Furthermore, choose
$\delta\geq\delta_1$ such that
$$w_m(x)\leq\frac{\varepsilon}{2M_\delta}\quad\text{for}\quad
x\in\mathbb{R}, |x|\geq\delta,$$ where
$M_\delta:=\underset{|x|\leq\delta}{\sup}|f(t)|$.  Then, for any
$x\in\mathbb{R}$, $|x|\geq\delta$,
\begin{equation*} \begin{split}
&|G_n(f)(x)|\leq
\Big(\frac{n}{4\pi}\Big)^{1/2}\int_{|t|\leq\delta}|f(t)|\exp\Big(-\frac{n}{4}(t-x)^2\Big)\dd t\\
&+\Big(\frac{n}{4\pi}\Big)^{1/2}\int_{|t|\geq\delta}|f(t)|\exp\Big(-\frac{n}{4}(t-x)^2\Big)\dd t\\
&\leq
M_\delta+\frac{\varepsilon}{2M}\Big(\frac{n}{4\pi}\Big)^{1/2}\int_{|t|\geq\delta}(1+e_m(t))
\exp\Big(-\frac{n}{4}(t-x)^2\Big)\dd t\\
&\leq M_\delta+\frac{\varepsilon}{2M}G_n(1+e_m)(x)
\end{split}
\end{equation*} so that $$w_m(x)|G_n(f)(x)|\leq\varepsilon$$ because of
the previous estimates.

This proves that $G_n(f)\in C_0^{w_m}(\mathbb{R})$.  Moreover, note
that, since $|f|\leq\norm{f}_{w_m}(1+e_m)$, then
$$w_m|G_n(f)|\leq\norm{f}_{w_m}w_mG_n(1+e_m)$$ and hence
$\norm{G_n(f)}_{w_m}\leq M\norm{f}_{w_m}$, that is, $\norm{G_n}\leq
M$ for any $n\geq 1$.

On the other hand, from Corollary 6.13, we infer that
$\{\textbf{1},e_1,e_2\}$ is a Korovkin subset of
$C_0^{w_m}(\mathbb{R})$, where $e_k(x):=x^k$ ($k=1,2$, $x\in\mathbb{R}$),
and for every $h\in\{\textbf{1},e_1,e_2\}$ clearly
$G_n(h)\rightarrow h$ as $n\rightarrow\infty$ with respect to
$\norm{\cdot}_{w_n}$ because of (\ref{4.48})-(\ref{4.50}).
Therefore, the result follows.
\end{proof}
\begin{remark}\rm  A generalization of Theorem 6.19 can be found in [21,
Theorem 3.1] and [23].
\end{remark}

\section{Korovkin-type theorems for the identity operator on $C(X)$, $X$
compact} The results of the previous section also apply when $X$ is
a compact metric space (and hence $C_0(X)=C(X)$). However, in this
particular case, some of them have a particular relevance and hence
are worth an explicit description.

From Theorem 6.5, by replacing $f_0$ with the constant function
$\textbf{1}$, the following result immediately follows.
\begin{theorem}\sl  If $M$ is a subset of $C(X)$ that
separates the points of $X$, then $\{\textbf{1}\}\cup M\cup M^2$ is
a Korovkin subset of $C(X)$.

Moreover, if $M=\{f_1,\dots,f_n\}$, $n\geq 1$, is finite, then
$\left\{\textbf{1},f_1,\dots,f_n, \sum\limits_{i=1}^nf_i^2\right\}$
is a Korovkin subset of $C(X)$.

In particular, if $f\in C(X)$ is injective, then $\{\textbf{1}, f,
f^2\}$ is a Korovkin subset of $C(X)$.
\end{theorem}

Theorem 7.1 extends both the Korovkin Theorems 3.1 and 4.3 as well
as Volkov's Theorem 4.2.  It was obtained by Freud ([59]) for finite
subsets $M$ and by Schempp ([106]) in the general case. By adapting
the proof of Proposition 6.8 and by using Proposition 6.4, it is not
difficult to show that:

\begin{proposition}\sl

\begin{itemize}
  \item [(1)] If $0<\lambda_1<\lambda_2$, then
  $\{\textbf{1},e_{\lambda
_{1}},e_{\lambda _{2}}\}$ is a Korovkin subset of
  $C([0,1])$, where $e_{\lambda _{k}}(x):=x^{\lambda _{k}}$
for every $\ x\in \lbrack 0,1\rbrack $ and $k=1,2$.  \item [(2)] If
$u\in C([0,1])$ is strictly convex, then
  $\{\textbf{1},e_1,u\}$ is a Korovkin subset of $C([0,1])$.
\end{itemize} \end{proposition}

Next, we shall discuss some applications of Theorem 7.1 by
considering a convex compact subset $K$ of a locally convex space
$E$. We denote by $$A(K)$$ the subspace of all real-valued
continuous affine functions on $K$. We recall that a function
$u:K\rightarrow \mathbb{R}$ is said to be affine if
\begin{equation}\label{7.1}u(\lambda x+(1-\lambda )y)=\lambda u(x)+(1-\lambda )u(y)
\end{equation}
for every $x,y\in K$ and $\lambda \in \lbrack 0,1],$ or, equivalently, if%
\begin{equation}\label{7.2}
u(\sum_{i=1}^{n}\lambda _{i}x_{i})=\sum_{i=1}^{n}\lambda
_{i}u(x_{i})
\end{equation}
for every $n\geq 2$, $x_{1},\ldots ,x_{n}\in K$   and
$\lambda _{1},\ldots ,\lambda _{n}\geq 0$ such that
$\sum_{i=1}^{n}\lambda _{i}=1.$

Note that $A(K)$ contains the constant functions as well as the
restrictions to $K$ of every continuous linear functional on $E$ so
that, by the Hahn-Banach theorem, $A(K)$ separates the points of
$K$. Therefore, from Theorem 7.1, we immediately get:

\begin{corollary}\sl
$A(K)\cup A(K)^2$ is a Korovkin subset of $C(K)$.
\end{corollary}
Consider now a continuous selection $(\widetilde{\mu}_x)_{x\in K}$
of probability Borel measures on $K$, i.e., for every $f\in C(K)$
the function $x\mapsto\int\limits_Kf\dd\widetilde{\mu}_x$ is
continuous on $K$. Such a function will be denoted by $T(f)$, that
is
\begin{equation}\label{7.3} T(f)(x):=\int_Kf\dd\widetilde{\mu}_x\quad (x\in
K).
\end{equation} Note that $T$ can be viewed as a positive linear
operator from $C(K)$ into $C(K)$ and $T(\textbf{1})=\textbf{1}$.

Conversely, by the Riesz representation theorem (see Theorem 11.3 of
the \Appendix), each positive linear operator $T:C(K)\longrightarrow
C(K)$ generates a continuous selection of probability Borel measures
on $K$ satisfying (\ref{7.3}).  Such a selection will be referred to
as the \dword{canonical continuous selection associated with} $T$.

From now on, we fix a continuous selection
$(\widetilde{\mu}_x)_{x\in K}$ of probability Borel measures on $K$
satisfying
\begin{equation}\label{7.4} \int_Ku\dd\widetilde{\mu}_x=u(x)\quad\text{for
every}\quad u\in A(K),
\end{equation} namely,
\begin{equation}\label{7.5}
T(u)=u\quad\text{for every}\quad u\in A(K).
\end{equation} Property
(\ref{7.4}) means that each $x\in K$ is the barycenter or the
resultant of $\widetilde{\mu}_x$ ([39, Vol.~II, Def.~26.2]).

For every $n\geq 1$ and $x\in K$, denote by
$$\widetilde{\mu}_x^{(n)}$$ the product measure
on $K^n$ of $\widetilde{\mu}_x$ with itself $n$ times (see, e.g.,
[30, \S  22]) and set, for every $f\in C(K)$,
\begin{equation}\label{7.6}
B_n(f)(x):=\int_{K^n}f\Big(\frac{x_1+\dots+x_n}{n}\Big)\dd\widetilde{\mu}_x^{(n)}(x_1,\dots,x_n).
\end{equation} By Fubini's theorem, we can also write
\begin{equation}\label{7.7}
B_n(f)(x)=\int_K\cdots\int_Kf\Big(\frac{x_1+\dots+x_n}{n}\Big)\dd\widetilde{\mu}_x(x_1)\cdots
d\widetilde{\mu}_x(x_n).
\end{equation} By the continuity property of the
product measure (see [39, Vol.~I, Proposition 13.12] and [30,
Theorem 30.8]) and the continuity of the selection, it follows that
$B_n(f)\in C(K)$.

The positive linear operator $B_n:C(K)\longrightarrow C(K)$ is
referred to as the \dword{$n$-th Bernstein-Schnabl operator}
associated with the given selection $(\widetilde{\mu}_x)_{x\in K}$
(or the given positive linear operator $T:C(K)\longrightarrow
C(K)).$

Bernstein-Schnabl operators were introduced for the first time in
1968 by Schnabl ([107-108]) in the context of the set of all
probability Radon measures on a compact Hausdorff space and, as we
shall see next, they generalize the classical Bernstein operators
(\ref{3.5}). Subsequently, Grossman ([66]) introduced the general
definition (\ref{7.6}) (or (\ref{7.7})). In the particular case of
Bauer simplices (see, e.g., [8, Section 1.5, p. 59]), these
operators have been extensively studied by Nishishiraho ([87-91]).
Another particular class of them has been also studied by Altomare
([5]). Their construction essentially involves positive projections
and, in this case, they satisfy many additional properties useful in
the study of evolution problems (see the end of Section 10). For
some shape preserving properties of these operators, we also refer
to Ra\c{s}a ([97-99]).

For a comprehensive survey on these operators, we refer to [8,
Chapter 6] and to the references contained in the relevant notes.
More recent results can be also found in [7], [18], [19], [24-25],
[100-102], [103].

Below, we discuss some examples.
\begin{examples}\rm
\begin{itemize}\item []
\item [1.]Consider $K=[0,1]$ and set
$\widetilde{\mu}_x:=x\delta_1+(1-x)\delta_0$ for all $x\in[0,1]$.
Then $(\widetilde{\mu}_x)_{0\leq x\leq 1}$ is a continuous selection
and the corresponding Bernstein-Schnabl operators turn into the
Bernstein operators (\ref{3.5}) (use an induction argument).  \item
[2.]Let $\alpha,\beta,\gamma\in C([0,1])$ and assume that
$0\leq\alpha\leq 1,\;0\leq\beta\leq 1,\;0\leq\gamma\leq
1,\;\alpha+\beta+\gamma=\textbf{1}$ and $\beta(x)/2+\gamma(x)=x$ for
every $x\in[0,1]$.  Put
$\widetilde{\mu}_x:=\alpha(x)\delta_0+\beta(x)\delta_{1/2}+\gamma(x)\delta_1$.
Then the Bernstein-Schnabl operators associated with
$(\widetilde{\mu}_x)_{0\leq x\leq 1}$ on $[0,1]$ are given by
\begin{equation}\label{7.8}
B_n(f)(x)=\sum_{h=0}^n\sum_{k=0}^n\binom{n}{h}\binom{n-h}{k}\alpha(x)^{n-h-k}\beta(x)^h
\gamma(x)^kf\Big(\frac{h+2k}{2n}\Big)
\end{equation} $(f\in C([0,1]),\;0\leq x\leq 1)$.
\item [3.]Let $a=a_0<a_1<\dots<a_p=b$, $p\geq 1$, be a subdivision of the compact real
interval $[a,b]$.  For every $x\in[a,b]$, set
\begin{equation}\label{7.9}
\widetilde{\mu}_x:=\frac{x-a_k}{a_{k+1}-a_k} \delta_{a_{k+1}} +
\frac{a_{k+1}-x}{a_{k+1}-a_k} \delta_{a_{k}},
\end{equation} whenever $x \in [a_k, a_{k+1}]$, $0 \leq k \leq p-1$.
Then the Bernstein-Schnabl operators are given by
\begin{equation}\label{7.10}
B_n(f)(x)=\frac{1}{(a_{k+1}-a_k)^n}\sum_{r=0}^n\binom{n}{r}(x-a_k)^r(a_{k+1}-x)^{n-r}
f\Big(\frac{r}{n}a_{k+1}+\frac{n-r}{n}a_k\Big)
\end{equation} whenever
$x\in[a_k,a_{k+1}],\;0\leq k\leq p-1\;(f\in C([a,b]),\;n\geq 1).$
\item [4.]Consider the $d$-dimensional simplex $K_d$ of
$\mathbb{R}^d$ defined by (4.38) and for every $x=(x_1,\dots,x_d)\in
K_d$, set
$$\widetilde{\mu}_x:=\Big(1-\sum_{i=1}^dx_i\Big)\delta_0+\sum_{i=1}^dx_i\delta_{a_i}$$
where $a_i=(\delta_{ij})_{1\leq j\leq d}$ for every
$i=1,\dots,d,\;\delta_{ij}$ being the Kronecker symbol.  Then the
Bernstein-Schnabl operators associated with
$(\widetilde{\mu}_x)_{x\in K_d}$ turn into the Bernstein operators
on the simplex $K_d$ defined by (\ref{4.38}).  \item [5.]Consider
the hypercube $Q_d:=[0,1]^d$ of $\mathbb{R}^d$ and for every
$x=(x_1,\dots,x_d)\in Q_d$, set
$$\widetilde{\mu}_x:=\sum_{h_1,\dots,h_d=0}^1x_1^{h_1}(1-x_1)^{1-h_1}\cdots
x_d^{h_d}(1-x_d)^{1-h_d}\delta_{b_{h_1,\dots,h_d}},$$ where
$b_{h_1,\dots,h_d}=(\delta_{h_11},\dots,\delta_{h_d1})\quad
(h_1,\dots,h_d\in\{0,1\}).$ Then the Bernstein-Schnabl operators are
the Bernstein operators on $Q_d$ defined by (\ref{4.41}).
\end{itemize} \end{examples} Many other
significant examples can be described in the framework of other
finite dimensional subsets such as balls, or more generally,
ellipsoids of $\mathbb{R}^d$, or in the setting of
infinite-dimensional Bauer simplices (see [8, Chapter 6], [24],
[103]).

In order to discuss the approximation properties of
Bernstein-Schnabl operators, note that for every $n\geq 1$
\begin{equation}\label{7.11}
B_n(\textbf{1})=\textbf{1}.
\end{equation} Moreover, if $u\in A(K)$, then
for every $x_1,\dots,x_n\in K$
$$u\Big(\frac{x_1+\dots+x_n}{n}\Big)=\frac{u(x_1)+\dots+u(x_n)}{n}$$ and
$$u^2\Big(\frac{x_1+\dots+x_n}{n}\Big)=\frac{1}{n^2}\Big[\sum_{i=1}^nu^2(x_i)+2\sum_{1\leq
i<j\leq n}u(x_i)u(x_j)\Big].$$ Therefore, because of (\ref{7.4}),
for every $x\in K$, we get
\begin{equation}\label{7.12} B_n(u)(x)=u(x)
\end{equation} and
\begin{equation}\label{7.13}
B_n(u^2)(x)=\frac{1}{n}T(u^2)(x)+\frac{n-1}{n}u^2(x).
\end{equation} In other words, for every $h\in A(K)\cup
A(K)^2,\;B_n(h)\rightarrow h$ uniformly on $K$ and hence by virtue
of Corollary 7.3, we obtain that:

\begin{theorem}\sl  For every $f\in C(K)$,
$$\lim_{n\rightarrow\infty}B_n(f)=f\quad\text{uniformly on}\quad K.$$
\end{theorem}

For additional properties of Bernstein-Schnabl operators, including
estimates of the rate of convergence, asymptotic formulae,
shape-preserving properties and, especially, their connections with
the approximation of positive semigroups as well as of the solutions
of evolution equations, we refer to the references we cited before
Examples 7.4.

\section{Korovkin-type theorems in weighted continuous function spaces and
in $L^p(X,\widetilde{\mu})$ spaces} By using Korovkin-type theorems
in spaces of continuous functions it is also possible to get some
results for $L^p(X,\widetilde{\mu})$-spaces, $1\leq p<+\infty$. In
this regard, other than the space $C_0(X)$, also weighted continuous
function spaces can be efficiently used. In this section, we shall
develop this method by also showing some additional Korovkin-type
results in weighted continuous function spaces that complement the
ones already discussed at the end of Section 6, and then by
presenting the corresponding ones in
$L^p(X,\widetilde{\mu})$-spaces.

The method relies on a fundamental characterization of Korovkin
subspaces which is due to Bauer and Donner ([31]) (see also [8,
Theorem 4.1.2]).  Its counterpart for Korovkin subspaces for
positive linear operators was obtained in [8, Theorem 3.1.4].  An
extension to locally convex function spaces was obtained in [10].
For the sake of brevity, we omit its proof.

As in the previous section, $X$ denotes a fixed locally compact
Hausdorff space with a countable base.
\begin{theorem}\sl  ([31]) Given a linear
subspace $H$ of $C_0(X)$, the following statements are equivalent:
\begin{itemize}
\item [(i)] $H$ is a Korovkin subspace of $C_0(X)$; \item [(ii)] For every $f\in C_0(X)$ and
for every $\varepsilon>0$, there exist finitely many functions
$h_0,\dots,h_n\in H$, $k_0,\dots,k_n\in H$ and $u,v\in C_0(X),
u,v\geq 0$ such that $\norm{u}\leq\varepsilon,
\norm{v}\leq\varepsilon$ and
$$\Norm{\underset{0\leq j\leq n}{\inf}k_j-\underset{0\leq i\leq
n}{\sup}h_i}\leq\varepsilon\quad\text{and}\quad\underset{0\leq i\leq
n}{\sup}h_i-u\leq f\leq\underset{0\leq j\leq n}{\inf}k_j+v.$$
\end{itemize}
\end{theorem}

\begin{remark}\rm
If $X$ is compact and $\textbf{1}\in H$, then Theorem 8.1 can be
stated in a simpler form (see [8, Theorem 4.1.4]).
\end{remark}

Let us mention two important consequences of Theorem 8.1.

\begin{theorem}\sl  ([31]).  Let $H$
be a Korovkin subspace of $C_0(X)$.  If $E$ is a Banach lattice and
if $S:C_0(X)\rightarrow E$ is a lattice homomorphism, then $H$ is a
Korovkin subspace of $C_0(X)$ for $S$.
\end{theorem}

\begin{proof}
Consider a sequence $(L_n)_{n\geq 1}$ of positive linear operators
from $C_0(X)$ into $E$ such that $M:=\underset{n\geq
1}{\sup}\norm{L_n}<+\infty$ and
$\lim_{n\rightarrow\infty}L_n(h)=S(h)$ for every $h\in H$.  Given
$f\in C_0(X)$ and $\varepsilon>0$, choose $\delta\in]0,\varepsilon[$
such that $\norm{S(g)}\leq\varepsilon$ for every $g\in C_0(X)$,
$\norm{g}\leq\delta$. By Theorem 8.1, there exist $h_0,\dots,h_p,
k_0,\dots,k_p\in H$ and $u,v\in C_0(X)$, $u,v\geq 0$, such that
$\norm{u}\leq\delta$, $\norm{v}\leq\delta$ and
$$
\norm{v}\leq\varepsilon$$ as well as
$$\Norm{\underset{0\leq j\leq n}{\inf}S(k_j)-\underset{0\leq i\leq
n}{\sup}s(h_i)}\leq\delta\;\text{and}\;\underset{0\leq i\leq
p}{\sup}h_i-u\leq f\leq\underset{0\leq j\leq p}{\inf}k_j+v.$$
Therefore,
$$\Norm{\underset{0\leq j\leq p}{\inf}S(k_j)-\underset{0\leq i\leq
p}{\sup}S(h_i)}\leq\varepsilon.$$ Moreover, for every $n\geq 1$,
$$\underset{0\leq i\leq p}{\sup}L_n(h_i)-L_n(u)\leq
L_n(f)\leq\underset{0\leq j\leq p}{\inf}L_n(k_j)+L_n(v)$$ and
$$\underset{0\leq i\leq p}{\sup}S(h_i)-S(u)\leq S(f)\leq\underset{0\leq
j\leq p}{\inf}S(k_j)+S(v)\,.$$ Accordingly,
\begin{eqnarray*}
L_n(f)-S(f)\leq\sum_{j=0}^p|L_n(k_j)-S(k_j)|+|\underset{0\leq j\leq
p}{\inf}S(k_j)-\underset{0\leq i\leq p}{\sup}S(h_i)| +L_n(v)+S(u)
\end{eqnarray*} and
\begin{eqnarray*}
S(f)-L_n(f)\leq\sum_{i=0}^p|L_n(h_i)-S(h_i)|+|\underset{0\leq j\leq
p}{\inf}S(k_j)-\underset{0\leq i\leq p}{\sup}S(h_i)|+
+L_n(u)+S(v),\end{eqnarray*} so that
\begin{eqnarray*}
|S(f)-L_n(f)|&\leq&\sum_{i=0}^p|L_n(h_i)-S(h_i)|+\sum_{j=0}^p|L_n(k_j)-S(k_j)|\\
&+&|\underset{0\leq j\leq p}{\inf}S(k_j)-\underset{0\leq i\leq
p}{\sup}S(h_i)|+L_n(u)+L_n(v)+S(u)+S(v)\\
&\leq&\sum_{i=0}^p|L_n(h_i)-S(h_i)|+\sum_{j=0}^p|L_n(k_j)-S(k_j)|+2M\varepsilon+2\varepsilon.
\end{eqnarray*}
It is now easy to conclude that $L_n(f)\rightarrow S(f)$ as
$n\rightarrow\infty$ because $L_n(h_i)\rightarrow S(h_i)$ and
$L_n(k_j)\rightarrow S(k_j)$ as $n\rightarrow\infty$ for every
$i,j=0,\dots,p$.
\end{proof}

\begin{remark}\rm  From Theorem 8.3, it follows that, if $M$ is a Korovkin subset of $C_0(X)$, then $M$ is a
Korovkin subset for the natural embedding from $C_0(X)$ into an
arbitrary closed lattice subspace $E$ of $B(X)$ containing $C_0(X)$
(for instance, $E=C_b(X)$ or $E=B(X)$).
\end{remark} Another consequence of Theorem 8.3 is indicated below.

\begin{proposition}\sl  ([120]).  Let $E$ be a Banach lattice and consider a
lattice homomorphism $S:C_0(X)\rightarrow E$ such that $S(C_0(X))$
is dense in $E$.  If $M$ is a Korovkin subspace of $C_0(X)$, then
$S(M)$ is a Korovkin subspace of $E$.
\end{proposition}

\begin{proof}
Consider an equibounded sequence $(L_n)_{n\geq 1}$ of positive
linear operators on $E$ and assume that $L_n(k)\rightarrow k$ as
$n\rightarrow\infty$ for every $k\in S(M)$. This means that
$L_n(S(h))\rightarrow S(h)$ as $n\rightarrow\infty$ for each $h\in
M$ and hence for each $h\in H:=\mathcal{L}(M)$.  So, by Theorem 8.3,
$L_n(S(f))\longrightarrow S(f)$ for every $f\in C_0(X)$.  The result
now follows from the assumption that $S(C_0(X))$ is dense in $E$ and
$(L_n)_{n\geq 1}$ is equibounded.
\end{proof}
\begin{remark}\rm A typical situation where
Proposition 8.5 can be applied concerns the case when $E$ is a
Banach lattice containing $C_0(X)$ as a dense sublattice and $S$ is
the natural embedding from $C_0(X)$ into $E$.  Thus, in this case,
\center{\it{every Korovkin subspace of $C_0(X)$ is a Korovkin
subspace of $E$ as well}}. \end{remark} \medskip
After these preliminaries, we
can now proceed to discuss some Korovkin-type results for
$L^p(X,\widetilde{\mu})$ spaces.

Consider a Borel measure $\widetilde{\mu}$ on $X$ and, given
$p\in[1,+\infty[$, consider the space
$$L^p(X,\widetilde{\mu})$$ endowed with its natural norm
$\norm{\cdot}_p$ (for more details, we refer to the \Appendix\
(formulae (\ref{11.11}) and (\ref{11.12}))).

Since $\widetilde{\mu}$ is regular, $K(X)$ is dense in
$L^p(X,\widetilde{\mu})$ with respect to the convergence in $p$-th
mean (see Theorem 11.4 of \Appendix). Therefore, we get the
following useful results.

\begin{corollary}\sl  Let $M$ be a Korovkin subset
of $C_0(X)$.  Furthermore, let $Y$ be a locally compact Hausdorff
space with a countable base and consider a Borel measure
$\widetilde{\mu}$ on $Y$ and $p\in[1,+\infty[$.  Let
$S:C_0(X)\longrightarrow L^p(Y,\widetilde{\mu})$ be a lattice
homomorphism such that $K(Y)\subset S(C_0(X))$.  Then $S(M)$ is a
Korovkin subset of $L^p(Y,\widetilde{\mu})$.

In particular, if $\widetilde{\mu}$ is a finite Borel measure on
$X$, then $M$ is also a Korovkin subset of $L^p(X,\widetilde{\mu})$.
\end{corollary}

\begin{corollary}\sl  Let $X$ be a compact subset of $\mathbb{R}^d$, $d\geq 1$,
and consider a finite Borel measure $\widetilde{\mu}$ on $X$.  Then
$\{\textbf{1},pr_1,\dots,pr_d,\sum\limits_{i=1}^dpr_i^2\}$ is a
Korovkin subset of $L^p(X,\widetilde{\mu})$ for every
$p\in[1,+\infty[$.

If, in addition, $X$ is contained in some sphere of $\mathbb{R}^d$,
then $\{\textbf{1},pr_1,\dots,pr_d\}$ is a Korovkin subset of
$L^p(X,\widetilde{\mu})$.
\end{corollary}

\begin{corollary}\sl  Consider $0<\lambda_1<\lambda_2<\lambda_3$ and
$p\in[1,+\infty[$.  Then
\begin{itemize}
\item [(1)] $\{f_{\lambda _{1}},
f_{\lambda _{2}}, f_{\lambda _{3}}\}$ is a Korovkin subset of
$L^p([0,+\infty[)$, where $f_{\lambda _{k}}(x):=\exp (-\lambda
_{k}x)$ for every $\ x\in \lbrack 0,+\infty \lbrack $ and $k=1,2,3$.
\item [(2)] If $u:\mathbb{R}\longrightarrow ]0,1[$ is a
strictly increasing continuous function satisfying
$\lim\limits_{x\rightarrow+\infty}u(x)=0$ and
$\lim\limits_{x\rightarrow-\infty}u(x)=1$, then
$$\{\Phi ,\Phi u^{\lambda _{1}},\Phi u^{\lambda _{2}}\}$$
is a Korovkin subset of $L^p(\mathbb{R})$, where $%
\Phi (x):=\exp (-x^{2})$ \ $(x\in \mathbb{R)}$.
\end{itemize} \end{corollary}
\begin{proof} (1) Consider the lattice homomorphism $S:C([0,1])\rightarrow
L^p([0,+\infty[)$ defined by $$S(f)(x):=\exp(-\lambda_1
x)f(\exp(-x))\quad (f\in C([0,1]),\, x\geq 0).$$ Then $S$ maps the
subset $M:=\{\textbf{1},e_{\lambda _{2}-\lambda _{1}},e_{\lambda
_{3}-\lambda _{1}}\}$ into
$$\{f_{\lambda _{1}},
f_{\lambda _{2}}, f_{\lambda _{3}}\}.$$ Hence the result follows
from Corollary 8.7 and Proposition 7.2, (1).

(2) A similar reasoning can be used by considering now the lattice
homomorphism $S:C([0,1])\rightarrow L^p(\mathbb{R})$ defined by
$$S(f)(x):=\exp(-x^2)f(u(x))\quad (f\in C([0,1]),\, x\in\mathbb{R}).$$
\vskip-\baselineskip
\end{proof}

We now proceed to state a result analogous to Corollary 8.7 for more
general weighted function spaces (see Section 6). Consider a
continuous weight $w$ on $X$ and the relevant weighted space
$C_0^w(X)$ defined by (\ref{6.8}). Note that, if $\widetilde{\mu}$
is a Borel measure on $X$ and if
\begin{equation}\label{8.1} w^{-1}\in L^p(X,\widetilde{\mu}) \qquad
\textrm{ for some } p\in[1,+\infty[. \end{equation}
 then
$C_0^w(X)\subset L^p(X,\widetilde{\mu})$ and it is dense with
respect to the convergence in the $p$-th mean.

\begin{corollary}\sl  Consider a Korovkin subset $M$ of $C_0^w(X)$, i.e., $wM$
is a Korovkin subset of $C_0(X)$ (see (\ref{6.11})), and a Borel measure
$\widetilde{\mu}$ on $X$. If (8.1) holds, then $M$ is a
Korovkin subset of $L^p(X,\widetilde{\mu})$.
\end{corollary}

\begin{proof} Consider the lattice homomorphism
$S:C_0(X)\longrightarrow L^p(X,\widetilde{\mu})$ defined by
$S(f):=w^{-1}f$ $(f\in C_0(X))$. Then $M=S(wM)$ and
$S(C_0(X))=C_0^w(X)$.  Therefore, the result follows from Corollary
8.7.
\end{proof}

The previous Corollary together with Corollaries 6.12 - 6.15 furnish a simple
but useful method to construct Korovkin subsets in $L^p(X,\widetilde{\mu})$-spaces.

By using similar methods, we can extend some of the previous results
to more general weighted function spaces which often occur in the
applications. Some of the subsequent results are taken from [22, Section
2].

In what follows, we shall assume that $X$ is noncompact and we shall
denote by
$$X_{\infty}:=X\cup\{\infty\}$$ the Alexandrov one-point compactification of
$X$ (see, e.g., [30, \S 27]). A function $f\in C(X)$ is said to be
convergent at infinity if there exists a (unique) $l\in\mathbb{R}$
such that for any $\varepsilon>0$ there exists a compact subset $K$
of $X$ such that $|f(x)-l|\leq\varepsilon$ for each $x\in
X\backslash K$.  In such a case, we also write
$\lim\limits_{x\rightarrow\infty}f(x)=l$. Similarly, we shall write
$\lim\limits_{x\rightarrow\infty}f(x)=\infty$ to mean that for every
$M\geq 0$ there exists a compact subset $K$ of $X$ such that
$f(x)\geq M$ for every $x\in X\backslash K$.

Given a weight $w\in C(X)$, consider the Banach lattice
\begin{equation}\label{8.1}
C_*^w(X):=\{f\in C(X)\mid wf\;\text{is convergent at infinity}\;\}
\end{equation} endowed with the natural (pointwise) order and the weighted
norm
\begin{equation}\label{8.2} \norm{f}_w:=\norm{wf}_{\infty}\qquad (f\in
C_*^w(X)).
\end{equation} For every $f\in C_*^w(X)$, denote by $T(f)$ the
function on $X_{\infty}$ defined by
\begin{equation}\label{8.3}
T(f)(z):=\begin{cases} w(z)f(z)&\text{if}\;z\in X,\\
\lim\limits_{x\rightarrow\infty}w(x)f(x)&\text{if}\;z=\infty.
\end{cases}\end{equation} Then $T(f)\in C(X_{\infty})$ and
$\norm{T(f)}_{\infty}=\norm{f}_{\infty}$.  Moreover, the linear
operator $T:C_*^w(X)\rightarrow C(X_{\infty})$ is a lattice
isomorphism whose inverse we shall denote by
$$S:C(X_{\infty})\rightarrow C_*^w(X).$$
Thus a subset $M$ of $C_*^w(X)$ is a Korovkin subset of $C_*^w(X)$
if and only if $T(M)$ is a Korovkin subset of $C(X_{\infty})$.

Consider now a Borel measure $\widetilde{\mu}$ on $X$ and assume
that (8.1) holds.
Then $C_*^w(X)\subset L^p(X,\widetilde{\mu})$ and $C_*^w(X)$ is
dense in $L^p(X,\widetilde{\mu})$ because $K(X)\subset C_*^w(X)$.
\begin{proposition}\sl  Under  assumption (8.1), each Korovkin subset of $C_*^w(X)$ is a Korovkin
subset of $L^p(X,\widetilde{\mu})$ as well.
\end{proposition}
\begin{proof} The above lattice
isomorphism $S$ can be considered as a lattice homomorphism from
$C(X_{\infty})$ into $L^p(X,\widetilde{\mu})$ and its range is
$C_*^w(X)$ which is dense in $L^p(X,\widetilde{\mu})$.  Since
$H:=T(M)$ is a Korovkin subset of $C(X_{\infty})$, then $M=S(H)$ is
a Korovkin subset of $L^p(X,\widetilde{\mu})$ by Corollary 8.7.
\end{proof} Below, we state some consequences of Proposition 8.11.
\begin{proposition}\sl  Let $M$ be a subset
of $C_*^w(X)$ and assume that
\begin{itemize}
\item [(i)] for every $x_0,y_0\in X$, $x_0\neq y_0$, there exists $h\in
\mathcal{L}(M)$, $h\geq 0$, such that $h(x_0)=0$ and $h(y_0)>0$;
\item [(ii)] for every $x_0\in X$ there exist positive functions
$h,k\in \mathcal{L}(M)$ such that $$h(x_0)=0,\quad
\lim_{x\rightarrow\infty}w(x)h(x)>0,$$ and $$k(x_0)>0,\quad
\lim_{x\rightarrow\infty}w(x)k(x)=0.$$
\end{itemize} Then $M$ is a
Korovkin subset of $C_*^w(X)$ and hence of $L^p(X,\widetilde{\mu})$
for every Borel measure $\widetilde{\mu}$ on $X$ and for every $p\in[1,+\infty[$ satisfying (8.1).
\end{proposition}
\begin{proof} Conditions (i) and (ii) mean that $T(M)$ satisfies the
hypotheses of Propositions 6.4 in $C(X_{\infty})$.  Therefore, the
result follows from Propositions 6.4 and 8.11.
\end{proof}
\begin{corollary}\sl
Consider a subset $M$ of $C(X)$
that separates the points of $X$ and a strictly positive function
$f_0\in C_0^w(X)$.  Further, assume that
\begin{itemize}
\item [(1)] $f_0M\cup f_0M^2\subset C_*^w(X)$ (see (\ref{6.1}) and
(\ref{6.2})), \item [(2)] there exists $g\in M$ such that
$$\lim\limits_{x\rightarrow\infty}w(x)f_0(x)g(x)=0 \textrm{ and }
\lim\limits_{x\rightarrow\infty}w(x)f_0(x)g^2(x)>0.$$
\end{itemize} Then
$\{f_0\}\cup f_0M\cup f_0M^2$ is a Korovkin subset of $C_*^w(X)$ and
hence of $L^p(X,\widetilde{\mu})$ for every Borel measure $\widetilde{\mu}$ on $X$ and for every $p\in[1,+\infty[$ satisfying (8.1).
\end{corollary}
\begin{proof} We shall verify
conditions (i) and (ii) of Proposition 8.12.  Given $x_0,y_0\in X,
x_0\neq y_0$, there exists $f\in M$ such that $f(x_0)\neq f(y_0)$.
Therefore, the function
$h:=f_0(f-f(x_0))^2\in\mathcal{L}(\{f_0\}\cup f_0M\cup f_0M^2)$
satisfies $h(x_0)=0$ and $h(y_0)>0$ and hence property (i) follows.

As regards condition (ii) of Proposition 8.12, considering a
function $g\in M$ satisfying assumption (2), then for every $x_0\in
X$ the functions $h:=f_0(g-g(x_0))^2$ and $k:=f_0$ satisfy condition
(ii) of Proposition 8.12.
\end{proof}
\begin{corollary}\sl  Consider
$f_1,\dots,f_n\in C(X), n\geq 1$, that separate the points of $X$,
and a strictly positive function $f_0\in C_0^w(X)$. Further, assume
that
\begin{itemize}
\item [(i)] $f_0f_i\in C_0^w(X)$ for every $i=1,\dots,n$, as well as
$f_0\sum\limits_{i=1}^nf_i^2\in C_*^w(X)$; \item [(ii)]
$\lim\limits_{x\rightarrow\infty}w(x)f_0(x)\sum\limits_{i=1}^nf_i^2(x)>0$.
\end{itemize} Then
$\left\{f_0,f_0f_1,\dots,f_0f_n,f_0\sum\limits_{i=1}^nf_i^2\right\}$
is a Korovkin subset $C_*^w(X)$ and hence in
$L^p(X,\widetilde{\mu})$ for every Borel measure $\widetilde{\mu}$ on $X$ and for every $p\in[1,+\infty[$ satisfying (8.1).
\end{corollary}
\begin{proof} The
proof is similar to that of Corollary 8.13, except that in this
case, the function $h$ must be chosen as
$h:=f_0\sum\limits_{i=1}^n(f_i-f_i(x_0))^2$.
\end{proof} The particular case of Corollary 8.14 where $f_0=\textbf{1}$ is
worth stating separately.
\begin{corollary}\sl  Consider
$f_1,\dots,f_n\in C(X), n\geq 1,$ that separate the points of $X$
and $w\in C_0(X)$ a weight such that
\begin{itemize}
\item [(i)] $f_i\in C_0^w(X)$ for every $i=1,\dots,n$,
as well as $\sum\limits_{i=1}^nf_i^2\in C_*^w(X)$; \item [(ii)]
$\lim\limits_{x\rightarrow\infty}w(x)\sum\limits_{i=1}^nf_i^2(x)>0$.
\end{itemize} Then $\left\{\textbf{1},f_1,\dots,f_n\sum\limits_{i=1}^nf_i^2\right\}$
is a Korovkin subset in $C_*^w(X)$ and hence in
$L^p(X,\widetilde{\mu})$ for every finite Borel measure $\widetilde{\mu}$ on $X$ and for every $p\in[1,+\infty[$ satisfying (8.1).
\end{corollary} A simple
situation when the assumptions of Corollaries 8.14 and 8.15 are
satisfied is indicated below.
\begin{corollary}\sl  Consider
$f_1,\dots,f_n\in C(X), n\geq 1,$ that separate the points of $X$ and assume
that
$\lim\limits_{x\rightarrow\infty}\sum\limits_{i=1}^nf_i^2(x)=+\infty.$
The following statements hold:
\begin{itemize}
\item [(1)] If $f_0\in C(X)$   is a strictly
positive function, then
$\left\{f_0,f_0f_1,\dots,f_0f_n,f_0\sum\limits_{i=1}^nf_i^2\right\}$ is a
Korovkin subset in $C_*^w(X)$, where $$w:=\frac{1}{f_0\Big(1+\sum\limits_{i=1}^nf_i^2\Big)}.$$
Therefore, if $\widetilde{\mu}$ is a
Borel measure on $X$ and $1\leq p<+\infty$, and if $f_0\in L^p(X,\widetilde{\mu})$ as well as $f_0f_i^2\in
L^p(X,\widetilde{\mu})$ for every $i=1,\dots,n$, then $\left\{f_0,f_0f_1,\dots,f_0f_n,f_0\sum\limits_{i=1}^nf_i^2\right\}$ is a
Korovkin subset in $L^p(X,\widetilde{\mu})$.
\item [(2)]In particular,
$\left\{\textbf{1}, f_1,\dots, f_n,\sum\limits_{i=1}^nf_i^2\right\}$ is a
Korovkin subset in $C_*^w(X)$, where $$w:=\frac{1}{1+\sum\limits_{i=1}^nf_i^2}.$$
Therefore, if $\widetilde{\mu}$ is a finite
Borel measure on $X$ and $1\leq p<+\infty$, and if $f_i^2\in
L^p(X,\widetilde{\mu})$ for every $i=1,\dots,n$, then $\left\{\textbf{1},f_1,\dots,f_n,\sum\limits_{i=1}^nf_i^2\right\}$ is a
Korovkin subset in $L^p(X,\widetilde{\mu})$.

\end{itemize} \end{corollary}

The following particular case of Corollary 8.16 will be used next.

\begin{corollary}\sl  Let $X$ be an unbounded closed subset
of $\mathbb{R}^d, d\geq 1$.
The following statements hold:
\begin{itemize}
\item [(1)] If $f_0\in C(X)$
is a strictly positive function,
then
$\left\{f_0,f_0pr_1,\dots,f_0pr_d,f_0\norm{\cdot}^2\right\}$ is a
Korovkin subset in $C_*^w(X)$, where $$w:=\frac{1}{f_0\Big(1+\norm{\cdot}^2\Big)}.$$
Therefore, if $\widetilde{\mu}$ is a Borel measure on $X$ and $1\leq p<+\infty$, and
if $f_0\in L^p(X,\widetilde{\mu})$ as well as $\norm{\cdot}^2 f_0\in L^p(X,\widetilde{\mu})$,
then $\left\{f_0,f_0pr_1,\dots,f_0pr_d,f_0\norm{\cdot}^2\right\}$ is a
Korovkin subset in $L^p(X,\widetilde{\mu})$.
\item [(2)] In particular,
$\left\{\textbf{1},pr_1,\dots,pr_d,\norm{\cdot}^2\right\}$ is a
Korovkin subset in $C_*^w(X)$, where $$w:=\frac{1}{1+\norm{\cdot}^2}.$$
Therefore, if $\widetilde{\mu}$ is a finite
Borel measure on $X$ and $1\leq p<+\infty$, and if $\norm{\cdot}^2\in
L^p(X,\widetilde{\mu})$, then $\left\{\textbf{1},pr_1,\dots,pr_d,\norm{\cdot}^2\right\}$ is a
Korovkin subset in $L^p(X,\widetilde{\mu})$.
\end{itemize}
\end{corollary}
\begin{proof} It is sufficient to apply Corollary 8.16 with
$f_i=pr_i$, $1\leq i\leq d$.
\end{proof} Below, we list
some additional examples.
\begin{examples}\rm  \begin{itemize}\item []
\item [1)] Let $I$ be a noncompact real interval and set $r_1:=\inf
I\in\mathbb{R}\cup\{-\infty\}$ and $r_2:=\sup
I\in\mathbb{R}\cup\{+\infty\}$.  Consider a strictly positive
injective function $f_0\in C(I)$.

Assume that for some strictly positive function $w\in C(I)$ the following properties hold:

\begin{itemize}
    \item [(i)]$\lim\limits_{x\rightarrow
    r_i}w(x)f_0(x)=\lim\limits_{x\rightarrow r_i}w(x)f_0^2(x)=0$ for
    every $i=1,2$ such that $r_i\not\in I$; \item [(ii)] there exists
    $l>0$ such that $\lim\limits_{x\rightarrow r_i}w(x)f_0^3(x)=l$ for
    every $i=1,2$ such that $r_i\not\in I$.
\end{itemize}
Then $\{f_0,f_0^2,f_0^3\}$ is a Korovkin subset in $C_*^w(I)$
and hence in $L^p(I,\widetilde{\mu})$ for every Borel measure
$\widetilde{\mu}$ on $I$ and $p\in[1,+\infty[$ such that $1/w\in L^p(I,\widetilde{\mu})$.

The result is a direct consequence of
Corollary 8.14 with $n=1$ and $f_1=f_0$.
\item [2)] Let $I$ be an unbounded closed real interval. The following statements hold:
\begin{itemize}
\item [(1)] If $f_0\in C(I)$
is a strictly positive function, then
$\left\{f_0,f_0e_1,f_0e_2\right\}$ is a
Korovkin subset in $C_*^w(I)$, where $$w:=\frac{1}{f_0\Big(1+e_2\Big)}.$$
Therefore, if $\widetilde{\mu}$ is a Borel measure on $I$ and $1\leq p<+\infty$, and
if $f_0\in L^p(I,\widetilde{\mu})$ as well as $f_0e_2\in L^p(X,\widetilde{\mu})$,
then $\left\{f_0,f_0e_1,f_0e_2\right\}$ is a
Korovkin subset in $L^p(X,\widetilde{\mu})$.
\item [(2)] In particular,
$\left\{\textbf{1},e_1,e_2\right\}$ is a
Korovkin subset in $C_*^w(I)$, where $$w:=\frac{1}{1+e_2}.$$
Therefore, if $\widetilde{\mu}$ is a finite
Borel measure on $I$ and $1\leq p<+\infty$, and if $e_2\in
L^p(X,\widetilde{\mu})$, then $\left\{\textbf{1},e_1,e_2\right\}$ is a
Korovkin subset in $L^p(X,\widetilde{\mu})$.
\end{itemize}
\end{itemize}
\end{examples}

For a rather complete list of references on Korovkin-type theorems
in $L^p$-spaces, we refer to [8, Appendix D.2.3].

Next, we discuss a simple application of the above results (more
precisely, of Corollary 8.17) concerning the operators $G_n, n\geq
1$, defined by (\ref{4.44}) as
\begin{equation}\label{8.5}
G_n(f)(x):=\Big(\frac{n}{4}\Big)^{d/2}\int_{\mathbb{R}^d}f(t)\exp\Big(-\frac{n}{4}\norm{t-x}^2\Big)\dd
t\quad (x\in\mathbb{R}^d)
\end{equation}
for every Borel measurable function
$f:\mathbb{R}^d\rightarrow\mathbb{R}$ for which the integral
(\ref{8.5}) is absolutely convergent $(d\geq 1)$.

Let $\varphi:\mathbb{R}^d\rightarrow\mathbb{R}$ be a Borel
measurable strictly positive function which we assume to be
integrable with respect to the $d$-dimensional Lebesgue measure
$\lambda_d$. We denote by
$\widetilde{\mu}_{\varphi}:=\varphi\lambda_d$ the finite Borel
measure on $\mathbb{R}^d$ having density $\varphi$ with respect to
$\lambda_d$, i.e.,
\begin{equation}\label{8.6}
    \widetilde{\mu}_{\varphi}(B)=\int_B\varphi(x)\dd x
\end{equation}
for every Borel subset $B$ of $\mathbb{R}^d$. If $p\in[1,+\infty[$
and $f\in L^p(\mathbb{R}^d,\widetilde{\mu}_\varphi)$, we shall set

\begin{equation}\label{8.7}
    \norm{f}_{\varphi,p}:=\Big(\int_{\mathbb{R}^d}|f(x)|^p\varphi(x)\dd x\Big)^{1/p}.
\end{equation}
\begin{theorem}\sl  Assume that
\begin{equation}\label{8.8} C_\varphi:=\underset{n\geq
     1,\, t\in\mathbb{R}^d}{\sup}\Big\{\frac{1}{\varphi(t)}\Big(\frac{n}{4\pi}\Big)^{d/2}
    \int_{\mathbb{R}^d}\varphi(x)\exp\Big(-\frac{n}{4}
    \norm{t-x}^2\Big)\dd x\Big\}<+\infty.
\end{equation}
Then for every $p\in[1,+\infty[$ and $n\geq 1$ and for every $f\in
L^p(\mathbb{R}^d,\widetilde{\mu}_\varphi)$, the integral (\ref{8.5})
converges absolutely for a.e.\  $x\in\mathbb{R}^d$.

Moreover $G_n(f)\in L^p(\mathbb{R}^d,\widetilde{\mu}_\varphi)$ and

\begin{equation}\label{8.9}
    \norm{G_n(f)}_{\varphi,p}\leq\max\{1,C_\varphi\}\norm{f}_{\varphi,p}.
\end{equation} Finally, if $pr_i^2\in
    L^p(\mathbb{R}^d,\widetilde{\mu}_\varphi)$ for every $i=1,\dots,d$,
    then
\begin{equation}\label{8.10}
    \lim_{n\rightarrow\infty}G_n(f)=f\quad\text{with respect to}\quad
    \norm{\cdot}_{\varphi,p}.
\end{equation}

\end{theorem}
\begin{proof} Setting, for every $n\geq 1$,
$$K_n(x,t):=\frac{1}{\varphi(t)}\Big(\frac{n}{4\pi}\Big)^{d/2}\exp\Big(-\frac{n}{4}\norm{t-x}^2\Big),\quad
(t,x\in\mathbb{R}^d),$$ we get
$$G_n(f)(x)=\int_{\mathbb{\mathbb{R}^d}}K_n(t,x)f(t)\dd \widetilde{\mu}_\varphi(t)\quad
(x\in\mathbb{R}^d).$$ Moreover, if $x\in\mathbb{R}^n$ is fixed, then
$$\int_{\mathbb{R}^d}K_n(t,x)\dd \widetilde{\mu}_\varphi(t)=\Big(\frac{n}{4\pi}\Big)^{d/2}
\int_{\mathbb{R}^d}\exp\Big(-\frac{n}{4}\norm{t-x}^2\Big)\dd t=1$$
and, for $t\in\mathbb{R}^d$ fixed,
$$\int_{\mathbb{R}^d}K_n(t,x)\dd \widetilde{\mu}_\varphi(x)\leq
C_\varphi.$$ Therefore, the first part of the statement follows from
Fubini's and Tonelli's theorems and from H\"{o}lder's inequality
(see also [58, Theorem 6.18]). As regards the final part, note that
each $pr_i$ belongs to $L^p(\mathbb{R}^d,\widetilde{\mu}_\varphi)$
too because $|pr_i|\leq 1+pr_i^2,\quad (i=1,\dots,d)$.

From (\ref{8.9}), it follows that the sequence $(G_n)_{n\geq 1}$ is
equibounded from $L^p(\mathbb{R}^d,\widetilde{\mu}_\varphi)$ into
itself.  Moreover, formulae (\ref{4.48})--(\ref{4.50}) imply that
\[
G_{n}(h)\rightarrow h\text{ in }L^{p}(\mathbb{R}^{d},\tilde{\mu}_{\varphi })
\]
for every $h\in\left\{\textbf{1},pr_1,\dots,p_d,\sum\limits_{i=1}^dpr_i^2\right\}$
and hence the result follows from Corollary 8.17, (2).
\end{proof}

\begin{remark}\rm  It is not difficult to show that condition
(\ref{8.8}) is satisfied, for instance, by the functions
$\varphi_m(x):=(1+\norm{x}^m)^{-1}$ or
$\varphi_m(x)=\exp(-m\norm{x})\;(x\in\mathbb{R}^d)$ for every $m\geq
1$.  For more details, we refer to the papers [21, Section 3], [22]
and [23, Section 4] where an extension of Theorem 8.19 has been
established for a sequence of more general integral operators and
where the operators $G_{n}$, $n\geq 1$, have been studied also in
weighted continuous function spaces.
\end{remark}

\section{Korovkin-type theorems and Stone-Weierstrass theorems} In this
section, we shall deepen the connections between Korovkin-type
theorems and Stone-Weierstrass theorems, which have already been
pointed out with Theorems 3.8 and 4.7.

We start by giving a new proof of a generalization of Weierstrass'
theorem due to Stone ([111]) by means of Theorem 6.6. We first need
the following result.

\begin{lemma}\sl  Every closed
subalgebra $A$ of $C_b(X)$ is a lattice subspace (i.e., $|f|\in A$
for every $f\in A$).

\end{lemma}
\begin{proof} It is well-known that
$$t^{1/2}=\sum_{n=0}^\infty\binom{1/2}{n}(t-1)^n=\lim_{n\rightarrow\infty}p_n(t)$$
uniformly with respect to $t\in[0,2]$, where
$$p_n(t):=\sum_{k=0}^{n}\binom{1/2}{k}(t-1)^k\quad(n\geq 1,\, t\in[0,2]).$$
Since $\lim\limits_{n\rightarrow\infty}p_n(0)=0$, we also get
$$t^{1/2}=\lim_{n\rightarrow\infty}q_n(t)\quad\text{uniformly with respect
to}\quad t\in[0,2],$$ where   $q_n:=p_n-p_n(0)$, $n\geq 1$.
Therefore, for any $f\in A,\, f\neq 0,$
$$|f|=\norm{f}\left(\frac{f^2}{\norm{f}^2}\right)^{1/2}=\norm{f}
\lim_{n\rightarrow\infty}q_n\left(\frac{f^2}{\norm{f}^2}\right)\in\overline{A}=A.$$
\vskip-\baselineskip
\end{proof}

Before stating the Stone approximation theorem, we recall that a
subset $M$ of $C_0(X)$ is said to \dword{separate strongly the
points of $X$} if it separates the points of $X$ and if for every
$x\in X$ there exists $f\in M$ such that $f(x)\neq 0$.

\begin{theorem}\sl  Let $X$ be a
locally compact Hausdorff space with a countable base and let $A$ be
a closed subalgebra of $C_0(X)$ that separates strongly the points
of $X$. Then $A=C_0(X)$.
\end{theorem}
\begin{proof} We first prove that there
exists a strictly positive function $f_0$ in $A$.  For this, note
that $A$ is separable since $C_0(X)$ is. Denoting by $\{h_n\mid
n\geq 1\}$ a dense subset of $A$, then for every $x\in X$ there
exists $n\geq 1$ such that $h_n(x)\neq 0$.  Hence the function
$f_0:=\sum\limits_{n=1}^\infty\frac{|h_n|}{2^n\norm{h_n}}$ lies in
$A$, because of Lemma 9.1, and is strictly positive on $X$.

Given such a function $f_0\in A$, we then have that $\{f_0\}\cup
f_0A\cup f_0A^2\subset A$ and hence, by Theorem 6.6, $A$ is a
Korovkin subspace in $C_0(X)$.

Considering $f\in C_0(X)$ and $\varepsilon>0$, by Theorem 8.1 there
exist $h_0,\dots,h_n\in A$, $k_0,\dots,k_n\in A$ as well as $u,v\in
C_0(X), u\geq 0, v\geq 0$, such that $\norm{u}\leq\varepsilon,
\norm{v}\leq\varepsilon$, and, finally,
$$\Norm{\underset{0\leq j\leq n}{\inf}k_j-\underset{0\leq i\leq n}
{\sup}h_i}\leq\varepsilon\;\text{ and
}\;\underset{0\leq i\leq n} {\sup}h_i-u\leq f\leq\underset{0\leq
j\leq n}{\inf}k_j+v.$$ Then $f-\underset{0\leq j\leq n}{\inf}k_j\leq
v$ and
$$\underset{0\leq j\leq
n}{\inf}k_j-f\leq|\underset{0\leq j\leq n}{\inf}k_j-\underset{0\leq
i\leq n}{\sup}h_i|+u.$$ Therefore, $$|f-\underset{0\leq j\leq
n}{\inf}k_j|\leq|\underset{0\leq i\leq n}{\sup}h_i-\underset{0\leq
j\leq n}{\inf}k_j|+u+v$$ and  then $$\norm{f-\underset{0\leq j\leq
n}{\inf}k_j}\leq 3\varepsilon$$ which shows that
$f\in\overline{A}=A$ (since, by Lemma 9.1, $\underset{0\leq j\leq
n}{\inf}k_j\in A)$ and this completes the proof.
\end{proof}

As in the case of the interval $[0,1]$ (see Theorem 3.8), we
actually show that Theorems 6.6 and 9.2 are equivalent.  The next
result will be very useful for our purposes.

\begin{theorem}\sl  Let $X$ be a locally
compact Hausdorff space and consider a subset $M$ of $C_0(X)$ such
that the linear subspace generated by it contains a strictly
positive function $f_0\in C_0(X)$.  Given $x_0\in X$, denote by
$A(M,x_0)$ the subspace of all functions $f\in C_0(X)$ such that
$\mu(f)=f(x_0)$ for every $\mu\in\mathcal{M}_b^+(X)$ satisfying
$\mu(g)=g(x_0)$ for every $g\in\{f_0\}\cup f_0M\cup f_0M^2.$ Then
$A(M,x_0)$ is a closed subalgebra of $C_0(X)$ which contains $M$.
\end{theorem}
\begin{proof}
 $A(M,x_0)$ is a closed subspace of $C_0(X)$.  Set
\begin{eqnarray*}
M(x_0)&:=&\{x\in X\mid g(x)=g(x_0)\quad\text{for every}\quad g\in M\}\\
&=&\{x\in X\mid f_0(x)(g(x)-g(x_0))^2=0\quad\text{for every}\quad
g\in M\}
\end{eqnarray*} and fix $f\in A(M,x_0)$ and $\mu\in\mathcal{M}_b^+(X)$
such that $\mu(g)=g(x_0)$ for every $g\in\{f_0\}\cup f_0M\cup
f_0M^2$.  In particular, $\mu(f_0(g-g(x_0))^2)=0$ for each $g\in M$.
On the other hand, if $x\in M(x_0)$, then $\delta_x(g)=g(x_0)$ for
every $g\in\{f_0\}\cup f_0M\cup f_0M^2$, so that $f(x)=f(x_0)$ and
hence $f_0(x_0)f^2=f^2(x_0)f_0$ on $M(x_0)$.

Therefore, by Corollary 11.6 in the \Appendix, we get that
$\mu(f_0(x_0)f^2)=f^2(x_0)\mu(f_0)=f^2(x_0)f_0(x_0)$ and hence
$\mu(f^2)=f^2(x_0)$.  Therefore, $f^2\in A(M,x_0)$ and hence
$A(M,x_0)$ is a subalgebra of $C_0(X)$.

Finally note that, if $h\in M,$ then $f_0(x_0)h=h(x_0)f_0$ on
$M(x_0)$. Therefore, if we again fix $\mu\in\mathcal{M}_b^+(X)$ such
that $\mu(g)=g(x_0)$ for each $g\in\{f_0\}\cup f_0M\cup f_0M^2$, by
applying Corollary 11.6 we at once obtain that
$f_0(x_0)\mu(h)=h(x_0)\mu(f_0)=h(x_0)f_0(x_0)$, so $\mu(h)=h(x_0)$
and hence $h\in A(M,x_0)$.
\end{proof}

With the help of the preceding theorem, it is easy to show this next
result.

\begin{theorem}\sl  The Korovkin-type Theorem 6.6 and the
Stone-Weierstrass Theorem 9.2 are equivalent.
\end{theorem}
\begin{proof}
In light of the proof of Theorem 9.2, we have only to show that
Theorem 9.2 implies Theorem 6.6.  So, consider a subset $M$ of
$C_0(X)$ that separates the points of $X$ and such that the linear
subspace generated by it contains a strictly positive function
$f_0\in C_0(X)$.  \\ \indent Given $x_0\in X$, consider the subspace
$A(M, x_0)$ defined in the preceding theorem. By virtue of Theorems
9.2 and 9.3, we then infer that $A(M, x_0)=C_0(X)$.  In other words,
we have shown that for every $\mu\in\mathcal{M}_b^+(X)$ and for
every $x_0\in X$ satisfying $\mu(g)=g(x_0)$ for every
$g\in\{f_0\}\cup f_0M\cup f_0M^2$, we also have $\mu(f)=f(x_0)$ for
every $f\in C_0(X)$ and this means that $\{f_0\}\cup f_0M\cup
f_0M^2$ is a Korovkin subset in $C_0(X)$ because of Theorem 6.1.
\end{proof}

For a further deepening of the relationship between Korovkin-type
theorems and Stone-Weierstrass theorems, we refer to [8, Section
4.4] and [12-13].

\section{Korovkin-type theorems for positive projections} In this last
section, we shall discuss some Korovkin-type theorems for a class of
positive linear projections by showing a nontrivial application of
the general theorem proved in Section 5.

Consider a compact metric space $X$ and  a \dword{positive linear
projection} $T:C(X)\rightarrow C(X)$, i.e., $T$ is a positive linear
operator such that $T(T(f))=T(f)$ for every $f\in C(X)$. We shall
denote by $H_T$ the range of $T$, i.e.,
\begin{equation}\label{10.1} H_T:=T(C(X))=\{h\in C(X)\mid T(h)=h\}\,.
\end{equation} We shall also assume that $\textbf{1}\in H_T$ (hence
$T(\textbf{1})=\textbf{1}$) and that $H_T$ separates the points of
$X$. In the sequel, we shall present some examples of such
projections.

For every $x\in X$ denote by $\mu_x\in\mathcal{M}^+(X)$ the Radon
measure defined by
\begin{equation}\label{10.2} \mu_x(f):=T(f)(x)\quad (f\in C(X))
\end{equation} and by $\widetilde{\mu}_x$ the unique probability Borel
measure on $X$ that corresponds to $\mu_x$ via the Riesz
representation theorem, i.e.,
\begin{equation}\label{10.3}
T(f)(x)=\mu_x(f)=\int_Xf\dd \widetilde{\mu}_x\quad(f\in C(X)).
\end{equation}
Thus, $(\widetilde{\mu}_x)_{x\in X}$ is the canonical continuous
selection associated with $T$.

By the Cauchy-Schwarz inequality (\ref{2.16}), for every $h\in H_T$,
we get
$$|h|=|T(h)|\leq\sqrt{T(\textbf{1})T(h^2)}=\sqrt{T(h^2)}$$ so that
\begin{equation}\label{10.4} h^2\leq T(h^2)\,.
\end{equation} Set
\begin{equation}\label{10.5} Y_T:=\{x\in X\mid T(f)(x)=f(x)\text{ for
every } f\in C(X)\}.
\end{equation} $Y_T$ is closed and it is actually
equal to the so-called Choquet boundary of $H_T$ and hence it is not
empty (see [8, Proposition 3.3.1 and Corollary 2.6.5]).

\begin{proposition}\sl  If $M$ is a subset of $H_T$ that separates the points
of $X$, then
\begin{equation}\label{10.6} Y_T=\{x\in
X\;|\;T(h^2)(x)=h^2(x)\text{ for every } h\in M\}.
\end{equation}
Moreover, if $(h_n)_{n\geq 1}$ is a finite or countable family of
$H_T$ that separates the points of $X$ and if the series
$u:=\sum\limits_{n=1}^\infty h_n^2$ is uniformly convergent on $X$,
then $u\leq T(u)$ and
\begin{equation}\label{10.7} Y_T:=\{x\in X\mid Tu(x)=u(x)\}.
\end{equation}
\end{proposition}

\begin{proof} Consider $x\in X$ such that
$T(h^2)(x)=h^2(x)$, i.e.,  $\mu_x(h^2)=h^2(x)$ for each $h\in M$.
Since $M\subset H_T$, we also get $\mu_x(h)=h(x)\;(h\in M)$ and
$\{\textbf{1}\}\cup M\cup M^2$ is a Korovkin subset of $C(X)$ by
Theorem 7.1.  According to Theorem 6.1, we then obtain that
$f(x)=\mu_x(f)=T(f)(x)$ for every $f\in C(X)$.

As regards the second part of the statement, note that
$T(u)=\sum\limits_{n=1}^\infty T(h_n^2)$ and hence $u\leq T(u)$ by
(\ref{10.4}).  Moreover, if $x\in X$ and $T(u)(x)=u(x)$, then
$$\sum\limits_{n=1}^\infty \Bigl( T(h_n^2)(x)-h_n^2(x)\Bigr)=0.$$ Hence, on account of
(\ref{10.4}), $T(h_n^2)(x)=h_n^2(x)$ for every $n\geq 1$ and so
$x\in Y_T$ by (\ref{10.6}).
\end{proof}
\begin{remark}\rm  Note that there always exists
a countable family $(h_n)_{n\geq 1}$ of $H_T$ that separates the
points of $X$ and such that the series $\sum\limits_{n=1}^\infty
h_n^2$ is uniformly convergent on $X$. Actually, since $C(X)$ is
separable, $H_T$ is separable as well so that, considering a
countable dense family $(\varphi_n)_{n\geq 1}$ of $H_T$, it is
enough to put $h_n:=\frac{\varphi_n}{2^n(1+\norm{\varphi_n})}\quad
(n\geq 1).$
\end{remark}
\begin{theorem}\sl  If $M$ is a subset of $H_T$ that separates the points of
$X$, then $$H_T\cup M^2\quad\text{is a Korovkin subset for}\quad
T.$$ Moreover, if $u\in C(X)$ satisfies $u\leq T(u)$ and if
(\ref{10.7}) holds, then $$H_T\cup\{u\}\quad\text{is a Korovkin
subset for}\quad T.$$

In particular, the above statement holds for
$u=\sum\limits_{n=1}^\infty h_n^2$, where $(h_n)_{n\geq 1}$ is an
arbitrary sequence in $H_T$ that separates the points of $X$ and
such that the series $\sum\limits_{n=1}^\infty h_n^2$ is uniformly
convergent on $X$.
\end{theorem}
\begin{proof} In order to apply Theorem 5.5, consider
$\mu\in\mathcal{M}^+(X)$ and $x\in X$ such that
$\mu(h)=T(h)(x)=(\mu_x(h))$ for every $h\in H_T\cup M^2$.
Accordingly, if $h\in M$, $$\mu(T(h^2)-h^2)=T(h^2)(x)-\mu_x(h^2)=0$$
and $T(h^2)-h^2\geq 0$.  Therefore, by using (\ref{10.6}) and by
applying Theorem 11.5 and Corollary 11.6 of the \Appendix, we see
that, for every $f\in C(X)$, $\mu(T(f))=\mu(f)$ because $Tf=f$ on
$Y_T$.  Hence $\mu(f)=\mu(T(f))=\mu_x(Tf)=Tf(x)$ and this finishes
the proof.

A similar reasoning can be used to show the second part of the
statement by using (\ref{10.7}) instead of (\ref{10.6}).
\end{proof}

Theorem 10.3 is due to Altomare ([3], [5]; see also [8, Section
3.3]).  For an extension to so called adapted spaces, we refer to
[9]. Below, we show some examples.

\begin{example}\rm  Consider the $d$-dimensional simplex $K_d$ of
$\mathbb{R}^d$, $d\geq 1$, defined by (\ref{4.38}) and the positive
projection $T_d:C(K_p)\rightarrow C(K_p)$ defined by
\begin{equation}\label{10.8} T_d(f)(x):=\Big(1-\sum_{i=1}^d
x_i\Big)f(0)+\sum_{i=1}^dx_if(a_i)
\end{equation} $(f\in C(K_d)$,
$x=(x_i)_{1\leq i\leq d}\in K_d)$, where $a_i:=(\delta_{ij})_{1\leq
j\leq d}$ for every $i=1,\dots,d$.

In this case, $H_{T_d}$ is the subspace $A(K_d)$ (see Section 7,
formula (\ref{7.1})) which in turn is generated by
$M:=\{\textbf{1},pr_1,\dots,pr_d\}$, and hence, by Theorem 10.3,
\begin{equation}\label{10.9}
\left\{\textbf{1},pr_1,\dots,pr_d,\sum_{i=1}^dpr_i^2\right\}\text{
\textsl{is a Korovkin subset for }} T_d.
\end{equation} More generally, if $K$ is
an infinite-dimensional Bauer simplex (see, e.g., [8, Section 1.5, p. 59]), then there exists a unique
linear positive projection $T$ on $C(K)$ whose range is $A(K)$.  In
this case, for every strictly convex $u\in C(K)$, we get that
$A(K)\cup\{u\}$ is a Korovkin subset for $T$ (see [8, Corollary
3.3.4]).
\end{example}
\begin{example}\rm
Consider the hypercube $Q_d:=[0,1]^d$ of $\mathbb{R}^d$, $d\geq 1$,
and the positive projection $S_d:C(Q_d)\rightarrow C(Q_d)$ defined
by
\begin{equation}\label{10.10}
S_d(f)(x):=\sum_{h_1,\dots,h_d=0}^{1}x_1^{h_1}(1-x_1)^{1-h_1}\cdots
x_d^{h_d}(1-x_d)^{1-h_d}f(b_{h_1,\dots,h_d})
\end{equation} where
$b_{h_1,\dots,h_d}:=(\delta_{h_i1})_{1\leq i\leq d}\quad
(h_1,\dots,h_d\in\{0,1\}).$

In this case, $H_{S_d}$ is the subspace of $C(Q_d)$ generated by
$\{\textbf{1}\}\cup\{\prod_{i\in J}pr_i \mid
J\subset\{1,\dots,d\}\}.$ Therefore, by Theorem 10.3, we obtain that
\begin{equation}\label{10.11}
\{\textbf{1}\}\cup\left\{\prod_{i\in J}pr_i\mid
J\subset\{1,\dots,d\}\right\}\cup\left\{\sum_{i=1}^dpr_i^2\right\}\quad\text{is
a Korovkin subset for} \, S_d.
\end{equation} For an extension of the
above result, see [8, Corollary 3.3.9].
\end{example}
\begin{example}\rm
Consider a bounded open subset $\Omega$ of $\mathbb{R}^d$, $d\geq
2$, which we assume to be regular in the sense of potential theory
(see, e.g., [68, Section 8.3] or [8, pp.~125-128]) (for instance,
each convex open subset of $\mathbb{R}^d$ is regular).  Denote by
$H(\Omega)$ the subspace of all $u\in C(\overline{\Omega})$ that are
harmonic on $\Omega$.

By the regularity of $\Omega$, it follows that for every $f\in
C(\overline{\Omega})$ there exists a unique $u_f\in H(\Omega)$ such
that $u_f|_{\partial\Omega}=f|_{\partial\Omega}$, i.e.,  $u_f$ is
the unique solution of the Dirichlet problem
\begin{equation}\label{10.12}
\begin{cases} \triangle u:=\sum\limits_{i=1}^d\frac{\partial^2 u}{\partial
x_i^2}=0&\text{on}\;\Omega,\\
u|_{\partial\Omega}=f|_{\partial\Omega}
\end{cases}\quad (u\in C(\overline{\Omega})\cap C^2(\Omega)).
\end{equation}
Then the \dword{Poisson operator}
$T:C(\overline{\Omega})\longrightarrow C(\overline{\Omega})$ defined
by
\begin{equation}\label{10.13}
T(f):=u_f\quad (f\in C(\overline{\Omega}))
\end{equation} is a positive
projection whose range is $H(\Omega)$.

Therefore, from Theorem 10.3, we get that
\begin{equation}\label{10.14}
H(\Omega)\cup\left\{\sum_{i=1}^dpr_i^2\right\}\textsl{ is a Korovkin
subset for}\ T.
\end{equation} \end{example} For further examples of positive
linear projections, we refer to [24] and [103].

We now discuss an application of Theorem 10.3. Consider a metrizable
convex compact subset $K$ of a locally convex space. For every $f\in
C(K)$, $z\in K$ and $\alpha\in[0,1]$, we denote by $f_{z,\alpha}\in
C(K)$ the function defined by
\begin{equation}\label{10.15} f_{z,\alpha}(x):=f(\alpha x+(1-\alpha)z)\quad
(x\in K).
\end{equation} Consider a positive linear projection
$T:C(K)\longrightarrow C(K)$, $T$ different from the identity
operator, and assume that
\begin{equation}\label{10.16} A(K)\subset H_T:=T(C(K)),
\end{equation} i.e.,
\begin{equation}\label{10.17} T(u)=u\text{ for
every } u\in A(K),
\end{equation} and
\begin{equation}\label{10.18}
h_{z,\alpha}\in H_T\text{ for every } h\in H_T,\, z\in K,\,
0\leq\alpha\leq 1.
\end{equation} Here, once again, we use the symbol $A(K)$ to designate the subspace of all continuous affine functions on $K$ (see (\ref{7.1})).
Moreover, note that in all the examples 10.4-10.6,
assumptions (\ref{10.16})--(\ref{10.18}) are satisfied.

Moreover, from Theorem 10.3, it also follows that
\begin{equation}\label{10.19} H_T\cup A(K)^2\text{ is a Korovkin subset
for } T.
\end{equation} Consider the sequence of Bernstein-Schnabl
operators $(B_n)_{n\geq 1}$ associated with $T$ according to
(\ref{7.6}) and (\ref{7.7}), i.e., associated with the continuous
selection $(\widetilde{\mu}_x)_{x\in X}$ defined by (\ref{10.3}).

Note that, if $T$ is the projection $T_p$ defined by (\ref{10.8})
(resp., the projection $S_p$ defined by (\ref{10.10})), then the
corresponding Bernstein-Schnabl operators are the Bernstein
operators (\ref{4.39}) (resp., the Bernstein operators
(\ref{4.41})). By Theorem 7.5, we already know that
\begin{equation}\label{10.20}
\lim_{n\rightarrow\infty}B_n(f)=f\text{ uniformly on }
K\quad (f\in C(K)).
\end{equation} In this particular case, we can also study the
iterates of the operators $B_n$.

For every $n,m\geq 1$, we set
\begin{equation}\label{10.21}
B_n^m:=\begin{cases} B_n&{\rm if}\; m=1,\\ B_n\circ B_n^{m-1}&{\rm
if}\;m\geq 2.
\end{cases}
\end{equation} From (\ref{10.18}), it follows that, if $h\in
H_T$, then $B_n(h)=h$ and hence
\begin{equation}\label{10.22}
B_n^m(h)=h=T(h).
\end{equation} On the other hand, if $u\in A(K)$, then by
applying an induction argument it is not difficult to obtain from
(\ref{7.13}) that
\begin{equation}\label{10.23}
B_n(u^2)=\Big(1-\Big(\frac{n-1}{n}\Big)^m\Big)T(u^2)+\Big(\frac{n-1}{n}\Big)^mu^2.
\end{equation} Now, we can easily prove the next result.
\pagebreak
\begin{theorem}\sl
Let $f\in C(K)$.  Then
\begin{itemize}
\item [(1)]$\lim\limits_{n\rightarrow\infty}B_n^m(f)=f$ uniformly on $K$
for every $m\geq 1$.  \item
[(2)]$\lim\limits_{m\rightarrow\infty}B_n^m(f)=T(f)$ uniformly on
$K$ for every $n\geq 1$.  \item [(3)]If $(k(n))_{n\geq 1}$ is a
sequence of positive integers, then
$$\lim_{n\rightarrow\infty}B_n^{k(n)}(f)=\begin{cases}
  f&\text{uniformly on}\; K\;\text{if}\;\;\frac{k(n)}{n}\rightarrow 0,\\
  T(f)&\text{uniformly on}\; K\;\text{if}\;\;\frac{k(n)}{n}\rightarrow +\infty.
\end{cases}$$ \end{itemize} \end{theorem}
\begin{proof} It is sufficient
to apply Corollary 7.3 and statement (\ref{10.19}) taking
(\ref{10.22}) and (\ref{10.23}) into account as well as the
elementary formula
$$\Big(\frac{n-1}{n}\Big)^m=\exp\Big(m\log\Big(1-\frac{1}{n}\Big)\Big)\quad
(n,m\geq 2).$$
\vskip-\baselineskip
\end{proof}

It is worth remarking that under some additional assumptions on $T$,
the sequence $(B_n^{k(n)}(f))_{n\geq 1}$ $(f\in\mathcal{C}(K))$
converges uniformly also when $\frac{k(n)}{n}\longrightarrow
t\in]0,+\infty[$. More precisely, for every $t\geq 0$ there exists a
positive linear operator $T(t):C(K)\longrightarrow C(K)$ such that
for every sequence $(k(n))_{n\geq 1}$ of positive integers
satisfying $\frac{k(n)}{n}\longrightarrow t$, and for every $f\in
C(K)$,
\begin{equation}\label{10.24}
T(t)f=\lim_{n\rightarrow\infty}B_n^{k(n)}(f).
\end{equation} Moreover, the
family $(T(t))_{t\geq 0}$ is a strongly continuous semigroup of
operators whose generator $(A, D(A))$ is the closure of the operator
$(Z, D(Z))$ where
\begin{equation}\label{10.25} D(Z):=\left\{u\in
C(K)\;|\;\lim_{n\rightarrow\infty}n(B_n(u)-u)\quad\text{exists
in}\quad C(K)\right\}
\end{equation} and, for every $u\in D(Z)\subset D(A)$,
\begin{equation}\label{10.26}
A(u)=Z(u)=\lim_{n\rightarrow\infty}n(B_n(u)-u).
\end{equation} If $K$ is a
subset of $\mathbb{R}^d$, $d\geq 1$, with nonempty interior and if
$T$ maps the subspace of all polynomials of degree $m$ into itself
for every $m\geq 1$, then $C^2(K)\subset D(Z)\subset D(A)$ and, for
every $u\in C^2(K)$,
\begin{equation}\label{10.27}
Au(x)=Zu(x)=\frac{1}{2}\sum_{i,j=1}^d\alpha_{ij}(x)\frac{\partial^2u(x)}{\partial
x_i\partial x_i}
\end{equation} $(x=(x_i)_{1\leq i\leq d})$, where for every
$i,j=1,\dots,d$,
\begin{equation}\label{10.28}
\alpha_{ij}(x):=T(pr_ipr_j)(x)-x_ix_j.
\end{equation} The differential
operator (\ref{10.27}) is an elliptic second order differential
operator which degenerates on the subset $Y_T$ defined by
(\ref{10.5}) (which contains all the extreme points of $K$).

Moreover, for every $u_0\in D(A)$, the (abstract) Cauchy problem
\begin{equation}\label{10.29}
\begin{cases} \frac{\partial u}{\partial
t}(x,t)=A(u(\cdot,t))(x)& x\in K, t\geq 0,\\ u(x,0)=u_0(x)& x\in K,\\
u(\cdot,t)\in D(A)& t\geq 0,
\end{cases}\end{equation} has a unique
solution $u:K\times[0,+\infty[\rightarrow\mathbb{R}$ given by
\begin{equation}\label{10.30}
u(x,t)=T(t)(u_0)(x)=\lim_{n\rightarrow\infty}B_n^{k(n)}(u_0)(x)\quad
(x\in K, t\geq 0)
\end{equation} and the limit is uniform with respect to $x\in
K$, where $\frac{k(n)}{n}\rightarrow t$.

These results were discovered by Altomare ([5], see also [8,
Sections 6.2 and 6.3]) and opened the door to a series of researches
whose main aims are, first, the approximation of the solutions of
initial boundary differential problem associated with degenerate
evolution equations (like (\ref{10.29})) by means of iterates of
positive linear operators (like (\ref{10.30})), and, second, both a
numerical and a qualitative analysis of the solutions by means of
formula (\ref{10.30}).

These researches are documented in several papers.  Here, we content
ourselves to cite, other than Chapter 6 of [8], also the papers [6],
[17] and [18] and the references therein.

\section{Appendix: A short review of locally compact spaces and of some continuous function
spaces on them} For the convenience of the reader, in this Appendix,
we collect some basic definitions and results concerning locally
compact Hausdorff spaces, some continuous function spaces and Radon
measures on them.  For more details, we refer the reader to Chapter
IV of [30] or to Chapter 3 of [54].

We start by recalling that a topological space $X$ is said to be
\dword{compact} if every open cover of $X$ has a finite subcover. A
subset of a topological space is said to be compact if it is compact
in the relative topology. A topological space is said to be
\dword{locally compact} if each of its points possesses a compact
neighborhood.

Actually, if $X$ is locally compact and \dword{Hausdorff} (i.e.,
for every pair of distinct points $x_1,x_2\in X$ there exist
neighborhoods $U_1$ and $U_2$ of $x_1$ and $x_2$, respectively, such
that $U_1\cap U_2=\emptyset)$, then each point of $X$ has a
fundamental system of compact neighborhoods.

Every compact space is locally compact.  The spaces $\mathbb{R}^d$,
$d\geq 1$, are fundamental examples of (noncompact) locally compact
spaces.  Furthermore, if $X$ is locally compact, then every open
subset of $X$ and every closed subset of $X$, endowed with the
relative topology, is locally compact.  More generally, a subset of
a locally compact Hausdorff space, endowed with the relative
topology, is locally compact if and only if it is the intersection
of an open subset of $X$ with a closed subset of $X$ (see [54,
Corollary 3.3.10]).  Therefore, every real interval is locally
compact.

A topological space $X$ is said to be \dword{metrizable} if its
topology is induced by a metric on $X$.  In this case, we say that
$X$ is \dword{complete} if such a metric is complete. Note that every
compact metrizable space is complete and \dword{separable}, i.e., it
contains a dense countable subset.

A special role in the measure theory on topological spaces (and in
the Korovkin-type approximation theory) is played by locally compact
Hausdorff spaces with a \dword{countable base} (or basis), i.e.,
with a countable family of open subsets such that every open subset
is the union of some subfamily of it.  Such spaces are metrizable,
complete and separable. Actually, a metrizable space has a countable
base if and only if it is separable. The spaces $\mathbb{R}^d$,
$d\geq 1$, and each open or closed subset of them are locally
compact Hausdorff spaces with a countable base.

From now on, $X$ will stand for a fixed locally compact Hausdorff
space.  We denote by $K(X)$ the linear space of all real-valued
continuous functions $f:X\longrightarrow\mathbb{R}$ whose (closed)
support $$\supp(f):=\overline{\{x\in X\;|\;f(x)\neq 0\}}$$ is
compact. $K(X)$ is a lattice subspace of $C_b(X)$ and it coincides
with $C(X)$ if $X$ is compact.

The next result shows that there are sufficiently many functions in
$K(X)$. (For a proof, see [30, Corollary 27.3].)
\begin{theorem}\sl  (\dword{Urysohn's lemma})
For every compact subset $K$ of $X$ and for every open subset $U$
containing $K$, there exists $\varphi\in K(X)$ such that
$0\leq\varphi\leq 1$, $\varphi=1$ on $K$ and $\supp(\varphi)\subset
U$ (and hence $\varphi=0$ on $X\backslash U$).
\end{theorem} Another fundamental function space is the
space $C_0(X)$ which is defined as the closure of $K(X)$ in $C_b(X)$
with respect to the sup norm, in symbols
$$C_0(X):=\overline{K(X)}.$$
Thus, $C_0(X)$ is a closed linear subspace of $C_b(X)$ and hence,
endowed with the sup-norm, is a Banach space.

By means of Urysohn's lemma, it is not difficult to prove the
following characterization of functions lying in $C_0(X)$.
\begin{theorem}\sl  Assume
that $X$ is not compact.  For a function $f\in C(X)$, the following
statements are equivalent:
\begin{itemize}
\item [(i)] $f\in C_0(X)$; \item [(ii)]$\{x\in
X\;|\;|f(x)|\geq\varepsilon\}$ is compact for every $\varepsilon>0$;
\item [(iii)] for every $\varepsilon>0$ there exists a compact subset $K$
of $X$ such that $|f(x)|\leq\varepsilon$ for every $x\in X\backslash
K$.
\end{itemize}
\end{theorem}

Because of the preceding theorem, the functions lying in $C_0(X)$
are said to vanish at infinity. If $X$ is compact, then
$C_0(X)=C(X)$.  Moreover, $C_0(X)$ is a lattice subspace of $C_b(X)$
and, endowed with the sup-norm, is separable provided $X$ has a
countable base.

Another characterization of functions in $C_0(X)$ involves sequences
of points of $X$ that converge to the point at infinity of $X$. More
precisely, assuming that $X$ is noncompact, a sequence $(x_n)_{n\geq
1}$ in $X$ is said to \dword{converge to the point at infinity} of
$X$ if for every compact subset $K$ of $X$ there exists
$\nu\in\mathbb{N}$ such that $x_n\in X\backslash K$ for every $n\geq
\nu$. For any such sequence and for every $f\in C_0(X)$, we then
have
\begin{equation}\label{11.1}
\lim_{n\rightarrow\infty}f(x_n)=0.
\end{equation}

Conversely, a function $f\in C(X)$ satisfying (\ref{11.1}) for every
sequence $(x_n)_{n\geq 1}$ converging to the point at infinity of
$X$ necessarily lies in $C_0(X)$ provided that $X$ is countable at
infinity, i.e.,  it is the union of a sequence of compact subsets of
$X$. Note also that $X$ is countable at infinity if and only if
there exists $f_0\in C_0(X)$ such that $f_{0}(x)>0$ for every $x\in
X$. Moreover, if $X$ has a countable basis, then $X$ is countable at
infinity.

A useful tool which plays an important role in Korovkin-type
approximation theory is given by \eword{Radon measures}.  Actually,
we shall only need to handle \dword{positive bounded Radon measures}
which are, by definition, positive linear functionals on $C_0(X)$.
The set of all of them will be denoted by $\mathcal{M}_b^+(X)$.

Every $\mu\in \mathcal{M}_b^+(X)$, that is, every positive linear
functional $\mu:C_0(X)\longrightarrow\mathbb{R}$, is continuous
(with respect to the sup-norm), and its norm
\begin{equation}\label{11.2}
\norm{\mu}:=\sup\{|\mu (f)|\mid f\in C_0(X), |f|\leq 1\}
\end{equation} is
also called the \dword{total mass} of $\mu$.

A simple example of bounded positive Radon measure is furnished by
the \dword{Dirac measure} at a point $a\in X$, which is defined by
\begin{equation}\label{11.3} \delta_a(f):=f(a)\quad (f\in C_0(X)).
\end{equation}
A positive linear combination of Dirac measures is called a
(positive) \dword{discrete measure}.

In other words, a Radon measure $\mu\in\mathcal{M}_b^+(X)$ is
discrete if there exist finitely many points $a_1,\dots,a_n\in X,
n\geq 1$, and finitely many positive real numbers $\lambda_1,\dots
,\lambda_n$ such that
\begin{equation}\label{11.4}
\mu=\sum_{i=1}^n\lambda_i\delta_{a_i},
\end{equation} i.e.,
\begin{equation}\label{11.5}
\mu(f)=\sum_{i=1}^n\lambda_if(a_i)\quad\text{for every}\quad f\in
C_0(X).
\end{equation}
In this case, $\norm{\mu}=\sum\limits_{i=1}^n\lambda_i$ and $\mu$ is
also said to be \dword{supported on} $\{a_1,\dots,a_n\}$.

There is a strong relationship between positive bounded Radon
measures and (positive) finite Borel measures on $X$. In order to
briefly describe it, we recall that the \dword{Borel
$\sigma$-algebra} in $X$ is, by definition, the $\sigma$-algebra
generated by the system of all open subsets of $X$.  It will be
denoted by $\mathcal{B}(X)$ and its elements are called Borel
subsets of $X$.  Open subsets, closed subsets and compact subsets of
$X$ are Borel subsets.

A \dword{Borel measure} $\widetilde{\mu}$ on $X$ is, by definition,
a measure $\widetilde{\mu}:\mathcal{B}(X)\rightarrow[0,+\infty]$
such that
\begin{equation}\label{11.6} \widetilde{\mu}(K)<+\infty\quad\text{for every
compact subset}\quad K\;\text{of}\;X.
\end{equation}
Every finite measure $\widetilde{\mu}$ on $\mathcal{B}(X)$, i.e.,
$\widetilde{\mu}(X)<+\infty$, is a Borel measure.

A measure $\widetilde{\mu}:\mathcal{B}(X)\longrightarrow[0,+\infty]$
is said to be \dword{inner regular} if
\begin{equation}\label{11.7}
\widetilde{\mu}(B)=\sup\{\widetilde{\mu}(K)\mid K\subset B,
K\;\text{compact}\}\text{ for every }B\in\mathcal{B}(X)
\end{equation}
and \dword{outer regular} if
\begin{equation}\label{11.8}
\widetilde{\mu}(B)=\inf\{\widetilde{\mu}(U)\mid B\subset U,
U\;\text{open}\}\text{ for every }B\in\mathcal{B}(X).
\end{equation} A measure $\widetilde{\mu}$ is said to be \dword{regular}
if it is both inner regular and outer regular. The Lebesgue-Borel
measure on $\mathbb{R}^d$, $d\geq 1$, is regular. Actually, if $X$
has a countable base, then every Borel measures on $X$ is regular
([30, Theorem 29.12]).

If $\widetilde{\mu}$ is a finite measure on $\mathcal{B}(X)$, then
every $f\in C_b(X)$ is $\widetilde{\mu}$-integrable.  Therefore, we
can consider the positive bounded Radon measure
$I_{\widetilde{\mu}}$ on $X$ defined by
\begin{equation}\label{11.9}
I_{\widetilde{\mu}}(f):=\int_Xf\dd \widetilde{\mu}\quad (f\in
C_0(X)).
\end{equation}
Then $\norm{I_{\widetilde{\mu}}}\leq\widetilde{\mu}(X)$, and, if
$\widetilde{\mu}$ is inner regular,
$\norm{I_{\widetilde{\mu}}}=\widetilde{\mu}(X)$.

As a matter of fact, formula (\ref{11.9}) describes all the positive
bounded Radon measures on $X$ as the following fundamental result
shows (see [30, Section 29]).
\begin{theorem}\sl  (\dword{Riesz representation theorem}) If $\mu\in \mathcal{M}_b^+(X)$,
then there exists a unique finite and regular Borel measure
$\widetilde{\mu}$ on $X$ such that $$\mu(f)=\int_Xf\dd
\widetilde{\mu}\quad\text{for every}\quad f\in C_0(X).$$ Moreover,
$\norm{\mu}=\widetilde{\mu}(X).$
\end{theorem}

Another noteworthy and useful result concerning regular Borel
measures is shown below. Consider a measure $\widetilde{\mu}$ on
$\mathcal{B}(X)$ and $p\in[1,+\infty[$.  As usual, we shall denote
by $\mathcal{L}^p(X,\widetilde{\mu})$ the linear subspace of all
$\mathcal{B}(X)$-measurable functions $f:X\longrightarrow\mathbb{R}$
such that $|f|^p$ is $\widetilde{\mu}$-integrable.

If $f\in\mathcal{L}^p(X,\widetilde{\mu})$, it is customary to set
\begin{equation}\label{11.10} N_p(f):=\left(\int_X|f|^p\dd \widetilde{\mu}\right)^{1/p}.
\end{equation} The functional
$N_p:\mathcal{L}^p(X,\widetilde{\mu})\rightarrow\mathbb{R}$ is a
seminorm and the convergence with respect to it is the usual
\dword{convergence in $p^{th}$-mean}.  Setting
\begin{eqnarray*}
\mathcal{N}&:=&\{f\in\mathcal{L}^p(X,\widetilde{\mu})\mid N_p(f)=0\}\\
&=&\{f\in F(X)\mid f\;\text{is}\;\mathcal{B}(X)\;\text{measurable
and}\;f=0\;\;\widetilde{\mu}\,\,\text{a.e.}\},
\end{eqnarray*} the quotient linear
space
\begin{equation}\label{11.11}
L^p(X,\widetilde{\mu}):=\mathcal{L}^p(X,\widetilde{\mu})/\mathcal{N}
\end{equation} endowed with the norm
\begin{equation}\label{11.12}
\norm{\widetilde{f}}_p:=N_p(f)\quad(\widetilde{f}\in
L^p(X,\widetilde{\mu})),
\end{equation} is a Banach space.  (Here,
$\widetilde{f}:=\{g\in\mathcal{L}^p(X,\widetilde{\mu})\mid
f=g\;\;\widetilde{\mu}\,\, \text{a.e.}\}.$) Note that, if
$\widetilde{\mu}$ is a Borel measure, then
$K(X)\subset\mathcal{L}^p(X,\widetilde{\mu})$ for every
$p\in[1,+\infty[$.

If $\widetilde{\mu}$  is also regular, we can say much more (for a
proof of the next result, see [30, Theorem 29.14]).
\begin{theorem}\sl  If
$\widetilde{\mu}$ is a regular Borel measure on $X$, then, for every
$p\in[1,+\infty[$, the space $K(X)$ is dense in
$\mathcal{L}^p(X,\widetilde{\mu})$ with the respect to convergence
in the $p^{th}$-mean (and hence in $L^p(X,\widetilde{\mu})$ with
respect to $\norm{\cdot}_p$).
\end{theorem}

Next, we discuss a characterization of discrete Radon measures.
Usually this characterization is proved by using the notion of
support of Radon measures (see, e.g., [36, Chapter III, Section 2]
and [39, Vol.  I, Section 11]). Below, we present a simple and
direct proof.

We start with the following result which is important in its own
right.
\begin{theorem}\sl  Let $\mu\in\mathcal{M}_b^+(X)$ and consider a closed subset
$Y$ of $X$ such that
\begin{equation}\label{11.13}
\mu(\varphi)=0\text{ for every } \varphi\in
K(X), \, \supp(\varphi)\subset X\backslash Y.
\end{equation} Then
$\mu(f)=\mu(g)$ for every $f,g\in C_0(X)$ such that $f=g$ on $Y$.
\end{theorem}
\begin{proof} It suffices to show that, if $f\in C_0(X)$ and
$f=0$ on $Y$, then $\mu(f)=0$.  Consider such a function $f\in
C_0(X)$ and, given $\varepsilon>0$, set $U:=\{x\in X\mid
|f(x)|<\varepsilon\}$.  Hence, by Theorem 11.2, $X\backslash U$ is
compact and $X\backslash U\subset X\backslash Y$.

By Urysohn's lemma (Theorem 11.1), there exists $\varphi\in K(X)$,
$0\leq\varphi\leq 1$, such that $\varphi=1$ on $X\backslash U$ and
$\supp(\varphi)\subset X\backslash Y$.  In particular,
$\supp(f\varphi)\subset \supp(\varphi)\subset X\backslash Y$ and
hence $\mu(f\varphi)=0$.

Therefore, $$|\mu
(f)|=|\mu(f)-\mu(f\varphi)|\leq\norm{\mu}\norm{f(1-\varphi)}\leq\norm{\mu}\varepsilon,$$
because $\norm{f(1-\varphi)}\leq\varepsilon$ as we now confirm.  For
every $x\in X$, we have, indeed, $|f(x)(1-\varphi(x))|=0$ if
$x\not\in U$ and, if $x\in U,
|f(x)(1-\varphi(x))|\leq|f(x)|\leq\varepsilon$. Since
$\varepsilon>0$ was arbitrarily chosen, we conclude that $\mu(f)=0$.
\end{proof}

We point out that there always exists a closed subset $Y$ of $X$
satisfying (\ref{11.13}).  The smallest of them is called the
support of the measure $\mu$ (see the references before Theorem
11.5). An important example of a subset $Y$ satisfying (\ref{11.13})
is given below.

\begin{corollary}\sl  Let $\mu\in\mathcal{M}_b^+(X)$ and consider an arbitrary
family $(f_i)_{i\in I}$ of positive functions in $C_0(X)$ such that
$\mu(f_i)=0$ for every $i\in I$. Then the subset $$Y:=\{x\in X\mid
f_i(x)=0\text{ for every } i\in I\}$$ satisfies
(\ref{11.13}).

Therefore, if $f,g\in C_0(X)$ and if
$$f(x)=g(x)\quad\text{for every}\quad x\in Y,$$ then
$\mu(f)=\mu(g)$.
\end{corollary}

\begin{proof}
Consider $\varphi\in K(X)$ such that
$$\supp(\varphi)\subset
X\backslash Y=\{x\in X \mid \mbox{ there exists $i\in I$ with }
f_i(x)>0\}.$$ By using a compactness argument, we then find a finite
subset $J$ of $I$ such that
$$\supp(\varphi)\subset\bigcup_{i\in J}\{x\in X\mid f_i(x)>0\}.$$ If we set
$\alpha:=\min\left\{\sum\limits_{i\in J}f_i(x)\mid x\in \supp(\varphi)\right\}>0$, we
then obtain $$|\varphi|\leq\frac{\norm{\varphi}}{\alpha}\sum_{i\in
J}f_i$$ and hence $\mu(f)=0$.
\end{proof} By means of the previous result, it is easy
to reach the announced characterization of discrete Radon measures.
\begin{theorem}\sl  Given $\mu\in\mathcal{M}_b^+(X)$ and different points
$a_1,\dots,a_n\in X, n\geq 1$, the following statements are
equivalent:
\begin{itemize}
\item [(i)] There exist $\lambda_1,\dots,\lambda_n\in[0,+\infty[$ such
that $\mu=\sum\limits_{i=1}^n\lambda_i\delta_{a_i}$ (see (\ref{11.5})); \item
[(ii)] if $\varphi\in K(X)$ and
$\supp(\varphi)\cap\{a_1,\dots,a_n\}=\emptyset$, then
$\mu(\varphi)=0$;
\item [(iii)] for every $x\in X\backslash\{a_1,\dots,a_n\}$ there exists
$f\in C_0(X)$, $f\geq 0$, such that $f(x)>0$, $f(a_i)=0$ for each
$i=1,\dots,n$, and $\mu(f)=0$.
\end{itemize}
\end{theorem}
\begin{proof}
(i)$\Rightarrow$(ii).  It is obvious.

(ii)$\Rightarrow$(iii).  If $x\in X\backslash\{a_1,\dots,a_n\}$, by
Urysohn's lemma, we can choose $\varphi\in K(X)$, $0\leq\varphi\leq
1$, such that $\varphi(x)=1$ and $\supp(\varphi)\subset
X\backslash\{a_1,\dots,a_n\}$ so that $\mu(f)=0$.

(iii)$\Rightarrow$(ii).  Consider $\varphi\in K(X)$ such that
$\supp(\varphi)\cap\{a_1,\dots,a_n\}=\emptyset$.  By hypothesis, for
every $x\in \supp(\varphi)$, there exists $f_x\in C_0(X)$, $f_x\geq
0$, such that $f_x(x)>0$, $f_x(a_i)=0$ for every $i=1,\dots,n$, and
$\mu(f_x)=0$.

Since $\supp(\varphi)$ is compact, there exist $x_1,\dots,x_p\in
\supp(\varphi)$ such that
$$\supp(\varphi)\subset\bigcup_{k=1}^p\{x\in X\;|\;f_{x_k}(x)>0\}.$$
Therefore, the function $f:=\sum\limits_{k=1}^pf_{x_k}\in C_0(X)$ is
positive, it does not vanish  at any point of $\supp(\varphi)$, and
$\mu(f)=0$.

If we set $m:=\min\{f(x)\mid x\in \supp(\varphi)\}>0$, it is
immediate to verify that $m|\varphi|\leq\norm{\varphi}f$ and hence
$\mu(\varphi)=0$.

(ii)$\Rightarrow$(i).  For every $j=1,\dots,n$, consider
$\varphi_j\in K(X)$ such that $0\leq\varphi_j\leq 1$,
$\varphi_j(a_j)=1$ and $\varphi_j(a_i)=0$ for each $i=1,\dots,n$,
$i\neq j$.

If $f\in C_0(X)$, then
$$f=\sum_{i=1}^nf(a_i)\varphi_i\quad\text{on}\quad\{a_1,\dots,a_n\}.$$ On
the other hand, the subset $\{a_1,\dots,a_n\}$ satisfies
(\ref{11.13}) and hence, by Theorem 11.5,
$$\mu(f)=\sum_{i=1}^nf(a_i)\mu(\varphi_i)=\sum_{i=1}^n\lambda_if(a_i)$$
where $\lambda_i:=\mu(\varphi_i)\geq 0\quad(i=1,\dots,n)$ and this
completes the proof.
\end{proof}

We end the Appendix by discussing some aspects of vague convergence
for Radon measures.  For more details, we refer, e.g., to [30, \S
30] or to [39, Vol.  I, \S 12].

A sequence $(\mu_n)_{n\geq 1}$ in $\mathcal{M}_b^+(X)$ is said to
\dword{converge vaguely} to $\mu\in\mathcal{M}_b^+(X)$ if
\begin{equation}\label{11.14}
\lim_{n\rightarrow\infty}\mu_n(f)=\mu(f)\quad\text{for every}\quad
f\in C_0(X).
\end{equation}
Thus, (\ref{11.14}) simply means that $\mu_n\rightarrow\mu$ weakly
in the dual space of $C_0(X)$.

If, in addition, $X$ has a countable base, then $C_0(X)$ is
separable and hence, by Banach's theorem, the unit ball of the dual
of $C_0(X)$ is weakly sequentially compact.  Therefore

\begin{theorem}\sl  If $X$ has a countable base, then every sequence
in $\mathcal{M}_b^+(X)$ that is bounded with respect to the norm
(\ref{11.2}), has a subsequence that converges vaguely to some
$\mu\in\mathcal{M}_b^+(X)$.
\end{theorem}

\noindent Francesco Altomare\\
Dipartimento di Matematica\\
Universit\'{a} degli Studi di Bari ``A. Moro"\\
Campus Universitario\\
Via E. Orabona, 4\\
70125 Bari - Italia\\
E-mail: altomare@dm.uniba.it\\
URL: http://www.dm.uniba.it/Members/altomare

\endddoc